\documentclass{gtart_h}

\def\ifplaintex{\expandafter\ifx\csname documentclass\endcsname\relax}

\def\ifplaintex{\expandafter\ifx\csname documentclass\endcsname\relax}

\ifplaintex 
\else
\fi

\expandafter\ifx\csname epsfbox\endcsname\relax\input epsf\fi

\def\gt{{\mathsurround=0pt\it $\cal G\mskip-2mu$eometry \&\ 
$\cal T\!\!$opology}}        

\def\gtp{{\mathsurround=0pt\it $\cal G\mskip-2mu$eometry \&\ 
$\cal T\!\!$opology $\cal P\!$ublications}}  

\def\lognumber#1{\def\thelognumber{#1}}
\def\volumenumber#1{\def\thevolumenumber{#1}}
\def\papernumber#1{\def\thepapernumber{#1}}
\def\volumeyear#1{\def\thevolumeyear{#1}}

\def\pagenumbers#1#2{\def\startpage{#1}\def\finishpage{#2}}
\def\published#1{\def\publishdate{#1}}
\def\proposed#1{\def\theproposer{#1}}
\def\seconded#1{\def\theseconders{#1}}
\def\received#1{\def\receiveddate{#1}}

\def\accepted#1{\def\accepteddate{#1}}

\def\shorttitle#1{\def\theshorttitle{#1}}

\let\\\par\let\thelognumber\relax
\let\thevolumenumber\relax\let\thepapernumber\relax
\let\thevolumeyear\relax\let\thesamplenumber\relax\let\startpage\relax
\let\finishpage\relax\let\publishdate\relax\let\receiveddate\relax
\let\reviseddate\relax\let\accepteddate\relax\let\theasciititle\relax
\let\theasciiauthors\relax
\let\theasciiabstract\relax
\let\theasciiemail\relax\let\theshortauthors\relax\let\theshorttitle\relax

\long\def\maketitlep{   

\count0=\startpage

\gt\hfill      
\hbox to 77pt{\vbox to 0pt{\vglue -15pt\epsfbox{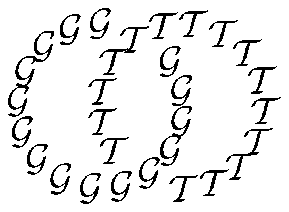}\vss}\hss}
\break
{\small\ifx\thesamplenumber\relax 
Volume \else Sample
\fi\thevolumenumber\ (\thevolumeyear)
\startpage--\finishpage\nl
Published: \publishdate}
\vglue 0.5truein plus 0.4fil minus 0.1truein

{\parskip=0pt\leftskip 0pt plus 1fil\def\\{\par\smallskip}{\ifplaintex\large
\else\Large\fi\bf\thetitle}\par\medskip}   

\vglue 0pt plus 0.1fil 

{\parskip=0pt\leftskip 0pt plus 1fil\def\\{\par}{\sc\theauthors}
\par\medskip}

\vglue 0pt plus 0.1fil 

{\small\parskip=0pt\let\newline\\
{\leftskip 0pt plus 1fil\def\\{\par}{\sl\theaddress}\par}
\expandafter\ifx\theemail\relax    
\relax\else\vglue 5pt plus 0.02fil minus 2pt\def\\{\stdspace{\rm 
and}\stdspace} 
\cl{Email:\stdspace\tt\theemail}\fi
\ifx\theurl\relax                  
\relax\else\vglue 5pt plus 0.02fil minus 2pt\def\\{\stdspace{\rm 
and}\stdspace}
\cl{URL:\stdspace\tt\theurl}\fi\par}

\vglue 7pt plus 0.3fil minus 3pt

{\bf Abstract}
\vglue 5pt plus 0.1fil minus 2pt

\theabstract

\vglue 7pt plus 0.3fil minus 3pt

{\bf AMS Classification numbers}\quad Primary:\quad \theprimaryclass

Secondary:\quad \thesecondaryclass

\vglue 5pt plus 0.3fil minus 2pt

{\bf Keywords:}\quad \thekeywords

\vglue 10pt plus 0.5fil minus 5pt

{\small  Proposed: \theproposer\hfill Received: \receiveddate\nl
Seconded: \theseconders\hfill 
\ifx\reviseddate\relax                         
Accepted: \accepteddate                        
\else
Revised: \reviseddate                          
\fi}
\eject
}       

\font\phead=cmsl9 scaled 950
\font=cmsl9 scaled 1050
\font\pnum=cmbx10 scaled 913
\font\lnum=cmbx10 
\font\pfoot=cmsl9 scaled 950
\font=cmsl9 scaled 1050
\ifplaintex
\headline{\vbox to 0pt{\vskip -4.5mm\line{\small\phead\ifnum
\count0=\startpage ISSN 1364-0380 (on line)
1465-3060 (printed) \hfill {\pnum\folio}\else\ifodd\count0\def\\{ }%
\ifx\theshorttitle\relax\thetitle\else\theshorttitle\fi\hfill{\pnum\folio}
\else\def\\{ and }{\pnum\folio}\hfill\ifx\theshortauthors\relax\theauthors
\else\theshortauthors\fi\fi\fi}\vss}}
\footline{\vbox to 0pt{\vglue 0mm\line{\small\pfoot\ifnum\count0=\startpage
\copyright\ \gtp\hfill\else
\gt, Volume \thevolumenumber\ (\thevolumeyear)\hfill\fi}\vss
}}
\else
\makeatletter
\def\@oddhead{{\small\lhead\ifnum\count0=\startpage ISSN 1364-0380 (on line)
1465-3060 (printed) \hfill {\lnum\number\count0}\else\ifodd\count0
\def\\{ }\ifx\theshorttitle\relax \thetitle \else\theshorttitle\fi\hfill
{\lnum\number\count0}\else\def\\{ and }{\lnum\number\count0}
\hfill\ifx\theshortauthors\relax 
\theauthors\else\theshortauthors\fi\fi\fi}}\def\@evenhead{\@oddhead}
\def\@oddfoot{\small\lfoot\ifnum\count0=\startpage\copyright\ \gtp\hfill\else
\gt, Volume \thevolumenumber\ (\thevolumeyear)\hfill\fi}
\def\@evenfoot{\@oddfoot}
\makeatother
\fi

\newwrite\gtoutfile
\long\gdef\makeheadfile{  
{\def\\{, }\def\s{ }
\immediate\openout\gtoutfile head.xxx
\immediate\write\gtoutfile{Proxy-for: \ifx\theasciiauthors\relax
\theauthors\else\theasciiauthors\fi\s<\ifx\theasciiemail\relax\theemail\else\theasciiemail\fi>}
\immediate\write\gtoutfile{\noexpand\\}
\immediate\write\gtoutfile{Authors: \ifx\theasciiauthors\relax
\theauthors\else\theasciiauthors\fi}
{\def\\{ }\immediate\write\gtoutfile{Title: \ifx\theasciititle\relax
\thetitle\else\theasciititle\fi}}
\immediate\write\gtoutfile{Subj-class: GT or SG or MG etc}
\immediate\write\gtoutfile{MSC-class: \theprimaryclass\ifx\thesecondaryclass\relax\else, \thesecondaryclass\fi}
\immediate\write\gtoutfile{Journal-ref: Geom. Topol. \thevolumenumber
(\thevolumeyear) \startpage-\finishpage}
\immediate\write\gtoutfile{Comments: Published by Geometry and Topology at}
\immediate\write\gtoutfile{\s\s http://www.maths.warwick.ac.uk/gt/GTVol\thevolumenumber/paper\thepapernumber.abs.html}
\immediate\write\gtoutfile{\noexpand\\}
\immediate\write\gtoutfile{}
\ifx\theasciiabstract\relax
\immediate\write\gtoutfile{\theabstract}\else
\immediate\write\gtoutfile{\theasciiabstract}\fi
\immediate\write\gtoutfile{}
\immediate\write\gtoutfile{\noexpand\\}
\immediate\write\gtoutfile{}
\immediate\closeout\gtoutfile}}  

\def\maketitlepage{\maketitlep\makeheadfile}
\let\maketitle\maketitlepage

\lognumber{582}

\received{18 March 2005}
\volumenumber{9}\papernumber{20}\volumeyear{2005}
\pagenumbers{833}{934} 
\published{23 May 2005\nl Revised: 13 December 2005 (see footnote 3
on page 834)}
\accepted{6 May 2005}
\proposed{Thomas Goodwillie}
\seconded{Ralph Cohen, Gunnar Carlsson}

\usepackage{amsmath,amssymb}
\usepackage[cmtip,arrow]{xy}
\usepackage{pb-diagram,pb-xy}
\usepackage[dvips]{color}
\usepackage{rotating}
\usepackage{graphicx}

\DeclareMathOperator*{\Barwedge}{\overline{\bigwedge}}

\numberwithin{equation}{section}

\theoremstyle{plain}

\newtheorem{thm}[equation]{Theorem}
\newtheorem{lemma}[equation]{Lemma}
\newtheorem{lem}[equation]{Lemma}
\newtheorem{cor}[equation]{Corollary}
\newtheorem{corollary}[equation]{Corollary}
\newtheorem{proposition}[equation]{Proposition}
\newtheorem{prop}[equation]{Proposition}

\theoremstyle{definition}
\newtheorem{definition}[equation]{Definition}
\newtheorem{remark}[equation]{Remark}

\newtheorem{example}[equation]{Example}
\newtheorem{examples}[equation]{Examples}

\newcommand{\isom}{\cong}
\newcommand{\basept}{\ast}
\newcommand{\union}{\cup}
\newcommand{\cat}[1]{\mathcal{#1}}
\newcommand{\ord}[1]{#1\textsuperscript{th}}
\newcommand{\smsh}{\wedge}

\newcommand{\wdge}{\vee}
\newcommand{\bigwdge}{\bigvee}
\newcommand{\varcirc}{\mathbin{\widehat{\circ}}}
\newcommand{\dual}{\mathbb{D}}
\newcommand{\spaces}[1][]{\mathcal{U}_{#1}}
\newcommand{\based}[1][]{\mathcal{T}_{#1}}
\newcommand{\spectra}[1][]{\mathcal{S}p}

\DeclareMathOperator*{\colim}{colim}

\DeclareMathOperator{\Map}{Map}
\DeclareMathOperator{\Hom}{Hom}
\DeclareMathOperator*{\hofib}{hofib}
\DeclareMathOperator*{\hocofib}{hocofib}

\newenvironment{fake}{\relax}{\relax}

\begin{document}

\title{Bar constructions for topological operads and\\the Goodwillie
derivatives of the identity}
\shorttitle{Bar constructions for topological operads}
\author{Michael Ching}
\address{Department of Mathematics, Room 2-089 \\
    Massachusetts Institute of Technology \\
    Cambridge, MA 02139, USA}
\email{mcching@math.mit.edu}

\begin{abstract}
We describe a cooperad structure on the simplicial bar construction
on a reduced operad of based spaces or spectra and, dually, an operad
structure on the cobar construction on a cooperad. We also show that if
the homology of the original operad (respectively, cooperad) is Koszul,
then the homology of the bar (respectively, cobar) construction is the
Koszul dual. We use our results to construct an operad structure on the
partition poset models for the Goodwillie derivatives of the identity
functor on based spaces and show that this induces the `Lie' operad
structure on the homology groups of these derivatives. We also extend
the bar construction to modules over operads (and, dually, to comodules
over cooperads) and show that a based space naturally gives rise to a
left module over the operad formed by the derivatives of the identity.
\end{abstract}

\begin{fake}\end{fake}

\maketitlepage

\section*{Introduction}
The motivation for this paper was an effort to construct an operad
structure on the derivatives (in the sense of Tom Goodwillie's homotopy
calculus \cite{goodwillie:1990,goodwillie:1991,goodwillie:2003}) of the
identity functor $I$ on the category of based spaces. Such an operad
structure has been `known' intuitively by experts for some time but, as
far as the author knows, no explicit construction has previously been
given. One piece of evidence for such a structure is the calculation, due
to various people, of the homology of these derivatives. This homology is
the suspension of the standard Lie operad and so is itself an operad. It
is reasonable to ask, therefore, if there is an operad structure on the
derivatives themselves\footnote{The Goodwillie derivatives of a homotopy
functor are a sequence of spectra with actions by the symmetric groups,
but are only defined up to homotopy. By an operad structure on these
derivatives, we mean choices of models for these spectra in a suitable
symmetric monoidal category, such as the category of $S$--modules of
EKMM \cite{elmendorf/kriz/mandell/may:1997}, together with an operad
structure on those models.} that induces this structure on the homology.

Our construction is based on the partition poset model for the
derivatives $\partial_*I$ described by Arone and Mahowald in
\cite{arone/mahowald:1999}. They show that the derivatives are the
dual spectra associated to certain finite complexes known as the
partition poset complexes. In the present work we notice that these
complexes are precisely the simplicial bar construction\footnote{See,
for example, \cite[Section~II.2.3]{markl/shnider/stasheff:2002} for
the general form of the two-sided simplicial bar construction.} on the
operad $P$ in based spaces with $P(n) = S^0$ for all $n$. Most of the
paper is concerned with showing that such a bar construction has a
natural cooperad structure.\footnote{\label{foot:new}After this paper
was written, the author learnt that this result had already been
proved in unpublished work of Salvatore \cite{salvatore:1998} using an
alternative definition of the bar construction on an operad. See
Remark~\ref{rem:BV}.} We do this by reinterpreting the bar
construction in terms of spaces of trees. The cooperad structure then
comes from a natural way to break trees apart. Taking duals, we get
the required operad structure on the derivatives of the identity. In
fact, we can view the derivatives of the identity as a cobar
construction on the cooperad $Q$ in spectra with $Q(n) = S$, the
sphere spectrum, for all $n$.

In the final part of the paper (Section~\ref{sec:alg}) we show that
by taking homology we do indeed recover the `Lie' operad structure
on $H_*(\partial_*I)$. We do this by introducing spectral sequences
for calculating the homology of the topological bar and cobar
constructions. The $E^1$ terms of these spectral sequences can be
identified with algebraic versions of the bar and cobar constructions,
which in turn are related to the theory of Koszul duality for operads
introduced by Ginzburg and Kapranov in \cite{ginzburg/kapranov:1994}. Our
main result on this connection is that if the homology of a topological
operad $P$ is Koszul, then the homology of the bar construction $B(P)$
is its Koszul dual cooperad. In our case of interest, we deduce that
the induced operad structure on the homology of the derivatives of the
identity is that of the Koszul dual of the cocommutative cooperad. This
is precisely the `Lie' operad structure referred to above.

\subsection*{Outline of the paper}
We now give a more detailed description of the paper. The first two
sections are concerned with preliminaries. In Section~\ref{sec:monoidal}
we recall the notions of symmetric monoidal and enriched categories and
specify the categories we will be working with in this paper. These
are symmetric monoidal categories that are enriched, tensored and
cotensored over the category $\based$ of based compactly-generated
spaces (where $\based$ is a symmetric monoidal category with respect
to the smash product). It is to operads in these categories that we
refer in the title when we say `topological operads'. We also require
an extra condition that relates the symmetric monoidal structure to the
tensoring over $\based$. This condition (see Definition \ref{def:axiom})
is crucial to our later constructions. The two main examples of
categories satisfying our requirements are: based spaces themselves,
and a suitable symmetric monoidal category of spectra, such as that of
EKMM \cite{elmendorf/kriz/mandell/may:1997}.

In Section~\ref{sec:operads} we recall the definitions of operads and
cooperads. We should stress that the constructions of this paper apply
only to what we call \emph{reduced} operads and cooperads. These are $P$
with $P(0) = \basept$ and $P(1) = S$ the unit of the symmetric monoidal
structure. The bar construction can still be defined for more general
operads, but the cooperad structure described here does not seem to extend
to such cases. In this section we also define modules and comodules over
operads and cooperads respectively.

The real substance of the paper starts in Section~\ref{sec:trees}. Here
we define the trees that will form the combinatorial heart of our
description of the bar and cobar constructions. It is not a coincidence
that these trees are the same species used by, for example, Getzler and
Jones in their work \cite{getzler/jones:1994} on the bar constructions
for algebraic operads and Koszul duality. We also describe what we call
a \emph{weighting} on a tree (Definition \ref{def:weighting}), that is,
a suitable assignment of lengths to the edges of the tree. The spaces
$w(T)$ of weightings are at the heart of everything we do in this paper.

In Section~\ref{sec:bardef} we give our description of the bar
construction on an operad in terms of such trees. If $P$ is an operad of
based spaces, we can think of a point in the bar construction $B(P)$ as a
weighted tree (that is, a tree with lengths assigned to the edges) with
vertices labelled by points coming from the spaces $P(n)$. See Definition
\ref{def:bar(operad)} for a precise statement and Definition
\ref{def:formal_bar} for a more formal approach. In
Section~\ref{sec:simpbar} we show that what we have defined is isomorphic
to the standard simplicial bar construction on an operad.

In Section~\ref{sec:cooperad} we concern ourselves with the cooperad
structure on $B(P)$. This is given by the process of `ungrafting' trees
(see Definition \ref{def:grafting} and beyond). This involves taking a
weighted, labelled tree and breaking it up into smaller trees. Finding the
right way to weight and label these smaller trees gives us the required
cooperad structure maps.

One of the advantages of the way we have set up the theory is that the
cobar construction on a cooperad is strictly dual to the bar construction
on an operad. In Section~\ref{sec:cobar} we go through the definitions and
results dual to those of Section~\ref{sec:bar}.

The short section Section~\ref{sec:dual} is devoted to a simple but key
result (Proposition \ref{prop:duality}) that relates the bar and cobar
constructions via a duality functor that reduces to Spanier--Whitehead
duality in the case of spectra. This result says that, under the right
circumstances, the dual of the bar construction on an operad $P$ is
isomorphic to the cobar construction on the dual of $P$. This allows
us, later on, to identify the derivatives of the identity as the cobar
construction on a cooperad of spectra.

Before turning to our main example and application, we deal in
Section~\ref{sec:bar(modules)} with the two-sided bar and cobar
constructions. These include the bar construction for a module over an
operad and, dually, the cobar construction for a comodule over a cooperad.
To describe these requires a fairly simple generalization of much of the
work we did in Sections~\ref{sec:trees}--\ref{sec:bar}, in particular,
a more general notion of tree (see Definition \ref{def:gen_trees}).

Finally, in Section~\ref{sec:application} we are able to complete the
main aim of this paper. We identify the partition poset complexes with a
bar construction and deduce the existence of an operad structure on the
derivatives of the identity functor (Corollary \ref{cor:operad}). We
also give examples of modules over the resulting operad, including,
in particular, a module $M_X$ naturally associated to a based space $X$.

The last section of the paper Section~\ref{sec:alg} is concerned with the
relationship of our work to the algebraic bar construction and Koszul
operads. As promised, we construct a spectral sequence (Proposition
\ref{prop:specseq}) relating the two and deduce the result on Koszul
duality (Proposition \ref{prop:koszul}).

\subsection*{Future Work}
The work of this paper raises various questions that seem to the author
to warrant further attention:
\begin{itemize}
\item What is the homotopy theory of the topological bar and cobar
constructions? In particular, how do they relate to known model structures
on the categories of operads and cooperads (see, for example,
Berger--Moerdijk \cite{berger/moerdijk:2003})?
\item Is there a deeper relationship between Goodwillie's homotopy
calculus and the theory of operads? The present paper does not do any
calculus, the only connection being via the partition poset complexes. One
might ask, for example, if the derivatives of other functors can be
described and/or treated using these ideas.
\item What object is described by an algebra or module over the operad
formed by the derivatives of the identity? In Remark \ref{rem:modules}
we show that a based space $X$ gives rise to such a module. How much of
(the homotopy theory of) the space $X$ is retained by this module?
\end{itemize}

\subsection*{Acknowledgements}
The work of this paper forms the author's PhD thesis written at
the Massachusetts Institute of Technology under the supervision of
Haynes Miller, to whom the greatest thanks are due for his constant
support, encouragement and advice. The idea that the derivatives of
the identity might be related to a cobar construction was suggested by
work of Kristine Bauer, Brenda Johnson and Jack Morava. The observation
that the partition poset complexes (and hence the derivatives of the
identity) can be described in terms of spaces of trees was mentioned to
the author by Tom Goodwillie, who heard it from Greg Arone. The work of
Benoit Fresse \cite{fresse:2004} on the algebraic side of the theory was
invaluable to the present paper. The author has also benefited greatly
from conversations with Mark Behrens and Andrew Mauer-Oats while writing
this paper, and finally would like to thank the referee for some helpful
comments and suggestions.

\section{Symmetric monoidal and enriched categories} \label{sec:monoidal}
On the one hand, the bar and cobar constructions are most easily defined
(and understood) in the category of based spaces. On the other hand,
our main application is in a category of spectra. We will develop the
theory in a general setting that encompasses both cases. This approach
will also allow us to appreciate more readily the duality between the
bar and cobar constructions.

In this section we recall the basic theory of symmetric monoidal and
enriched categories (see \cite[Section~6]{borceux:1994(2)} for a detailed
account). We state precisely (Definition \ref{def:axiom}) the structure we
will require of a category to make the bar and cobar constructions in it.
The only material in this chapter that is not standard is the definition
of enriched symmetric monoidal categories or `symmetric monoidal
$\cat{V}$--categories' as we have called them (Definition
\ref{def:axiom}). The `distributivity' morphism described there is a
key component of the constructions made later in the paper and so we
draw the reader's attention to it now.

\begin{definition}[(Symmetric monoidal categories)] \label{def:sym mon}
A \emph{monoidal category} consists of
\begin{itemize}
\item a (locally small) category $\cat{V}$,
\item a functor $- \smsh - \co \cat{V} \times \cat{V} \to \cat{V}$,
\item a \emph{unit} object $I$ in $\cat{V}$ together with natural
isomorphisms $X \smsh I \isom X \isom I \smsh X$,
\item a natural \emph{associativity} isomorphism $X \smsh (Y \smsh Z)
\isom (X \smsh Y) \smsh Z$,
\end{itemize}
such that the appropriate three coherence diagrams commute
\cite[Section~VII]{maclane:1971}. A \emph{symmetric monoidal category}
is a monoidal category together with
\begin{itemize}
\item a natural \emph{symmetry} isomorphism $X \smsh Y \isom Y \smsh X$,
\end{itemize}
such that four additional coherence diagrams also commute. We will
denote such a symmetric monoidal category by $(\cat{V},\smsh,I)$, or
just $\cat{V}$ with the rest of the structure understood.
\end{definition}

\begin{remark} \label{rem:symmetry}
We will not give names to the associativity and symmetry isomorphisms in a
symmetric monoidal category. When we write unbracketed expressions such as
\[ X \smsh Y \smsh Z \]
or unordered expressions such as
\[ \bigwedge_{a \in A} X_a \]
we mean any one particular choice of ordering and bracketing. Different
choices are related by the appropriate associativity and commutativity
isomorphisms between them. A map to or from a particular choice
determines a map to or from any other choice by composing with the
relevant isomorphism.
\end{remark}

\begin{definition} \label{def:closed sym mon}
A \emph{closed symmetric monoidal category} is a symmetric monoidal
category $(\cat{V},\smsh,I)$ together with a functor
\[ \cat{V}^\text{op} \times \cat{V} \to \cat{V}; \: (X,Y) \mapsto \Map(X,Y) \]
and a natural isomorphism of sets
\[ \Hom_{\cat{V}}(X \smsh Y,Z) \isom \Hom_{\cat{V}}(X,\Map(Y,Z)), \]
where $\Hom_{\cat{V}}(X,Y)$ is the set of morphisms from $X$ to $Y$
in the category $\cat{V}$.
\end{definition}

\begin{remark}
The natural isomorphism of sets in Definition \ref{def:closed sym mon}
can be made into an isomorphism within $\cat{V}$. That is, in any closed
symmetric monoidal category there is a natural isomorphism
\[ \Map(X \smsh Y, Z) \isom \Map(X, \Map(Y,Z)). \]
See \cite[Section~6.5.3]{borceux:1994(2)} for details.
\end{remark}

\begin{definition}[(Enriched categories)]
Let $(\cat{V},\smsh,I)$ be a given closed symmetric monoidal category. A
\emph{$\cat{V}$--category} or \emph{category enriched over $\cat{V}$}
consists of
\begin{itemize}
\item a class $\cat{C}$,
\item for each pair of elements $C,D \in \cat{C}$, an object
$\Map_{\cat{V}}(C,D)$ of $\cat{V}$,
\item composition morphisms
\[ \Map_{\cat{V}}(C,D) \smsh \Map_{\cat{V}}(D,E) \to \Map_{\cat{V}}(C,E)
\]
for each $C,D,E \in \cat{C}$,
\item identity morphisms
\[ I \to \Map_{\cat{V}}(C,C) \]
for each $C \in \cat{C}$,
\end{itemize}
that satisfy the appropriate conditions \cite[Section~6.2.1]{borceux:1994(2)}. We
will denote such a $\cat{V}$--category by $\cat{C}$ with the rest of
the structure understood.
\end{definition}

\begin{remark} \label{rem:enriched}
We include some basic observations about enriched categories from
\cite[Section~6.2]{borceux:1994(2)}.
\begin{enumerate}
\item Let $(\mathsf{Set},\times,\basept)$ be the symmetric monoidal
category of sets under cartesian product. A $\mathsf{Set}$--category is
then the same thing as a (locally small) category.
\item A $\cat{V}$--category $\cat{C}$ has an underlying category whose
objects are the elements of $\cat{C}$ and whose morphisms $C \to D$ are
the elements of the set $\Hom_{\cat{V}}(I,\Map_{\cat{V}}(C,D))$, where $I$
is the unit object of $\cat{V}$. We often therefore think of a
$\cat{V}$--category $\cat{C}$ as a normal category with extra structure
given by the objects $\Map_{\cat{V}}(C,D)$.
\item A closed symmetric monoidal category $\cat{V}$ is enriched over
itself with
\[ \Map_{\cat{V}}(X,Y) := \Map(X,Y). \]
\end{enumerate}
\end{remark}

\begin{definition}[(Tensoring and cotensoring)]
Let $\cat{C}$ be a $\cat{V}$--category. A \emph{tensoring} of $\cat{C}$
over $\cat{V}$ is a functor
\[ \cat{V} \times \cat{C} \to \cat{C}; \; (X,C) \mapsto X \otimes C \]
together with a natural isomorphism
\[ \Map_{\cat{V}}(X \otimes C, D) \isom \Map(X,\Map_{\cat{V}}(C,D)). \]
A category $\cat{C}$ \emph{tensored} over $\cat{V}$ is a
$\cat{V}$--category together with a chosen tensoring.

A \emph{cotensoring} of $\cat{C}$ over $\cat{V}$ is a functor
\[ \cat{V}^\text{op} \times \cat{C} \to \cat{C}; \; (X,D) \mapsto
\Map_{\cat{C}}(X,D) \]
together with a natural isomorphism
\[ \Map_{\cat{V}}(C, \Map_{\cat{C}}(X,D)) \isom
\Map(X,\Map_{\cat{V}}(C,D)). \]
A category $\cat{C}$ \emph{cotensored} over $\cat{V}$ is a
$\cat{V}$--category together with a chosen cotensoring.
\end{definition}

\begin{remark} \label{rem:tensor}
Here are some basic observations about tensorings and cotensorings.
\begin{enumerate}
\item A closed symmetric monoidal category $(\cat{V},\smsh,I)$ is
tensored and cotensored over itself with $X \otimes Y := X \smsh Y$
and $\Map_{\cat{V}}(X,Y) := \Map(X,Y)$.
\item If $\cat{C}$ is tensored over $\cat{V}$, we have natural
isomorphisms
\[ (X \smsh Y) \otimes C \isom X \otimes (Y \otimes C) \]
for $X,Y \in \cat{V}$ and $C \in \cat{C}$. If $\cat{C}$ is cotensored
over $\cat{V}$, we have natural isomorphisms
\[ \Map_{\cat{C}}(X \smsh Y, C) \isom \Map_{\cat{C}}(X,
\Map_{\cat{C}}(Y,C)) \]
for $X,Y \in \cat{V}$ and $C \in \cat{C}$.
\end{enumerate}
\end{remark}

\begin{proposition} \label{prop:dual1}
Let $\cat{C}$ be a $\cat{V}$--category. Then $\cat{C}^\text{op}$ has
a natural enrichment over $\cat{V}$.\footnote{Here $\cat{C}^\text{op}$
denotes the opposite category of the category underlying $\cat{C}$
described in Remark \ref{rem:enriched}(2).} If $\cat{C}$ is tensored,
then $\cat{C}^\text{op}$ is naturally cotensored and vice versa.
\end{proposition}
\begin{proof}
We define an enrichment on $\cat{C}^\text{op}$ by
\[ \Map_{\cat{V}}(C^\text{op},D^\text{op}) := \Map_{\cat{V}}(D,C) \]
where $C^\text{op}$ is the object in $\cat{C}^\text{op}$ corresponding to $C \in
\cat{C}$. If $- \otimes -$ is a tensoring for $\cat{C}$ then we get a
cotensoring for $\cat{C}^\text{op}$ by setting
\[ \Map_{\cat{C}^\text{op}}(X,D^\text{op}) := (X \otimes D)^\text{op}. \]
The required natural isomorphism comes from
\[ \begin{split}
    \Map_{\cat{V}}(C^\text{op}, \Map_{\cat{C}^\text{op}}(X,D^\text{op})) &=
    \Map_{\cat{V}}(X \otimes D, C) \\
        &\isom \Map(X,\Map_{\cat{V}}(D,C)) \\
        &= \Map(X,\Map_{\cat{V}}(C^\text{op},D^\text{op})). \\
\end{split} \]
The vice versa part is similar.
\end{proof}

We are interested in categories that both are themselves symmetric
monoidal categories and are enriched over another symmetric monoidal
category. The following definition contains the properties of these that
we require in this paper.

\begin{definition} \label{def:axiom}
Let $(\cat{V},\smsh,I)$ be a closed symmetric monoidal category. A
\emph{symmetric monoidal $\cat{V}$--category} consists of
\begin{itemize}
\item a symmetric monoidal category $(\cat{C},\barwedge,S)$ with $\cat{C}$
enriched, tensored and cotensored over $\cat{V}$,
\item a natural transformation
\[ d\co (X \smsh Y) \otimes (C \barwedge D) \to (X \otimes C) \barwedge
(Y \otimes D) \]
\end{itemize}
satisfying the following axioms:
\begin{itemize}
\item (Associativity) The diagram
\[ \begin{diagram} \dgARROWLENGTH=2.4em
    \node{(X \smsh Y \smsh Z) \otimes (C \barwedge D \barwedge E)}
    \arrow{e,t}{d} \arrow{s,l}{d}
        \node{((X \smsh Y) \otimes (C \barwedge D)) \barwedge (Z \otimes
        E)} \arrow{s,r}{id \, \barwedge \, d} \\
    \node{(X \otimes C) \barwedge ((Y \smsh Z) \otimes (D \barwedge E))}
    \arrow{e,t}{id \, \barwedge \, d}
        \node{(X \otimes C) \barwedge (Y \otimes D) \barwedge (Z
        \otimes E)}
\end{diagram} \]
commutes for all $X,Y,Z \in \cat{V}$ and $C,D,E \in \cat{C}$.
\item (Unit)
The composite
\[ X \otimes C  \isom (X \smsh I) \otimes (C \barwedge S) \arrow{e,t}{d}
(X \otimes C) \barwedge (I \otimes S) \isom X \otimes C \]
is the identity, for any $X \in \cat{V}$ and $C \in \cat{C}$. Recall
that $I,S$ are the units of the symmetric monoidal structures on $\cat{V},
\cat{C}$ respectively.
\end{itemize}
\end{definition}

The transformation $d$ (for `distribute') is our way of relating the
symmetric monoidal structures in the two categories. It will be essential
in constructing the cooperad structure on the bar construction of an
operad (see Definition \ref{def:formal_cooperad_maps}).

\begin{remark}
A closed symmetric monoidal category $\cat{V}$ is itself a symmetric
monoidal $\cat{V}$--category with the transformation $d$ given by the
symmetry and associativity isomorphism:
\[ (X \smsh Y) \smsh (C \smsh D) \isom (X \smsh C) \smsh (Y \smsh D) \]
\end{remark}

\begin{prop} \label{prop:dual}
Let $\cat{C}$ be a symmetric monoidal $\cat{V}$--category. Then
$\cat{C}^\text{op}$ is naturally also a symmetric monoidal $\cat{V}$--category.
\end{prop}
\proof
We already know from Proposition \ref{prop:dual1} that $\cat{C}^\text{op}$
is enriched, tensored and cotensored over $\cat{V}$ and there is a
canonical symmetric monoidal structure on $\cat{C}^\text{op}$ given by
that on $\cat{C}$. It therefore only remains to construct the map
$d$. The tensoring in $\cat{C}^\text{op}$ is given by the cotensoring in
$\cat{C}$. Therefore $d$ for $\cat{C}^\text{op}$ corresponds to the following
map in $\cat{C}$:
\[ \Map_{\cat{C}}(X,C) \barwedge \Map_{\cat{C}}(Y,D) \to \Map_{\cat{C}}(X
\smsh Y, C \barwedge D) \]
This is adjoint to a map
\[ (X \smsh Y) \otimes (\Map_{\cat{C}}(X,C) \barwedge \Map_{\cat{C}}(Y,D))
\to C \barwedge D \]
constructed by first using $d$ for $\cat{C}$ to get to
\[ (X \otimes \Map_{\cat{C}}(X,C)) \barwedge (Y \otimes
\Map_{\cat{C}}(Y,D)) \]
and then using the evaluation maps
\[ X \otimes \Map_{\cat{C}}(X,C) \to C \text{ and } Y \otimes
\Map_{\cat{C}}(Y,D) \to D. \eqno{\Box}\]

An important property of the categories that we work with in this paper
is that they are \emph{pointed}, that is, they have a \emph{null} object
$\basept$ that is both initial and terminal. The following proposition
describes how null objects interact with symmetric monoidal structures
and enrichments.

\begin{prop} \label{prop:null}
Let $(\cat{V},\smsh,I)$ be a closed symmetric monoidal category that is
pointed with null object $\basept$. Then
\[ \basept \smsh X \isom \basept \isom \Map(\basept,X) \isom
\Map(X,\basept) \]
for all $X \in \cat{V}$.

Moreover, let $\cat{C}$ be a category enriched over $\cat{V}$. If
$\cat{C}$ is tensored then $\basept \otimes C$ is an initial object
in $\cat{C}$ for all $C \in \cat{C}$. If $\cat{C}$ is cotensored then
$\Map_{\cat{C}}(\basept,D)$ is a terminal object in $\cat{C}$ for all
$D \in \cat{C}$.

Finally, if $\cat{C}$ is both tensored and cotensored over $\cat{V}$,
then the initial and terminal objects are isomorphic and so $\cat{C}$
is itself pointed.
\end{prop}
\begin{proof}
We observe that
\[ \Hom_{\cat{V}}(\basept \smsh X, Y) \isom
\Hom_{\cat{V}}(\basept,\Map(X,Y)) \]
which has one element for any $X,Y$. This tells us that $\basept \smsh X$
is initial and hence isomorphic to $\basept$. The other isomorphisms in
the first part of the proposition are similar.

Next, the tensoring functor $- \otimes C \co \cat{V} \to \cat{C}$ is a left
adjoint so preserves an initial object. Dually, the cotensoring functor
$\Map_{\cat{C}}(-,D) \co \cat{V}^\text{op} \to \cat{C}$ is a right adjoint so
preserves the terminal object. This gives us the second part.

Finally, if $\cat{C}$ is both tensored and cotensored, we get a map from
the terminal object to the initial object by
\[ \Map_{\cat{C}}(\basept,D) \to I \otimes \Map_{\cat{C}}(\basept,D)
\to \basept \otimes \Map_{\cat{C}}(\basept,D). \]
The first map here is an example of a general isomorphism $C \to I \otimes
C$ where $I$ is the unit object of $\cat{V}$. The second map comes from
$I \to \basept$. A map from a terminal object to an initial object must
be an isomorphism. Therefore $\cat{C}$ is pointed.
\end{proof}

\begin{examples} \label{ex:categories}
The categories with which we will mainly be concerned in this paper are
the following.
\begin{enumerate}
\item Let $\based$ be the category of compactly generated based spaces and
basepoint-preserving continuous maps of \cite{lewis/may/steinberger:1986}.
Then $\based$ is a pointed closed symmetric monoidal category under the
usual smash product $\smsh$, with unit the $0$--sphere $S^0$ and
$\Map(X,Y)$ equal to the space of basepoint-preserving maps $X \to Y$.
\item Let $\spectra$ be the category of $S$--modules of EKMM
\cite{elmendorf/kriz/mandell/may:1997}. Then $(\spectra,\smsh_S,S)$ is a
symmetric monoidal $\based$--category, where $S$ is the sphere spectrum
and $\smsh_S$ is the smash product of $S$--modules
\cite[Section~II.1.1]{elmendorf/kriz/mandell/may:1997}. The
enrichment, tensoring and cotensoring are described in
\cite[Section~VII.2.8]{elmendorf/kriz/mandell/may:1997}. For the distributivity
map $d$ we have a natural isomorphism
\[ d\co (X \smsh Y) \smsh (E \smsh_S F) \arrow{e,t}{\isom} (X \smsh E)
\smsh_S (Y \smsh F) \]
given by the fact that $X \smsh E \isom (X \smsh S) \smsh_S E$ (see
\cite[Section~II.1.4]{elmendorf/kriz/mandell/may:1997}).
\end{enumerate}
We will usually work with a general symmetric monoidal $\based$--category
denoted $(\cat{C},\barwedge,S)$, but these examples will be foremost in
our minds.
\end{examples}

\section{Operads and cooperads} \label{sec:operads}

In this section $(\cat{C},\barwedge,S)$ denotes a pointed symmetric
monoidal category with null object $\basept$. We will assume that
$\cat{C}$ has all necessary limits and colimits and write the coproduct
in $\cat{C}$ as a wedge product using $\vee$.

\begin{definition}[(Symmetric sequences)] \label{def:symseq}
A \emph{symmetric sequence} in $\cat{C}$ is a functor $F$ from the
category of nonempty finite sets and bijections to $\cat{C}$. For
each nonempty finite set $A$, the symmetric group $\Sigma_A$ acts on
$F(A)$. We will write $F(n)$ for $F(\{1,\dots,n\})$. Note that our
symmetric sequences (and hence our operads) do not have an $F(0)$ term
because our indexing sets are nonempty. We will often write `finite set'
when we mean `nonempty finite set' and these will usually be labelled
$A,B,\dots$. We write $\cat{C}^{\Sigma}$ for the category of symmetric
sequences in $\cat{C}$ (whose morphisms are the natural transformations).
\end{definition}

There are several different but equivalent ways to define operads (see
Markl--Shnider--Stasheff \cite{markl/shnider/stasheff:2002} for a
comprehensive guide). We will use the following definition.

\begin{definition}[(Operads)] \label{def:operad}
An \emph{operad} in the symmetric monoidal category
$(\cat{C},\barwedge,S)$ is a symmetric sequence $P$ together with
\emph{partial composition maps}
\[ - \circ_a - \co P(A) \barwedge P(B) \to P(A \cup_a B) \]
for each pair of finite sets $A,B$, and each $a \in A$ (where $A \cup_a B$
denotes $(A \setminus \{a\}) \amalg B$), and a \emph{unit map}
\[ \eta \co S \to P(1). \]
The composition maps must be natural in $A$ and $B$ and must satisfy
the following four axioms:
\begin{enumerate}
\item  The diagram
\[ \begin{diagram}
    \node{P(A) \barwedge P(B) \barwedge P(C)} \arrow{e,t}{id \, \barwedge
    \, \circ_b} \arrow{s,l}{\circ_a \, \barwedge \, id} \node{P(A)
    \barwedge P(B \cup_b C)} \arrow{s,r}{\circ_a} \\
    \node{P(A \cup_a B) \barwedge P(C)} \arrow{e,t}{\circ_b} \node{P(A
    \cup_a B \cup_b C)}
\end{diagram} \]
commutes for all $a \in A$ and $b \in B$. (Notice that $(A \cup_a B)
\cup_b C = A \cup_a (B \cup_b C)$.)
\item  The diagram
\[ \begin{diagram}
    \node{P(A) \barwedge P(B) \barwedge P(C)} \arrow{s,l}{\isom}
    \arrow{e,t}{\circ_a \, \barwedge \, id}  \node{P(A \cup_{a} B)
    \barwedge P(C)} \arrow[2]{s,r}{\circ_{a'}} \\
    \node{P(A) \barwedge P(C) \barwedge P(B)} \arrow{s,l}{\circ_{a'} \,
    \barwedge \, id} \\
    \node{P(A \cup_{a'} C) \barwedge P(B)} \arrow{e,t}{\circ_a} \node{P(A
    \cup_a B \cup_{a'} C)}
\end{diagram} \]
commutes for all $a \neq a' \in A$. (Notice that $(A \cup_a B) \cup_{a'}
C = (A \cup_{a'} C) \cup_a B$.)
\item  The diagram
\[ \begin{diagram}
    \node{P(A)} \arrow{e,t}{\eta \, \barwedge \, id} \arrow{se,b}{id}
    \node{P(1) \barwedge P(A)} \arrow{s,r}{\circ_1} \\
    \node[2]{P(\{1\} \cup_{1} A)}
\end{diagram} \]
commutes for all $A$.
\item  The diagram
\[ \begin{diagram}
    \node{P(A)} \arrow{e,t}{id \, \barwedge \, \eta} \arrow{se,b}{\isom}
    \node{P(A) \barwedge P(1)} \arrow{s,r}{\circ_a} \\
    \node[2]{P(A \cup_a \{1\})}
\end{diagram} \]
commutes for all $a \in A$. (The diagonal map here is induced by the
obvious bijection $A \to A \cup_a \{1\}$.)
\end{enumerate}
A \emph{morphism} of operads $P \to P'$ is a morphism of symmetric
sequences that commutes with the composition and unit maps.
\end{definition}

\begin{definition} \label{def:reduced}
An \emph{augmentation} of an operad $P$ is a map $\varepsilon\co P(1)
\to S$ such that the composite
\[ \begin{diagram} \node{S} \arrow{e,t}{\eta} \node{P(1)}
\arrow{e,t}{\varepsilon} \node{S} \end{diagram} \]
is the identity on $S$. An \emph{augmented operad} is an operad together
with an augmentation. An operad $P$ is \emph{reduced} if the unit map
$\eta\co S \to P(1)$ is an isomorphism. A reduced operad has a unique
augmentation given by the inverse of the unit map. A \emph{morphism}
of augmented operads is a morphism of operads that commutes with the
augmentation.
\end{definition}

\begin{remark} \label{rem:operad}
Operads are a generalization of monoids for the symmetric monoidal
category $(\cat{C},\barwedge,S)$. A monoid $X$ in $\cat{C}$ gives rise
to an operad $P_X$ with $P_X(1) = X$ and $P_X(n) = \basept$ for $n >
1$. Conversely, given an operad $P$ in the symmetric monoidal category
$\cat{C}$, $P(1)$ is a monoid in $\cat{C}$.
\end{remark}

An alternative definition of an operad is based on a monoidal structure
on the category of symmetric sequences. We define this monoidal structure
now.

\begin{definition}[(Composition product of symmetric sequences)]
\label{def:compprod}
Let the \emph{composition product} of the two symmetric sequences
$M,N$ be the symmetric sequence $M \circ N$ with
\[ (M \circ N)(A) := \bigvee_{A = \coprod_{j \in J} A_j} M(J) \barwedge
\Barwedge_{j \in J} N(A_j). \]
The coproduct here is taken over all unordered partitions of $A$ into
a collection of nonempty subsets $\{A_j\}_{j \in J}$. The particular
choice of indexing set is not important in the sense that we do not sum
over different $J$ that index the same partition. A bijection $A \to A'$
determines a bijection between partitions of $A$ and partitions of $A'$
in an obvious way. Thus we match up the terms in the coproducts that
define $(M \circ N)(A)$ and $(M \circ N)(A')$. If $J$ and $J'$ index
two corresponding partitions of $A$ and $A'$ respectively, then we get a
natural choice of bijection $J \to J'$. Moreover, if $j \in J$ and $j'
\in J'$ correspond under this bijection then we get a bijection $A_j
\to A'_{j'}$ by restricting the bijection $A \to A'$. The actions of $M$
and $N$ on these bijections together give us an isomorphism
\[ (M \circ N)(A) \to (M \circ N)(A'). \]
Thus $M \circ N$ becomes a symmetric sequence in $\cat{C}$.
\end{definition}

\begin{definition}
The \emph{unit symmetric sequence} in the pointed symmetric monoidal
category $(\cat{C},\barwedge,S)$ is the symmetric sequence $I$ given by
\[ I(A) := \begin{cases} S & \text{if $|A| = 1$}; \\ \basept &
\text{otherwise}; \end{cases} \]
where $\basept$ is the null object of $\cat{C}$.
\end{definition}

\begin{lemma} \label{lem:circ_unit}
Let $(\cat{C},\barwedge,S)$ be a pointed symmetric monoidal category. Then
for any symmetric sequence $M$ there are natural isomorphisms
\[ M \circ I \isom M \isom I \circ M. \]
\end{lemma}
\begin{proof}
For the finite set $A$, the only term that contributes to $(M \circ I)(A)$
comes from the partition of $A$ into singleton subsets. This makes it
clear that $M \circ I \isom M$. The only term that contributes to $(I
\circ M)(A)$ comes from the trivial partition of $A$ into one subset,
that is $A$ itself. From this we see that $I \circ M \isom M$.
\end{proof}

To get a monoidal structure on the category of symmetric sequences, we
also need an associativity isomorphism. This does not exist in general,
although it does in the case of the following lemma.

\begin{lemma} \label{lem:circ_assoc}
Let $(\cat{C},\barwedge,S)$ be a pointed symmetric monoidal category
in which $\barwedge$ commutes with finite coproducts. Then there are
natural isomorphisms
\[ L \circ (M \circ N) \isom (L \circ M) \circ N \]
for symmetric sequences $L,M,N$ in $\cat{C}$.
\end{lemma}
\begin{proof}
Using the hypothesis that $\barwedge$ commutes with finite coproducts,
it is not hard to see that each side is naturally isomorphic to the
symmetric sequence $(L \circ M \circ N)$ given by
\[ (L \circ M \circ N)(A) := \bigvee_{A = \coprod_{b \in B} A_b, \; B =
\coprod_{c \in C} B_c} L(C) \barwedge \Barwedge_{c \in C} M(B_c)
\barwedge \Barwedge_{b \in B} N(A_b). \]
The coproduct here is over all partitions of $A$ into nonempty subsets
indexed by some set $B$, together with a partition of $B$ into subsets
indexed by some $C$. Equivalently, the coproduct is indexed of pairs
of partitions of $A$, one (indexed by $B$) a refinement of the other
(indexed by $C$).
\end{proof}

The following description of operads is due to Smirnov. See
\cite[Theorem~1.68]{markl/shnider/stasheff:2002} for further details.

\begin{proposition} \label{prop:comp_monoidal}
Let $(\cat{C},\barwedge,S)$ be a pointed symmetric monoidal category in
which $\barwedge$ commutes with finite coproducts. Then the composition
product $\circ$ is a monoidal product on the category of symmetric
sequences in $\cat{C}$ with unit object $I$ and unit and associativity
isomorphisms given by Lemmas \ref{lem:circ_unit} and \ref{lem:circ_assoc}
respectively. In this case, an operad in $\cat{C}$ is precisely a monoid
for this monoidal product.
\end{proposition}
\begin{proof}
One can easily check that the axioms for a monoidal structure are
satisfied. If $P$ is an operad in $\cat{C}$, the operad compositions
make up a map
\[ P \circ P \to P \]
and the unit map $\eta$ gives a map of symmetric sequences
\[ I \to P. \]
The operad axioms then translate into associativity and unit axioms that
give $P$ the structure of a monoid under $\circ$.
\end{proof}

\begin{remark} \label{rem:triple_product}
If $\cat{C}$ is a \emph{closed} symmetric monoidal category then
$\barwedge$ has a right adjoint and so preserves all colimits. In
particular, the hypothesis of Lemma \ref{lem:circ_assoc} holds and so
we get a true monoidal structure on the symmetric sequences in $\cat{C}$.

Unfortunately, even when $\cat{C}$ is closed symmetric monoidal, its
opposite category $\cat{C}^\text{op}$ (with the standard symmetric monoidal
structure) is unlikely to be closed. Since we will want to dualize most
of the results of this paper to be able to deal with cooperads as well
as operads, we need to get round this hypothesis. For this, we notice
that in general there are natural maps of symmetric sequences
\[ (L \circ M \circ N) \to L \circ (M \circ N) \]
and
\[ (L \circ M \circ N) \to (L \circ M) \circ N \]
where $(L \circ M \circ N)$ is defined as in the proof of Lemma
\ref{lem:circ_assoc}. In general these are not isomorphisms so we do not
get a monoidal structure on the category of symmetric sequences. However, it is possible to define monoids in this more general case (see \cite{ching:2005c} for more details), and we get the
following alternative characterization of an operad.
\end{remark}

\begin{proposition} \label{prop:comp_general}
Let $(\cat{C},\barwedge,S)$ be a pointed symmetric monoidal category. An
operad in $\cat{C}$ is equivalent to a symmetric sequence $P$ together
with maps
\[ m \co P \circ P \to P ; \; \eta\co I \to P \]
of symmetric sequences such that the following diagrams commute:
\begin{enumerate}
\item Associativity:
\[ \begin{diagram} \node[2]{(P \circ P) \circ P} \arrow{e,t}{m \, \circ \,
id} \node{P \circ P} \arrow{se,t}{m} \\
            \node{(P \circ P \circ P)} \arrow{ne} \arrow{se} \node[3]{P}
            \\
            \node[2]{P \circ (P \circ P)} \arrow{e,t}{id \, \circ \, m}
            \node{P \circ P} \arrow{ne,t}{m}
\end{diagram} \]
where the two initial arrows are the maps mentioned in
Remark~\ref{rem:triple_product}.
\item Left unit:
\[ \begin{diagram} \node{P} \arrow{se,b}{id} \arrow{e,t}{id \, \circ \,
\eta} \node{P \circ P} \arrow{s,r}{m} \\ \node[2]{P} \end{diagram} \]
\item Right unit:
\[ \begin{diagram} \node{P} \arrow{se,b}{id} \arrow{e,t}{\eta \, \circ \,
id} \node{P \circ P} \arrow{s,r}{m} \\ \node[2]{P} \end{diagram} \]
\end{enumerate}
\end{proposition}

\begin{remark}
We will refer to an operad $P$ as a \emph{monoid} with respect to
the composition product, even when we do not in fact have a monoidal
structure. There are similarly defined notions of an object with a right
or left action of a monoid in this generalized setting. These give us
right and left modules over our operads which we now define.
\end{remark}

\begin{definition}[(Modules over operads)] \label{def:module}
A \emph{left module} over the operad $P$ is a symmetric sequence $M$
together with a left action of the monoid $P$, that is, a map
\[ l\co P \circ M \to M \]
such that the diagrams
\[ \begin{diagram} \node[2]{(P \circ P) \circ M} \arrow{e,t}{m \, \circ \,
id} \node{P \circ M} \arrow{se,t}{l} \\
            \node{(P \circ P \circ M)} \arrow{se} \arrow{ne} \node[3]{M}
            \\
            \node[2]{P \circ (P \circ M)} \arrow{e,t}{id \, \circ \, l}
            \node{P \circ M} \arrow{ne,t}{l}
\end{diagram} \]
and
\[ \begin{diagram} \node{M} \arrow{e,t}{\eta \, \circ \, id}
\arrow{se,b}{id} \node{P \circ M} \arrow{s,r}{l} \\ \node[2]{M}
\end{diagram} \]
commute.

A \emph{right module} over $P$ is a symmetric sequence $M$ with a right
action of $P$, that is a map
\[ M \circ P \to M \]
satisfying corresponding axioms. A \emph{$(P,P)$--bimodule} is a symmetric
sequence $M$ that is both a right and a left module over $P$ such that
\[ \begin{diagram} \node[2]{(P \circ M) \circ P} \arrow{e} \node{M \circ
P} \arrow{se} \\
            \node{P \circ M \circ P} \arrow{ne} \arrow{se} \node[3]{M} \\
            \node[2]{P \circ (M \circ P)} \arrow{e} \node{P \circ M}
            \arrow{ne}
\end{diagram} \]
commutes. Clearly, $P$ itself is a $(P,P)$--bimodule.
\end{definition}

\begin{remark}
It's useful to have a slightly more explicit description of a module over
an operad. The action map for a left $P$--module $M$ consists of maps
\[ P(r) \barwedge M(A_1) \barwedge \dots \barwedge M(A_r) \to M(A) \]
for every partition $A = \coprod_{i = 1}^{r} A_i$ of a finite set $A$ into
nonempty subsets. Conversely, giving maps of this form that satisfy
appropriate conditions uniquely determines a left $P$--module. Similarly,
a right module structure consists of maps of the form
\[ M(r) \barwedge P(A_1) \barwedge \dots \barwedge P(A_r) \to M(A). \]
\end{remark}

\begin{remark}
In the same way that operads are a generalization of monoids in $\cat{C}$,
modules over those operads are generalization of modules over the
monoids. A module $M$ over the monoid $X$ gives rise to a module $P_M$
over the operad $P_X$ described in Remark \ref{rem:operad}, with $P_M(n)
= \basept$ if $n > 1$ and $P_M(1) = M$.
\end{remark}

\begin{remark}
An augmentation for the operad $P$ is equivalent to either a left or
right module structure on the unit symmetric sequence $I$.
\end{remark}

The standard notion of an algebra over an operad is closely related to
that of a module. We briefly describe how this works.

\begin{definition}[(Algebras over an operad)]
An \emph{algebra} over the operad $P$ is an object $C \in \cat{C}$
together with maps
\[ P(A) \barwedge \Barwedge_{a \in A} C \to C \]
that satisfy appropriate naturality, associativity and unit axioms.
\end{definition}

The following result allows us to construct a left $P$--module from a
$P$--algebra.\footnote{There is a more basic way to view algebras over an
operad as modules. This requires us to introduce an $M(0)$ term to our
modules (that is, our symmetric sequences become functors from the
category of all finite sets, not just nonempty finite sets). With a
corresponding generalization of the composition product, and hence of the
notion of module, a $P$--algebra is equivalent to a left $P$--module
concentrated in the $M(0)$ term. The reason we do not allow our modules to
have this extra term is that the comodule structure on the bar
construction (see Section~\ref{sec:two-sided}) would not then exist
in general.}

\begin{lemma} \label{lem:algebra}
Let $C$ be an algebra over the operad $P$. Then there is a natural left
$P$--module structure on the constant symmetric sequence $\underline{C}$
with $\underline{C}(A) = C$ for all finite sets $A$.\footnote{The obvious
converse to this Lemma is not true. That is, a constant symmetric sequence
together with a left $P$--module structure need not arise from a
$P$--algebra. The construction given in the proof of this lemma forces
different components of the module structure map to be the same which
need not be the same in general.}
\end{lemma}
\proof
The components of the module structure map $P \circ \underline{C} \to
\underline{C}$ are given by the algebra structure maps as follows:
\[ P(r) \barwedge \underline{C}(A_1) \barwedge \dots \barwedge
\underline{C}(A_r) = P(r) \barwedge C^{\barwedge r} \to C =
\underline{C}(A) \eqno{\Box} \]

\begin{definition}[(Cooperads)] \label{def:cooperad}
The notion of a cooperad is dual to that of an operad. That is,
a \emph{cooperad} in $\cat{C}$ is an operad in the opposite category
$\cat{C}^\text{op}$ with the canonical symmetric monoidal structure determined
by that in $\cat{C}$. More explicitly, a cooperad consists of a symmetric
sequence $Q$ in $\cat{C}$ together with \emph{cocomposition maps}
\[ Q(A \cup_a B) \to Q(A) \barwedge Q(B) \]
and a \emph{counit map}
\[ Q(1) \to S \]
satisfying axioms dual to (1)--(4) of Definition \ref{def:operad}. A
\emph{morphism of cooperads} is a morphism of symmetric sequences that
commutes with the cocomposition and counit maps. A \emph{coaugmentation}
for a cooperad is a map $S \to Q(1)$ left inverse to the counit map. A
cooperad $Q$ is \emph{reduced} if the counit map is an isomorphism.
\end{definition}

\begin{remark} \label{rem:cooperads}
The description of an operad as a monoid for the composition product
of symmetric sequences naturally dualizes to cooperads. We define the
\emph{dual composition product} $\varcirc$ of two symmetric sequences
by replacing the coproduct in Definition \ref{def:compprod} with a
product. That is:
\[ M \varcirc N (A) := \prod_{A = \coprod_{j \in J} A_j} M(J) \barwedge
\Barwedge_{j \in J} N(A_j). \]
If $\barwedge$ commutes with finite products (which is in general not
likely) this is a monoidal product of symmetric sequences (the result dual
to Proposition \ref{prop:comp_monoidal}) and a cooperad is precisely a
comonoid for this product. In general we can define the triple product
$(L \varcirc M \varcirc N)$ by replacing coproduct with product in the
definition given in the proof of Lemma \ref{lem:circ_assoc}. We then
have natural maps
\[ (L \varcirc M) \varcirc N \to (L \varcirc M \varcirc N)
\quad\mbox{and}\quad
L \varcirc (M \varcirc N) \to (L \varcirc M \varcirc N) \]
which allow us to say what we mean by a comonoid in general. Thus we get
the result dual to Proposition \ref{prop:comp_general}, that a cooperad
in $\cat{C}$ is a symmetric sequence $Q$ together with maps
\[ Q \to Q \varcirc Q \quad\mbox{and}\quad Q \to I \]
such that the corresponding diagrams commute. In particular we have a
coassociativity diagram:
\[ \begin{diagram} \node[2]{Q \varcirc Q} \arrow{e} \node{(Q \varcirc Q)
\varcirc Q} \arrow{se} \\
            \node{Q} \arrow{ne} \arrow{se} \node[3]{(Q \varcirc Q \varcirc
            Q)} \\
            \node[2]{Q \varcirc Q} \arrow{e} \node{Q \varcirc (Q \varcirc
            Q)} \arrow{ne}
\end{diagram} \]
\end{remark}

\begin{remark}
In \cite{getzler/jones:1994} Getzler and Jones define a cooperad to be
a comonoid for the composition product $\circ$. In their case, $\circ$
and $\varcirc$ are equal because finite products are isomorphic to finite
coproducts in the category of chain complexes.
\end{remark}

\begin{definition}[(Comodules over a cooperad)] \label{def:comodule}
A \emph{left comodule} $C$ over the cooperad $Q$ is a left module over
$Q$ considered as an operad in $\cat{C}^\text{op}$. More explicitly, $C$ is
a symmetric sequence together with a left coaction of the comonoid $Q$,
that is, a map $C \to Q \varcirc C$. Equivalently, we have a suitable
collection of \emph{cocomposition maps}
\[ C(A) \to Q(r) \barwedge C(A_1) \barwedge \dots \barwedge C(A_r) \]
for partitions $A = \coprod_{i=1}^{r} A_i$. Similarly a \emph{right
comodule} is a symmetric sequence $C$ with a right coaction $C \to C
\varcirc Q$, or equivalently, \emph{cocomposition maps}
\[ C(A) \to C(r) \barwedge Q(A_1) \barwedge \dots \barwedge Q(A_r). \]
A \emph{bicomodule} is a symmetric sequence with compatible left and right
comodule structures. The cooperad $Q$ is itself a $(Q,Q)$--bicomodule.

A \emph{coalgebra} over a cooperad is the dual concept of an algebra over
an operad and the constant symmetric sequence with value equal to a
$Q$--coalgebra is a left $Q$--comodule.
\end{definition}

\section{Spaces of trees} \label{sec:trees}
As mentioned in the introduction to the paper, the key to finding a
cooperad structure on the bar construction on an operad is its
reinterpretation in terms of trees. These are the same sorts of trees used
in many other places to work with operads. See Getzler--Jones
\cite{getzler/jones:1994}, Ginzburg--Kapranov
\cite{ginzburg/kapranov:1994} and Markl--Shnider--Stasheff
\cite{markl/shnider/stasheff:2002} for many examples.

\begin{figure}[ht!]
\begin{center}
\input{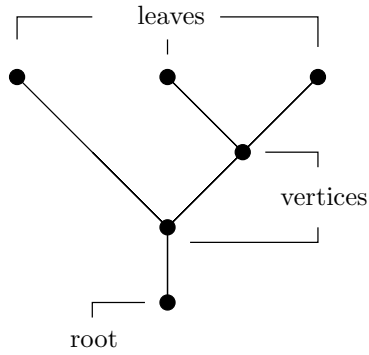}
\caption{Terminology for trees}
\label{fig:trees}
\end{center}
\end{figure}

\begin{definition}[(Trees)] \label{def:tree}
A typical tree of the sort we want is shown in Figure \ref{fig:trees}. It
has a root element at the base, a single edge attached to the root, and
no other vertices with only one incoming edge. We encode these geometric
requirements in the following combinatorial definition. A \emph{tree}
$T$ is a finite poset satisfying the following conditions:
\begin{enumerate}
\item $T$ has at least two elements: an initial (or minimal) element $r$,
the \emph{root}, and another element $b$ such that $b \leq t$ for all
$t \in T$, $t \neq r$.
\item For any elements $t,u,v \in T$, if $u \leq t$ and $v \leq t$,
then either $u \leq v$ or $v \leq u$.
\item For any $t < u$ in $T$ with $t \neq r$, there is some $v \in T$
such that $t < v$ but $u \nleq v$.
\end{enumerate}
We picture a tree by its \emph{graph}, whose vertices are the elements
of $T$ with an edge between $t$ and $u$ if $t < u$ and there is no $v$
with $t < v < u$. An \emph{incoming} edge to a vertex $t$ is an edge
corresponding to some relation $t < u$. Condition (1) above ensures that
the tree has a root $r$ with exactly one incoming edge (that connects
it to $b$). The second condition ensures that this graph is indeed a
tree in the usual sense. The third condition ensures that no vertices
except the root have exactly one incoming edge.

More terminology: the maximal elements of the tree $T$ will be called
\emph{leaves}. From now on, by a \emph{vertex}, we mean an element
other than the root or a leaf (see Figure \ref{fig:trees}). A tree
is \emph{binary} if each vertex has precisely two incoming edges. The
\emph{root edge} is the edge connected to the root element. The \emph{leaf
edges} are the edges connected to the leaves. The other edges in the tree
are \emph{internal edges}. Given a vertex $v$ of a tree, we write $i(v)$
for the set of incoming edges of the vertex $v$. We generally denote
trees with the letters $T,U,\dots$.
\end{definition}

\begin{remark} \label{rem:reduced_trees}
We stress that our trees are \emph{not} allowed to have vertices
with only one incoming edge, as guaranteed by condition (3) of
the definition. This reflects the fact that we will deal only with
\emph{reduced} operads in this paper.
\end{remark}

\begin{definition}[(Labellings)] \label{def:labelling}
A \emph{labelling} of the tree $T$ by a finite set $A$ is a bijection
between $A$ and the set of leaves of $T$. An \emph{isomorphism} of
$A$--labelled trees is an isomorphism of the underlying trees that
preserves the labelling. We denote the set of isomorphism classes
of $A$--labelled trees by $\mathsf{T}(A)$. For a finite set $A$,
$\mathsf{T}(A)$ is also finite. For a positive integer $n$, we write
$\mathsf{T}(n)$ for the set $\mathsf{T}(\{1,\dots,n\})$.
\end{definition}

\begin{figure}[ht!]
\begin{center}
\input{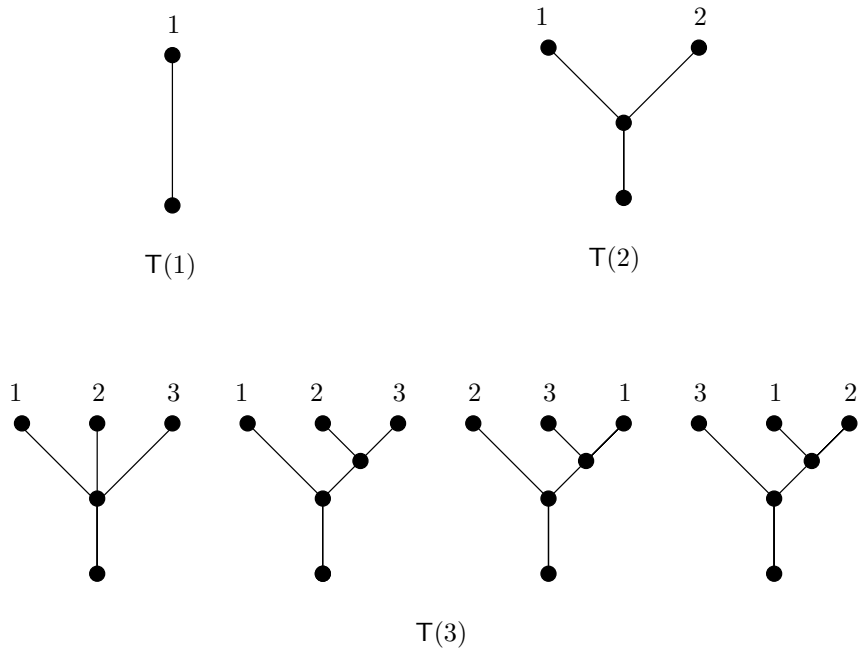}
\caption{Labelled trees with three or fewer leaves}
\label{fig:tree_examples}
\end{center}
\end{figure}

\begin{example}
There is up to isomorphism only one tree with one leaf. It has a single
edge whose endpoints are the root and the leaf. Thus $\mathsf{T}(1)$
has one element. It is easy to see that $\mathsf{T}(2)$ also only has
one element: the tree with one vertex that has two input edges. Figure
\ref{fig:tree_examples} shows $\mathsf{T}(1),\mathsf{T}(2),\mathsf{T}(3)$.
\end{example}

\begin{figure}[ht!]
\begin{center}
\input{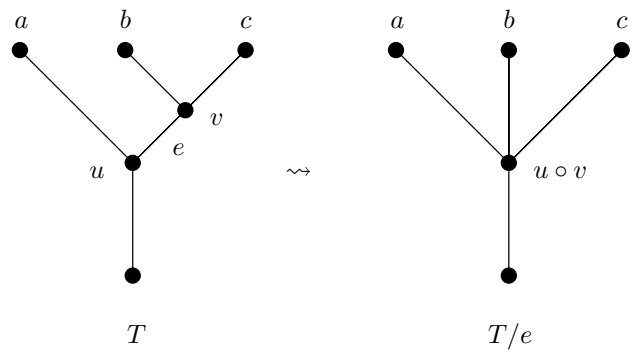}
\caption{Edge collapse of labelled trees}
\label{fig:edge_collapse}
\end{center}
\end{figure}

\begin{definition}[(Edge collapse)]
Given a tree $T$ and an internal edge $e$, denote by $T/e$ the tree
obtained by collapsing the edge $e$, identifying its endpoints. (In
poset terms, this is equivalent to removing from the poset the element
corresponding to the upper endpoint of the edge.) If $u$ and $v$ are those
endpoints, write $u \circ v$ for the resulting vertex of $T/e$. Note
that $T/e$ has the same leaves as $T$ so retains any labelling. See
Figure \ref{fig:edge_collapse} for an example.
\end{definition}

\begin{definition}
The process of collapsing edges gives us a partial order on the set
$\mathsf{T}(A)$ of isomorphism classes of $A$--labelled trees. We say
that $T \leq T'$ if $T$ can be obtained from $T'$ be collapsing a sequence
of edges. We think of the resulting poset as a category.
\end{definition}

We now give our trees topological significance by introducing `weightings'
on them.

\begin{definition} \label{def:weighting}
A \emph{weighting} on a tree $T$ is an assignment of nonnegative
`lengths' to the edges of $T$ in such a way that the `distance' from
the root to each leaf is exactly $1$. The set of weightings on a tree
$T$ is a subset of the space of functions from the set of edges of $T$
to the unit interval $[0,1]$ and we give it the subspace topology. We
denote the resulting space by $w(T)$. A tree together with a weighting
is a \emph{weighted tree}.
\end{definition}

\begin{example}
There is only one way to weight the unique tree $T \in \mathsf{T}(1)$ (the
single edge must have length $1$), so $w(T) = \basept$. For any $n$,
$\mathsf{T}(n)$ contains a tree $T_n$ with a single vertex that has $n$
incoming edges. For this tree we have $w(T_n) = \Delta^1$ the topological
$1$--simplex or unit interval. Figure \ref{fig:tree_examples} displays
another shape of tree with three leaves, one that has two vertices. For
such a tree $U$, we have $w(U) = \Delta^2$, the topological
$2$--simplex. Not all spaces of weightings are simplices, but we do have
the following result.
\end{example}

\begin{lemma} \label{lem:W}
Let $T$ be a tree with $n$ (internal) vertices. Then $w(T)$ is
homeomorphic to the $n$--dimensional disc $D^n$. If $n \geq 1$, the
boundary $\partial w(T)$ is the subspace of weightings for which at
least one edge has length zero.
\end{lemma}
\begin{proof}
Suppose $T$ has $l$ leaves. Then it has $n+l$ total edges and using the
lengths of the edges as coordinates we can think of $w(T)$ as a subset of
$\mathbb{R}^{n+l}$. For each leaf $l_i$ of $T$ there is a condition on
the lengths of the edges in a weighting that translates into an affine
hyperplane $H_i$ in $\mathbb{R}^{n+l}$. Then $w(T)$ is the intersection
of all these hyperplanes with $[0,1]^{n+l}$.

Now these hyperplanes all pass through the point that corresponds to the
root edge having length $1$ and all other edges length zero. Therefore
their intersection is another affine subspace of $\mathbb{R}^{n+l}$. To
see that they intersect transversely, we check that each $H_i$ does not
contain the intersection of the $H_j$ for $j \neq i$. Consider the point
$p_i$ in $\mathbb{R}^{n+l}$ that assigns length $1$ to each leaf edge
except that corresponding to leaf $l_i$, and length $0$ to all other edges
(including the leaf edge for $l_i$). Since the equation for the hyperplane
$H_j$ contains the length of exactly one leaf edge, this point $p_i$ is in
\[ \bigcap_{j \neq i} H_j \]
but not in $H_i$. This shows that the $H_i$ do indeed intersect
transversely. Therefore their intersection is an $n$--dimensional affine
subspace $V$ of $\mathbb{R}^{n+l}$.

Finally, notice that, as long as $n > 0$, $V$ passes through an interior
point of $[0,1]^{n+l}$, for example, the point where all edges except
the leaf edges have length $\varepsilon$ for some small $\varepsilon >
0$ and the leaf edges then have whatever lengths they must have to obtain
a weighting.
It then follows that $w(T) = V \cap [0,1]^{n+l}$ is homeomorphic to
$D^n$. If $n = 0$, there is only one tree and its space of weightings
is a single point, that is, $D^0$.

For the second statement, notice that the boundary of $w(T)$ is the
intersection of $V$ with the boundary of the cube $[0,1]^{n+l}$. If
a weighting includes an edge of length zero, it lies in this
boundary. Conversely, a weighting in this boundary must have some edge
with length either $0$ or $1$. If the root edge has length $1$, all
other edges must have length $0$. If some other edge has length $1$,
the root edge must have length $0$. In any case, some edge has length $0$.
\end{proof}

\begin{definition} \label{def:w_functor}
For each finite set $A$, the assignment $T \mapsto w(T)$ determines
a functor
\[ w(-)\co \mathsf{T}(A) \to \spaces \]
where $\spaces$ is the category of unbased spaces. To see this we must
define maps
\[ w(T/e) \to w(T) \]
whenever $e$ is an internal edge in the $A$--labelled tree $T$. Given a
weighting on $T/e$ we define a weighting on $T$ by giving edges in $T$
their lengths in $T/e$ with the edge $e$ having length zero. This is an
embedding of $w(T/e)$ as a `face' of the `simplex' $w(T)$. It's easy to
check that this defines a functor as claimed.

Let $w_0(T)$ be the subspace of $w(T)$ containing weightings for which
either the root edge or some leaf edge has length zero. We set
\[ \overline{w}(T) := w(T)/w_0(T). \]
This is a based space with basepoint given by the point to which $w_0(T)$
has been identified. If $T$ is the tree with only one edge then $w_0(T)$
is empty. We use the convention that taking the quotient by the empty
set is equivalent to adjoining a disjoint basepoint. So in this case,
$\overline{w}(T) = S^0$.

The maps $w(T/e) \to w(T)$ clearly map $w_0(T/e)$ to $w_0(T)$ and so
give us maps
\[ \overline{w}(T/e) \to \overline{w}(T). \]
For each finite set $A$, these form a functor
\[ \overline{w}(-) \co \mathsf{T}(A) \to \based \]
where $\based$ is the category of based spaces.
\end{definition}

\begin{figure}[ht!]
\begin{center}

\input{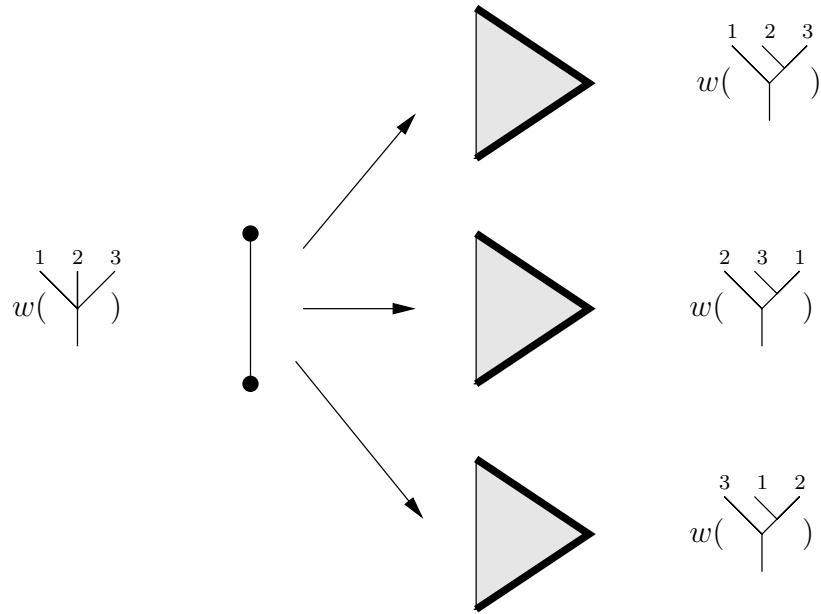}

\caption{Spaces of weightings of trees with three leaves}
\label{fig:w_functor}
\end{center}
\end{figure}

\begin{example} \label{ex:weightings}
Figure \ref{fig:w_functor} displays the spaces $w(T)$ for $T \in
\mathsf{T}(3)$ and how the functor $w(-)$ fits them together. Recall
that the poset $\mathsf{T}(3)$ has four objects: one minimal object
(the tree with one vertex and three incoming edges) and three maximal
objects (three binary trees with two vertices). As the picture shows, the
functor $w(-)$ embeds a 1--simplex for the minimal object as one of the
1--dimensional faces of a 2--simplex for each of the maximal objects. The
subspaces $w_0(T)$ are outlined in bold. Collapsing these we get the
functor $\overline{w}(-)$ which embeds $S^1$ (for the minimal object) as the
boundary of $D^2$ (for each maximal object).
\end{example}

\section{Bar constructions for reduced operads} \label{sec:bar}

This section forms the heart of the paper. We show that by giving an
explicit description of the simplicial bar construction in terms of trees,
we can construct a cooperad structure on it. In Section~\ref{sec:bardef}
we give our definition of the bar construction $B(P)$ for an operad $P$ in
$\cat{C}$. In Section~\ref{sec:simpbar} we show that this is isomorphic to
the standard simplicial reduced bar construction on $P$. Then in
Section~\ref{sec:cooperad} we prove the main result of this paper:
that $B(P)$ admits a natural cooperad structure.

We will work in a fixed symmetric monoidal $\based$--category
$(\cat{C},\barwedge,S)$ where $\based$ is the category of based
compactly-generated spaces and basepoint preserving maps. Since $\based$
is pointed, Proposition \ref{prop:null} implies that $\cat{C}$ too is
pointed. We denote the null object in $\cat{C}$ also by $\basept$. We
assume that $\cat{C}$ has all limits and colimits. The examples to bear in
mind are $\cat{C} = \based$ itself and $\cat{C} = \spectra$, which we take
to be the category of $S$--modules of EKMM
\cite{elmendorf/kriz/mandell/may:1997}, although other categories of
spectra could be used. We will use the notation developed in Section~1
for the enrichment, tensoring and cotensoring of $\cat{C}$ over $\based$.

Before we start we should stress that the constructions in this
paper only apply to \emph{reduced} operads and cooperads. That is,
those for which the unit (or counit) map is an isomorphism. This is
reflected in several places, most notably in the fact that our trees
are not allowed to have vertices with only one incoming edge (see Remark
\ref{rem:reduced_trees}). It is a necessary condition for our construction
of the cooperad structure on $B(P)$.

\subsection{Definition of the bar construction} \label{sec:bardef}

We give two definitions of the bar construction for an operad. The first
is somewhat informal and relies on $\cat{C}$ being the category of based
spaces, but captures how we really think about these objects. The second
is a precise formal definition as a coend in the category $\cat{C}$.

\begin{definition} \label{def:bar(operad)}
Let $P$ be a reduced operad in $\based$. The \emph{bar construction}
on $P$ is the symmetric sequence $B(P)$ defined as follows. A general
point $p$ in $B(P)(A)$ consists of
\begin{itemize}
\item an isomorphism class of $A$--labelled trees: $T \in \mathsf{T}(A)$,
\item a weighting on $T$ and,
\item for each (internal) vertex $v$ of $T$, a point $p_v$ in the based
space $P(i(v))$ (recall that $i(v)$ is the set of incoming edges of the
vertex $v$),
\end{itemize}
subject to the following identifications:
\begin{itemize}
\item If $p_v$ is the basepoint in $P(i(v))$ for any $v$ then $p$ is
identified with the basepoint $\basept \in B(P)(A)$.
\item If the internal edge $e$ has length zero, we identify $p$ with
the point $q$ given by
\begin{itemize}
\item the tree $T/e$,
\item the weighting on $T/e$ in which an edge has the same length
as the corresponding edge of $T$ in the weighting that makes up
$p$,\footnote{This is the inverse image under the injective map
\[ w(T/e) \to w(T) \]
of the weighting corresponding to $p$. The condition that $e$ has
length zero says precisely that the weighting for $p$ is in the image
of this map.}
\item $q_{u \circ v}$ given by the image under the composition map
\[ P(i(u)) \smsh P(i(v)) \to P(i(u \circ v)) \]
of $(p_u,p_v)$ (notice that $i(u \circ v) = i(u) \circ_v i(v)$),
\item $q_t = p_t$ for the other vertices $t$ of $T/e$.
\end{itemize}
\item If a root or leaf edge has length zero, $p$ is identified with
$\basept \in B(P)(A)$.
\end{itemize}
A bijection $\sigma\co A \to A'$ gives us an isomorphism $\sigma_*\co B(P)(A)
\to B(P)(A')$ by relabelling the leaves of the underlying trees. In this
way, $B(P)$ becomes a symmetric sequence in $\based$.
\end{definition}

\begin{example}
Consider $B(P)(1)$. There is only one tree with a single leaf and only
one weighting on it. It has no vertices so $B(P)(1)$ does not depend at
all on $P$. With the basepoint (which is disjoint in this case because
nothing is identified to it) we get $B(P)(1) = S^0$.

Next consider $B(P)(2)$. Again there is only one tree, but this time it
has a vertex (with two incoming edges) and the space of ways to weight the
tree is the $1$--simplex $\Delta^1$. Making all the identifications we
see that
\[ B(P)(2) = \Sigma P(2), \]
the reduced suspension of $P(2)$.
\end{example}

\begin{definition}[(The functors $P_A$)] \label{def:P}
A key ingredient of the general definition of the bar construction is
that an operad $P$ in $\cat{C}$ determines a functor
\[ P_A(-)\co \mathsf{T}(A)^\text{op} \to \cat{C}. \]
where $\mathsf{T}(A)$, as always, is the poset of isomorphism classes of
$A$--labelled trees ordered by edge collapse. For a tree $T$ we define
\[ P_A(T) := \Barwedge_{\text{vertices $v$ in $T$}} P(i(v)) \]
where we recall that $i(v)$ is the set of incoming edges to the vertex
$v$. If $e$ is an internal edge in $T$ with endpoints $u$ and $v$ then
there is a partial composition map
\[ P(i(u)) \barwedge P(i(v)) \to P(i(u \circ v)). \]
Using this we get a map
\[ P_A(T) \to P_A(T/e). \]
The associativity axioms for the operad $P$ ensure that these maps make
$P_A(-)$ into a functor as claimed.
\end{definition}

Recall from Definition \ref{def:w_functor} that we have a functor
\[ \overline{w}(-) \co \mathsf{T}(A) \to \based \]
given by taking the space of weightings on a tree, modulo those for
which a root or leaf edge has length zero.

\begin{definition}[(Formal definition of the bar construction)]
\label{def:formal_bar}
Let the \emph{bar construction} of the reduced operad $P$ be the symmetric
sequence $B(P)$ defined by
\[ B(P)(A) := \overline{w}(-) \otimes_{\mathsf{T}(A)} P_A(-) = \int^{T \in
\mathsf{T}(A)} \overline{w}(T) \otimes P_A(T). \]
This is the coend in $\cat{C}$ of the bifunctor
\[ \overline{w}(-) \otimes P_A(-) \co \mathsf{T}(A) \times
\mathsf{T}(A)^\text{op}
\to \cat{C}. \]
(See \cite{maclane:1971} for the theory of coends.) The definition of
the coend is a colimit over a category whose objects are morphisms in
$\mathsf{T}(A)$ and we will write the coend above as
\[ \colim_{T \leq T' \in \mathsf{T}(A)} \overline{w}(T) \otimes P_A(T') \]
when we need to manipulate it as such.

A bijection $A \to A'$ induces an isomorphism of categories $\mathsf{T}(A)
\to \mathsf{T}(A')$ by the relabelling of trees. If $T \mapsto T'$
under this isomorphism then $P_A(T) = P_A(T')$ and $\overline{w}(T) =
\overline{w}(T')$. Therefore we get an induced isomorphism $B(P)(A) \to
B(P)(A')$. This makes $B(P)$ into a symmetric sequence in $\cat{C}$.
\end{definition}

\begin{remark}
To see that our two definitions of the bar construction are equivalent
when $\cat{C} = \based$, recall that the coend is a quotient of the
coproduct
\[ \bigwdge_{T \in \mathsf{T}(A)} \overline{w}(T) \otimes P_A(T). \]
That is, a point consists of a weighted tree together with elements
of the $P(i(v))$ for vertices $v$ subject to some identifications. The
maps $P_A(T) \to P_A(T/e)$ and $\overline{w}(T/e) \to \overline{w}(T)$ encode the
identifications made in Definition \ref{def:bar(operad)}.
\end{remark}

\begin{remark}
Our definition of the bar construction is rather reminiscent of the
geometric realization of simplicial sets or spaces. This line of thought
leads to the definition of an \emph{arboreal object} in $\cat{C}$ as
a functor
\[ \mathsf{T}(A)^\text{op} \to \cat{C} \]
in which $\mathsf{T}(A)$ plays the role of the simplicial indexing
category $\Delta$ for simplicial sets. With the spaces of weightings
$\overline{w}(T)$ playing the role of the topological simplices, the bar
construction $B(P)$ can be thought of as the geometric realization of
the arboreal object $P_A(-)$. We will formalize and extend these ideas
in future work \cite{ching:2005b}.
\end{remark}

\begin{remark} \label{rem:BV}
The $W$--construction of Boardman and Vogt (also sometimes called the
bar construction) is defined in a very similar manner to $B(P)$. It uses
slightly different spaces of trees and produces an operad instead of a
cooperad. See \cite{vogt:2003} for details. Benoit Fresse has noticed
a relationship between $W(P)$ and $B(P)$, namely that
\[ B(P) = \Sigma\operatorname{Indec}(W(P)) \]
where $\Sigma$ is a single suspension (that is, tensoring with $S^1$)
and $\operatorname{Indec}$ denotes the `operadic indecomposables functor'. It is the cooperad structure on $\Sigma\operatorname{Indec}(W(P))$, corresponding to that on $B(P)$, that was described by Salvatore in \cite{salvatore:1998}.
\end{remark}

\begin{example} \label{ex:bars}
Let $\mathcal{A}ss$ be the operad for associative monoids in unbased
spaces. This is given by
\[ \mathcal{A}ss(n) := \Sigma_n \]
(with the discrete topology and regular $\Sigma_n$--action). The
composition maps are the inclusions given by identifying
\[ \Sigma_r \times \Sigma_{n_1} \times \dots \times \Sigma_{n_r} \]
with a subgroup of $\Sigma_{n_1 + \dots + n_r}$. We obtain an operad
$\mathcal{A}ss_+$ in $\based$ by adding a disjoint basepoint to each of
the terms of $\mathcal{A}ss$. Let us calculate $B(\mathcal{A}ss_+)$.

The points $p_v \in \mathcal{A}ss_+(i(v))$ required by Definition
\ref{def:bar(operad)} can be thought of as determining an order on
the incoming edges to vertices of a tree. This allows us to identify
a point in $B(\mathcal{A}ss_+)(n)$ with a \emph{planar} weighted tree
with leaves labelled $1,\dots,n$. This breaks $B(\mathcal{A}ss_+)(n)$
up into a wedge of $n!$ terms, each corresponding to an ordering of the
leaves of the trees involved.

As we now show, each of these terms is an $(n-1)$--sphere. Think of
constructing a planar weighted tree with leaves labelled in a fixed order
(say, $1,\dots,n$) by the following method. Connect the first leaf to
the root with an edge of length $1$. Then attach the second leaf at
some point along the edge already drawn. Attach the third leaf at some
point along the path from the second leaf to the root, and so on. The
space of choices made in doing all this is $[0,1]^{n-1}$ and we obtain
precisely the planar weighted trees we want in this manner (see Figure
\ref{fig:planar}). The root edge or a leaf edge will have length zero if
and only if at least one of our choices was either $0$ or $1$. Hence the
space we want is obtained by identifying the boundary of $[0,1]^{n-1}$
to a basepoint. This gives $S^{n-1}$.

\begin{figure}[ht!]
\begin{center}
\input{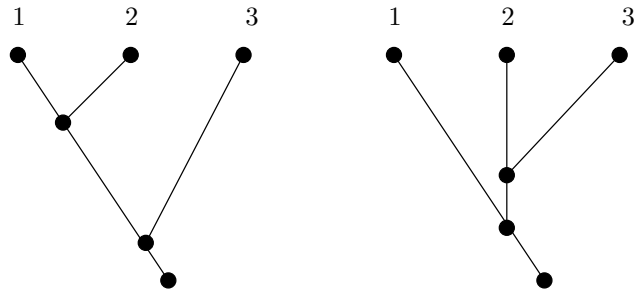}
\caption{Constructing planar weighted trees}
\label{fig:planar}
\end{center}
\end{figure}

Therefore we have
\[ B(\mathcal{A}ss_+)(n) \isom S^{n-1} \smsh (\Sigma_n)_+ \]
where $\Sigma_n$ acts trivially on the $S^{n-1}$ term and by translation
on the non-basepoints of $(\Sigma_n)_+$.

We can also picture what happens for $n = 3$ in terms of sticking
together the spaces $\overline{w}(T) \smsh \mathcal{A}ss_3(T)_+$ for $T \in
\mathsf{T}(3)$. The $\overline{w}(T)$ are the quotients of the spaces pictured
in Figure \ref{fig:w_functor} by the subspaces outlined in bold. To
make up $B(\mathcal{A}ss_+)(3)$ we need six copies of the 1--simplex
(corresponding to the points in $\mathcal{A}ss(3)$) and twelve copies
of the 2--simplex. (There are four points in $\mathcal{A}ss(2) \times
\mathcal{A}ss(2)$ and three trees of this type.) These fit together to
form six disjoint copies of the space of Figure \ref{fig:ass}, one for
each permutation of $1,2,3$. The type of tree used to form each part
is shown.
\begin{figure}[ht!]
\begin{center}
\input{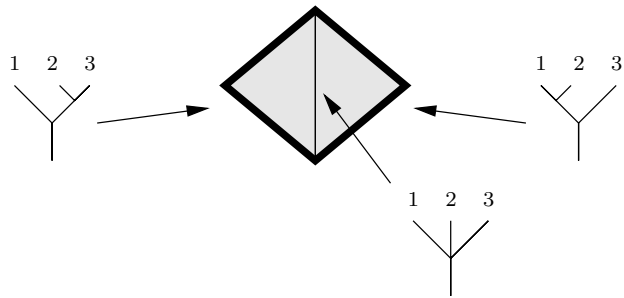}
\caption{One sixth of $B(\mathcal{A}ss_+)(3)$}
\label{fig:ass}
\end{center}
\end{figure}
When we collapse the bold subspaces to the basepoint we get a wedge of
six copies of $S^2$ as expected.
\end{example}

\subsection{Relation to the simplicial bar construction}
\label{sec:simpbar}

In this section we show that $B(P)$ is isomorphic to the geometric
realization of the standard simplicial bar construction on
the reduced operad $P$. This simplicial bar construction can be
defined for any augmented monoid in a monoidal category.\footnote{See
\cite[Section~II.2.3]{markl/shnider/stasheff:2002} for a discussion of different
forms of the simplicial bar construction.} We have seen (Proposition
\ref{prop:comp_monoidal}) that under the right conditions an operad is
just a monoid for the monoidal product on the category of symmetric
sequences given by the composition product $\circ$. To define the
simplicial bar construction in general (that is, without the assumption
that $\barwedge$ commutes with finite coproducts) we must say what
we mean by higher iterates of $\circ$. For this we use the following
natural extension of the three-way product introduced in the proof of
Lemma \ref{lem:circ_assoc}.

\begin{definition}[(Iterated composition product)] \label{def:iterated}
The \emph{composition product} of the symmetric sequences $M_1,\dots,M_r$
is the symmetric sequence given by
\[ (M_1 \circ \dots \circ M_r)(A) :=\!\! \bigvee_{A_i = \coprod_{a \in
A_{i-1}} A_{i,a}}\!\! M_1(A_1) \barwedge \Barwedge_{a \in A_1}
M_2(A_{2,a}) \barwedge \dots \barwedge \Barwedge_{a \in A_{r-1}}
M_r(A_{r,a}) \]
for each finite set $A = A_r$. Here we are taking the coproduct over
partitions of $A$ into subsets $A_{r,a}$ indexed over $a \in A_{r-1}$,
partitions of $A_{r-1}$ indexed over $A_{r-2}$, and so on. Equivalently
we can view this coproduct as indexed over sequences of $r-1$ partitions
of $A$, each a refinement of the next.
\end{definition}

\begin{remark}
There is a natural map from $(M_1 \circ \dots \circ M_r)$ to any of
the symmetric sequences obtained by choosing ways to bracket this
expression. All the `obvious' diagrams relating these maps commute. If
$\barwedge$ commutes with finite coproducts in $\cat{C}$ then all these
maps are isomorphisms and reflect the associativity isomorphisms of the
monoidal product $\circ$.
\end{remark}

\begin{definition}[(Simplicial bar construction)]
\label{def:simplicial_bar}
Let $P$ be a reduced operad in $\cat{C}$. The \emph{simplicial bar
construction} $\cat{B}_{\bullet}(P)$ is the simplicial object in the
category of symmetric sequences on $\cat{C}$ with
\[ \cat{B}_k(P) = \underset{k}{\underbrace{P \circ \dots \circ P}}. \]
For $i = 1,\dots,k-1,$ face maps
\[ d_i\co \underset{k}{\underbrace{P \circ \dots \circ P}} \to
\underset{k-1}{\underbrace{P \circ \dots \circ P}} \]
are given by
\[ \dots \circ P \circ P \circ \dots \to \dots \circ (P \circ P) \circ
\dots \to \dots \circ P \circ \dots \]
where we are using the operad composition $P \circ P \to P$ to compose
the \ord{$i$} and \ord{$i{+}1$} factors. The maps $d_0$ and $d_k$ are given
by applying the augmentation map $P \to I$ to the first and last copies
of $P$ respectively. Degeneracy maps
\[ s_j\co \underset{k}{\underbrace{P \circ \dots \circ P}} \to
\underset{k+1}{\underbrace{P \circ \dots \circ P}} \]
are given for $j = 0,\dots, k$ by using the unit map $I \to P$ to insert
a copy of $P$ between the \ord{$j$} and \ord{$(j+1)$} factors:
\[ \cdots \circ P \circ P \circ \cdots \isom \cdots \circ P \circ I
\circ P \circ \cdots \to \cdots \circ P \circ P \circ P \circ \cdots. \]
\end{definition}

\begin{remark}
It is sufficient for this definition that $P$ be augmented. However,
we need $P$ to be reduced to make the following identification of the
simplicial bar construction with $B(P)$ as defined previously.
\end{remark}

\begin{prop} \label{prop:simp}
Let $P$ be a reduced operad in $\cat{C}$. Then the geometric
realization\footnote{The geometric realization of a simplicial symmetric
sequence is defined pointwise: $|X|(A) = |X(A)|$. Note that a simplicial
symmetric sequence is the same thing as a symmetric sequence of simplicial
objects.} of $\cat{B}_{\bullet}(P)$ is isomorphic to the bar construction
$B(P)$.
\end{prop}
\begin{proof}
We give the proof for $\cat{C} = \based$ (which is the only case we
require in this paper) based on the informal description of $B(P)$ in
Definition \ref{def:bar(operad)}. The same idea could be used to write a
proof that works for any $\cat{C}$ using the formal definition of $B(P)$
as a coend.

The idea is that the iterated composition products that make up the
simplicial bar construction can be thought of in terms of sequences of
partitions which in turn are related to trees of the type we are using
to define $B(P)$.

We first give an explicit description of the $n$--simplices in
$\cat{B}_{\bullet}(P)(A)$. These are given by the object
\[ \underset{n}{\underbrace{P \circ \dots \circ P}}(A). \]
Enlarging on the last sentence of Definition \ref{def:iterated}, we can
write this as a coproduct over all sequences of partitions
\[ \widehat{0} = \lambda_0 \leq \lambda_1 \leq \dots \leq \lambda_{n-1}
\leq \lambda_n = \widehat{1} \]
of the set $A$,
where $\lambda \leq \mu$ if $\lambda$ is \emph{finer} than $\mu$ (if
two elements of $A$ are in the same block in $\lambda$, they are also
in the same block in $\mu$) and $\widehat{0},\widehat{1}$ are the minimal and
maximal partitions with respect to this order. The terms in the coproduct
are appropriate smash products of the $P(r)$. We get a factor of $P(r)$
every time one of the blocks of one of the partitions breaks up into $r$
blocks in the next partition along.

A point in the geometric realization $|\cat{B}_{\bullet}(P)|$ can be
represented by a point in the topological $n$--simplex $\Delta^n$ together
with a choice of sequence of partitions as described above and a point
in the appropriate smash product of the spaces $P(r)$.

A sequence of partitions determines an $A$--labelled tree $T$
as follows. (See Figure \ref{fig:simplicial} for an example when
$n = 3$.) Take a vertex for each block of each $\lambda_i$ for $i =
1,\dots,n$. Add a root, and a leaf for each element of $A$. Two vertices
are joined by an edge if they come from consecutive partitions of
the sequence and the block for one is contained in the block for the
other. Finally we add a root edge from the $\lambda_n$ vertex to the
root and a leaf edge from each leaf to the corresponding $\lambda_1$
vertex. (Notice that vertices in this tree might have only one input
edge -- let's allow this for the moment.)
\begin{figure}[ht!]
\begin{center}
\input{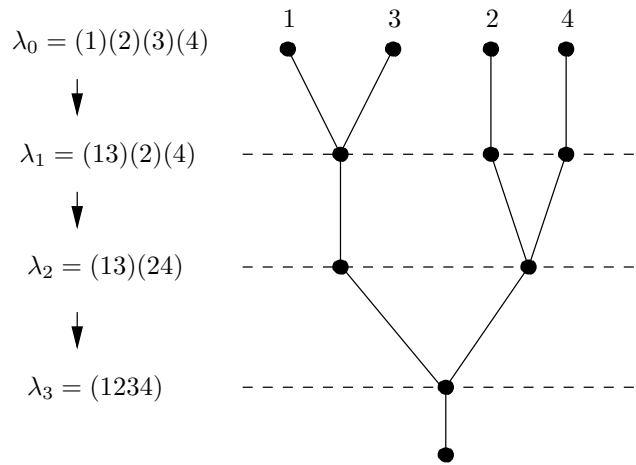}
\caption{Producing trees from sequences of partitions}
\label{fig:simplicial}
\end{center}
\end{figure}

A point in $\Delta^n$ determines a weighting on the tree we have just
constructed. Thinking of $\Delta^n$ as the subspace of $\mathbb{R}^{n+1}$
with $x_0 + \dots + x_n = 1$ and $x_i \geq 0$, we get a weighting by
giving the root edge length $x_0$, the edges connecting the vertices for
$\lambda_{i-1}$ to the vertices for $\lambda_i$ length $x_i$ and the leaf
edges length $x_n$. We can now remove the vertices with only one input
edge, connecting their input and output edges. This gives us a point in
$w(T)$ for some tree $T$ in the sense of Definition \ref{def:tree}.

Finally notice that because $P(1) = S^0$ (as $P$ is reduced), the smash
product of spaces $P(r)$ determined by the sequence of partitions is
precisely $P_A(T)$. Therefore we actually obtain a point in $B(P)(A)$.

It remains to show that this process sets up a homeomorphism between the
geometric realization $|\cat{B}_{\bullet}(P)(A)|$ and $B(P)(A)$. There
are a couple of key steps. Firstly the degeneracy maps in the simplicial
bar construction are isomorphisms on terms in the coproduct. These
correspond to inserting lots of vertices with one input edge in our
trees, which are then removed by our construction. So we only have to
worry about the identifications made by the face maps. The face maps are
given by removing partitions from the sequences, which corresponds to
edge collapse. Hence the identifications made in defining $B(P)$ are the
same as those in defining the realization of $\cat{B}_{\bullet}(P)$. This
completes the proof.
\end{proof}

\subsection{Cooperad structure on the bar construction}
\label{sec:cooperad}

Up to this point, all we have done is identify the simplicial bar
construction on a reduced operad in terms of trees. The main point of
this paper is that this identification allows us to see that there is a
cooperad structure on the bar construction. In this section we describe
that structure. The key to getting the cooperad cocomposition maps is
the process of grafting (or rather \emph{un}grafting) trees.

\begin{definition}[(Tree grafting)] \label{def:grafting}
Let $T$ be an $A$--labelled tree, let $U$ be a $B$--labelled tree and let
$a$ be an
element of $A$. We define the \emph{grafting of $U$ onto $T$ at $a$}
to be the tree $T \cup_a U$ obtained by identifying the root edge of $U$
to the leaf edge of $T$ corresponding to $a$. Figure \ref{fig:grafting}
below illustrates this process.

\begin{figure}[ht!]
\begin{center}

\input{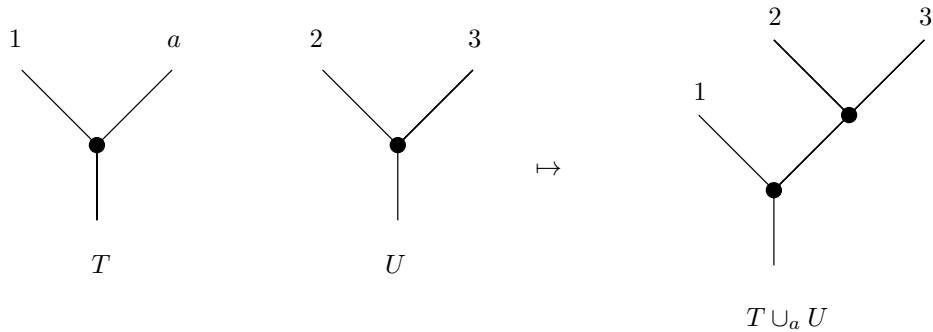}

\caption{Tree grafting}
\label{fig:grafting}
\end{center}
\end{figure}

We denote the newly identified edge by $e_a$. Every other edge of $T
\cup_a U$ comes either from $T$ or from $U$. The vertices of $T \cup_a U$
are the vertices of $T$ together with the vertices of $U$ (and they have
the same number of incoming edges). Finally there is a natural $A \cup_a
B$--labelling of $T \cup_a U$, given by combining the labellings of $T$
and of $U$.

We say that an $A \cup_a B$--labelled tree is \emph{of type $(A,B)$} if it
is of the form $T \cup_a U$ for an $A$--labelled tree $T$ and a
$B$--labelled tree $U$. The next lemma says that an $A \cup_a B$--labelled
tree is a grafting in at most one way. This is trivial but crucial to
the construction of the cooperad structure maps below.
\end{definition}

\begin{lemma} \label{lem:ungrafting}
For any $A \cup_a B$--labelled tree $V$ there is at most one pair $(T,U)$
such that $V = T \cup_a U$.
\end{lemma}
\begin{proof}
In the grafted tree $T \cup_a U$ the `upper' endpoint of the edge $e_a$
is a vertex whose `parent leaves' are labelled precisely by the elements
of $B$. There can be at most one such vertex $v$ in $V$ and cutting along
the edge immediately below $v$ produces the trees $T,U$ that make up $V$.
\end{proof}

\begin{definition} \label{def:cooperad_maps}
To give $B(P)$ a cooperad structure we have to define maps
\begin{equation} B(P)(A \cup_a B) \to B(P)(A) \barwedge B(P)(B)
\label{eq:cooperad_maps} \end{equation}
for finite sets $A,B$ and $a \in A$. A point $p$ in $B(P)(A \cup_a B)$
consists of a weighted tree $V$ labelled by $A \cup_a B$ together with
elements of $p_v \in P(i(v))$ for vertices $v$ of $V$. We treat two cases:
\begin{enumerate}
\item If $V$ is not of the form $T \cup_a U$ for an $A$--labelled tree $T$
and a $B$--labelled tree $U$, then we will map $p$ to the basepoint on
the right-hand side of (\ref{eq:cooperad_maps}).
\item If $V$ is of this form (that is, it is of type $(A,B)$)
then things are more interesting. Below we describe how the map
(\ref{eq:cooperad_maps}) is defined in this case.
\end{enumerate}
Since $V$ is of type $(A,B)$, Lemma \ref{lem:ungrafting} tells us that
there is a unique $A$--labelled tree $T$ and a unique $B$--labelled
tree $U$ such that $V = T \cup_a U$. We use these trees as the basis for
elements $q \in B(P)(A)$ and $r \in B(P)(B)$ respectively. What remains
to be seen is how the weighting and vertex labels of $V$ determine
weightings and vertex labels for $T$ and $U$.

The vertex labels are easy because the vertices of $T \cup_a U$ consist
of the vertices of each of $T$ and $U$ with the same numbers of input
edges. Therefore we take
\[ q_v := p_v \in P(i(v)) \]
for vertices $v$ of $T$ and
\[ r_u := p_u \in P(i(v)) \]
for vertices $u$ of $U$.

The way in which a weighting on $T \cup_a U$ determines weightings on $T$
and $U$ is the key part of our construction. This comes about via a map
\begin{equation} \overline{w}(T \cup_a U) \to \overline{w}(T) \smsh
\overline{w}(U)
\label{eq:key} \end{equation}
(recall that $\overline{w}(-)$ is the space of weightings on a tree with
those that have zero length root or leaf edges identified to a basepoint).

So take a weighting of $T \cup_a U$. Define a weighting on $T$ by giving
the edges the same lengths they had in $T \cup_a U$ and giving the leaf
edge for $a$ the necessary length to make the root-leaf distances equal to
$1$.\footnote{Intuitively, we have collapsed the $U$ part of the tree to a
single edge with the same overall length.} Next define a weighting on $U$
by taking the lengths from $T \cup_a U$ and scaling up by a constant
factor to make the root-leaf distances equal to $1$ (the length of the
root edge of $U$ comes from the length of the edge $e_a$ in $T \cup_a U$).
The scaling factor is the inverse of the total length of the $U$ part of
$T \cup_a U$. The only time this doesn't work is if all the $U$--edges in
$T \cup_a U$ (including $e_a$) are of length zero. However in that case
the weighting we just defined on $T$ has a leaf edge of length zero and
so is the basepoint in $\overline{w}(T)$. This is almost enough to define a
map of the form (\ref{eq:key}). The only thing left to check is that if a
leaf or root edge of $T \cup_a U$ is of length zero then the same is true
of either of the chosen weightings on $T$ and $U$. This is clear. Figure
\ref{fig:key} illustrates a particular case of the map (\ref{eq:key}).

\begin{figure}[ht!]
\begin{center}
\input{key.pstex_t}
\caption{The map $\overline{w}(T \cup_a U) \to \overline{w}(T) \smsh
\overline{w}(U)$}
\label{fig:key}
\end{center}
\end{figure}

This completes the definition of the cooperad structure maps
(\ref{eq:cooperad_maps}):
\[ B(P)(A \cup_a B) \to B(P)(A) \barwedge B(P)(B) \]
given, in summary, by:
\[ p = (V,\{p_v\}) \mapsto \begin{cases} q = (T,\{p_v\}_{v \in T}), \;
r = (U,\{p_v\}_{v \in U}) & \text{if $V = T \cup_a U$}; \\ \basept &
\text{otherwise}. \end{cases} \]
with the weightings on $T,U$ given by the map (\ref{eq:key}) just
constructed.

We still have to check that these maps are well-defined. To see this we
have to look at the identifications made in the definition of $B(P)(A
\cup_a B)$:
\begin{itemize}
\item If $p_v$ equals the basepoint in $P(i(v))$ for any vertex $v \in
V$ then the same will be true of the corresponding vertex in either $T$
or $U$. Hence such a $p$ maps to the basepoint.
\item If an interior edge $e$ of the tree $V$ underlying the point $p$
is of length zero, $p$ is identified with another point $p'$ as described
in Definition \ref{def:bar(operad)}. We have various possibilities:
\begin{enumerate}
\item $V$ is not of the form $T \cup_a U$ in which case neither is $V/e$
and both $p$ and $p'$ map to the basepoint.
\item $V = T \cup_a U$ and $e$ corresponds to an internal edge of $T$. In
this case, the points $q$ and $q'$ will be identified via the collapse
of that edge, and the points $r$ and $r'$ will be equal. So $p$ and $p'$
map to the same element of $B(P)(A) \barwedge B(P)(B)$.
\item $V = T \cup_a U$ and $e$ corresponds to an internal edge of
$U$. This is similar to case (2).
\item $V = T \cup_a U$ and $e$ is the edge $e_a$ obtained from identifying
the root edge of $U$ with the $a$--leaf edge of $T$. In this case
$V/e$ is no longer of the form $T \cup_a U$ and so $p'$ maps to the
basepoint. But in the weighting on $U$ determined by that on $T \cup_a
U$ the root edge has length scaled up from the length of $e_a$ which is
therefore zero. So the point $r$ is the basepoint in $B(P)(B)$ and so $p$
also maps to the basepoint.
\end{enumerate}
\item We have already checked in the definition of the map (\ref{eq:key})
that if a root or leaf edge in $p$ is of length zero, then the same is
true of at least one of $q$ and $r$. Therefore such a $p$ maps to the
basepoint in $B(P)(A) \barwedge B(P)(B)$.
\end{itemize}
This completes the check that our maps (\ref{eq:cooperad_maps}) are
well-defined. The final piece of the cooperad structure for $B(P)$ is a
counit map $B(P)(1) \to S^0$. But we already saw that $B(P)(1) \isom S^0$
(in the based space case) so our counit is this isomorphism. Note that
this means $B(P)$ turns out to be a reduced cooperad.
\end{definition}

\begin{example}
The map
\[ B(P)(\{1,2,3\}) \to B(P)(\{a,3\}) \smsh B(P)(\{1,2\}) \]
is pictured in Figure \ref{fig:key}. The left-hand side (with
vertices labelled by elements of $P(2)$) represents a point $p$ of
$B(P)(\{1,2,3\})$. The two trees on the right-hand side (with vertices
labelled by those same elements in the obvious way) represent the image
of $p$ in $B(P)(\{a,3\}) \barwedge B(P)(\{1,2\})$. In this example,
all points that are based on trees of shapes other than that shown are
mapped to the basepoint.
\end{example}

We will save for later the task of checking that these maps do
indeed give us a cooperad structure. First we translate Definition
\ref{def:cooperad_maps} into the category-theoretic language needed to
define the cocomposition maps for a general $\cat{C}$. To do this, we
notice that the `ungrafting' process more-or-less makes our categories
$\mathsf{T}(A)$ into a cooperad of categories. To make this precise,
we describe an `add a disjoint basepoint' functor for categories.

\begin{definition}[(Categories with initial objects)]
Write $\mathsf{Cat}_+$ for the category in which an object is a (small)
category $\mathsf{C}_+$ together with an initial object $\basept$
such that $\Hom_{\mathsf{C}_+}(X,\basept)$ is empty for all $X \neq
\basept$. The morphisms in $\mathsf{Cat}_+$ are functors that preserve
the initial objects.

There is a functor from the category $\mathsf{Cat}$ of all (small)
categories to $\mathsf{Cat}_+$ given by adding an initial object
with the correct morphisms to a category $\mathsf{C}$ to obtain
$\mathsf{C}_+$. Note that every object in $\mathsf{Cat}_+$ can be
obtained in this way, but not every morphism in $\mathsf{Cat}_+$ is
given by adding an initial object to a morphism in $\mathsf{Cat}$.

Define a symmetric monoidal product $\smsh$ on $\mathsf{Cat}_+$ by
\[ \mathsf{C}_+ \smsh \mathsf{D}_+ := \mathsf{C}_+ \times
\mathsf{D}_+/\mathsf{C}_+ \wdge \mathsf{D}_+, \]
where the wedge product is the disjoint union with the initial objects
identified and the quotient identifies this wedge product to the initial
object of the smash product. Notice that if $\mathsf{C},\mathsf{D}
\in \mathsf{Cat}$ then
\[ \mathsf{C}_+ \smsh \mathsf{D}_+ = (\mathsf{C} \times \mathsf{D})_+. \]
The unit for this product is the category $1_+$ with two objects and a
single morphism between them.
\end{definition}

In particular we write $\mathsf{T}(A)_+$ for the category formed by adding
an initial object to our poset of $A$--labelled trees $\mathsf{T}(A)$. The
reason for making all these new definitions is then the following result.

\begin{prop} \label{prop:category_cooperad}
The categories $\mathsf{T}(A)_+$ form a reduced cooperad in
$\mathsf{Cat}_+$.
\end{prop}
\begin{proof}
The cocomposition maps have the form
\[ \mathsf{T}(A \cup_a B)_+ \to \mathsf{T}(A)_+ \smsh \mathsf{T}(B)_+ =
(\mathsf{T}(A) \times \mathsf{T}(B))_+ \]
and are given by `ungrafting' trees. Take $V \in \mathsf{T}(A \cup_a
B)$. If $V$ is a tree of type $(A,B)$ we map it to the pair $(T,U)$
where $T,U$ are the unique trees that graft together to give $V$ (see
Lemma \ref{lem:ungrafting}). If $V$ is not of type $(A,B)$ (or is the
initial object) we map it to the initial object of the right-hand side.

First we must check that we have indeed given a functor here. Suppose
that $V \leq V'$ in $\mathsf{T}(A \cup_A B)$. The only interesting
case is when $V$ is of type $(A,B)$, so maps to a pair $(T,U)$ on the
right-hand side. We have to show two things: that $V'$ is also of type
$(A,B)$ with decomposition $(T',U')$ and then that $T \leq T'$ and $U
\leq U'$. Well, let $e_a$ be the edge in $V$ at which the grafting took
place. Since $V$ is obtained from $V'$ by a sequence of edge collapses,
$e_a$ must come from an edge $e_{a'}$ in $V'$ that is not collapsed in
this sequence. This edge breaks $V'$ into two parts and we can write $V'
= T' \cup_{a'} U'$ for some trees $T',U'$ with some labellings (a priori,
not necessarily by $A$ and $B$). But it is now clear that $U'$ must yield
$U$ after undergoing some edge collapses. So $U' \in \mathsf{T}(B)$
and $U \leq U'$. Similarly, $T' \in \mathsf{T}(A)$ and $T \leq T'$
(after relabelling $a'$ by $a$).

Notice that $\mathsf{T}(1)_+$ is isomorphic to the unit $1_+$ for the
symmetric monoidal structure on $\mathsf{Cat}_+$. We take as unit map the
(unique) isomorphism $1_+ \to \mathsf{T}(1)_+$.

It still remains to check that the cooperad axioms do indeed hold for
our cocomposition maps. This is simple and we leave it to the reader.
\end{proof}

\begin{remark}
The original categories $\mathsf{T}(A)$ in fact already form an
\emph{operad} in $\mathsf{Cat}$ with composition maps given by grafting
rather than ungrafting. This operad structure is effectively what is used
by Boardman and Vogt to define their $W$--construction.
\end{remark}

The next step is to show that the bar construction can be defined as a
coend in $\mathsf{T}(A)_+$ instead of $\mathsf{T}(A)$.

\begin{lemma}
Let $P$ be a reduced operad in $\cat{C}$. The functors $\overline{w}(-)$
and $P_A(-)$ on $\mathsf{T}(A)$ naturally extend to functors
\[ \overline{w}(-)\co \mathsf{T}(A)_+ \to \based \]
and
\[ P_A(-)\co (\mathsf{T}(A)_+)^\text{op} \to \cat{C} \]
and we have
\[ B(P)(A) = \int^{T \in \mathsf{T}(A)_+} \overline{w}(T) \otimes P_A(T). \]
\end{lemma}
\begin{proof}
We set $\overline{w}(\basept) = \basept_{\based}$ and $P_A(\basept) =
\basept_{\cat{C}}$ with the necessary definition on morphisms (given
by the fact that $\basept_{\based}$ is an initial object in $\based$
and $\basept_{\cat{C}}$ is a terminal object in $\cat{C}$). It is then
clear that $\basept \in \mathsf{T}(A)_+$ does not contribute anything to
the coend which therefore reduces to the previous definition of $B(P)(A)$.
\end{proof}

The maps (\ref{eq:key}) of Definition \ref{def:cooperad_maps} are still
the key ingredients in constructing the cooperad maps for $B(P)$.

\begin{lemma} \label{lem:trans1}
The maps
\[ \overline{w}(T \cup_a U) \to \overline{w}(T) \smsh \overline{w}(U) \]
previously defined form part of a natural transformation
\[ \begin{diagram} \node{\mathsf{T}(A \cup_a B)_+} \arrow[2]{s}
\arrow{ese,t}{\overline{w}(-)} \\
    \node[2]{\Downarrow} \node{\based} \\
    \node{\mathsf{T}(A)_+ \smsh \mathsf{T}(B)_+} \arrow{ene,b}{\overline{w}(-)
    \smsh \overline{w}(-)} \\
\end{diagram}. \]
\end{lemma}
\begin{proof}
The bottom functor here is defined in the obvious way on $\mathsf{T}(A)
\times \mathsf{T}(B)$ and sends $\basept$ to $\basept$. For $V \in
\mathsf{T}(A \cup_a B)$ not of type $(A,B)$, the corresponding part of
the natural transformation is
\[ \overline{w}(V) \to \basept. \]

The only really interesting naturality square comes from $V \leq V'$
with $V'$ of type $(A,B)$ and $V$ not. The square that must commute in
this case is
\[ \begin{diagram} \node{\overline{w}(V)} \arrow{s} \arrow{e}
\node{\overline{w}(V')} \arrow{s} \\
            \node{\basept} \arrow{e} \node{\overline{w}(T') \smsh
\overline{w}(U').}
\end{diagram} \]
This is the content of part (4) of the checking we did towards the end of
Definition \ref{def:cooperad_maps}: from any weighting on $V$, the
weighting we get on $V'$ will have length zero for the edge connecting the
$T'$--part to the $U'$--part. Hence the root edge of the corresponding
weighting on $U'$ will have length zero. So we map into the basepoint
of $\overline{w}(T') \smsh \overline{w}(U')$.
\end{proof}

We have a corresponding result for the functors $P_A(-)$ of Definition
\ref{def:formal_bar}.

\begin{lemma} \label{lem:trans2}
Let $P$ be a reduced operad in $\cat{C}$. Then there is a natural
transformation
\[ \begin{diagram} \node{\mathsf{T}(A \cup_a B)_+^\text{op}} \arrow[2]{s}
\arrow{ese,t}{P_{A \cup_a B}(-)} \\
    \node[2]{\Downarrow} \node{\cat{C}} \\
    \node{(\mathsf{T}(A)_+ \smsh \mathsf{T}(B)_+)^\text{op}}
    \arrow{ene,b}{P_A(-) \barwedge P_B(-)}
\end{diagram} \]
\end{lemma}
\begin{proof}
In other words, given $V \in \mathsf{T}(A \cup_a B)$ we have maps
\[ P_{A \cup_a B}(V) \to P_A(T) \barwedge P_B(U) \]
when $V = T \cup_a U$. There are obvious isomorphisms that we take for
these maps. The naturality squares are easily seen to commute. Again
the only one that seems like it might be interesting is for $V \leq V'$
with $V'$ of type $(A,B)$ and $V$ not. But in fact this square just
turns out to be
\[ \begin{diagram} \node{P_{A \cup_a B}(V')} \arrow{s} \arrow{e}
\node{P_{A \cup_a B}(V)} \arrow{s} \\
        \node{P_A(T') \barwedge P_B(U')} \arrow{e} \node{\basept}
\end{diagram} \]
which is not so interesting after all.
\end{proof}

Finally, we can give the formal construction of the cocomposition maps
for the cooperad $B(P)$.

\begin{definition} \label{def:formal_cooperad_maps}
Let $P$ be a reduced operad in $\cat{C}$ and let $B(P)$ be the symmetric
sequence of Definition \ref{def:formal_bar}. The cocomposition map
\[ B(P)(A \cup_a B) \to B(P)(A) \barwedge B(P)(B) \]
is given by the following sequence of maps:
\begin{equation} \begin{split}
B(P)&(A \cup_a B)  = \colim_{V \leq V' \in \mathsf{T}(A \cup_a B)_+}
\overline{w}(V) \otimes P_{A \cup_a B}(V') \\
  \longrightarrow & \colim_{(T,U) \leq (T',U') \in \mathsf{T}(A)_+ \smsh
  \mathsf{T}(B)_+} (\overline{w}(T) \smsh \overline{w}(U)) \otimes (P_A(T')
  \barwedge P_B(U')) \\
  \longrightarrow & \colim_{(T,U) \leq (T',U') \in \mathsf{T}(A)_+ \smsh
  \mathsf{T}(B)_+} (\overline{w}(T) \otimes P_A(T')) \barwedge
  (\overline{w}(U) \otimes P_B(U')) \\
  \longrightarrow & \left(\colim_{T \leq T' \in \mathsf{T}(A)_+} \overline{w}(T)
  \otimes P_A(T') \right) \barwedge \left(\colim_{U \leq U' \in
  \mathsf{T}(B)_+} \overline{w}(U) \otimes P_B(U') \right) \\
  = & \; B(P)(A) \barwedge B(P)(B). \\
  \end{split} \label{eq:cocomp} \end{equation}
The first map here comes from combining the natural transformations of
Lemmas \ref{lem:trans1} and \ref{lem:trans2}. The second is given by the
transformation $d$ of Definition \ref{def:axiom}. It is for precisely this
reason that the axiom giving us $d$ is necessary. The third map is given
by universal properties of colimits. This completes the construction of
the cooperad structure maps for $B(P)$.
\end{definition}

The next task is to check that the maps we have described actually do
make $B(P)$ into a cooperad. That is, we must check the duals of axioms
(1)--(4) from Definition \ref{def:operad}. The key step is to see that
the maps (\ref{eq:key}) satisfy corresponding conditions.

\begin{lemma} \label{lem:keyassoc}
Let $T,U,V$ be $A$--, $B$-- and $C$--labelled trees respectively and let
$a,a' \in A$ and $b \in B$. Let $I$ denote the unique $\basept$--labelled
tree. Recall that $\overline{w}(I) = S^0$. Then the following diagrams commute:
\begin{enumerate}
\item
\( \qquad\qquad\begin{diagram}
    \node{\overline{w}(T \cup_a U \cup_b V)} \arrow{s} \arrow{e}
    \node{\overline{w}(T \cup_a U) \smsh \overline{w}(V)} \arrow{s} \\
    \node{\overline{w}(T) \smsh \overline{w}(U \cup_b V)} \arrow{e}
\node{\overline{w}(T)
    \smsh \overline{w}(U) \smsh \overline{w}(V)}
\end{diagram} \)
\item
\( \qquad\qquad\begin{diagram}
    \node{\overline{w}(T \cup_a U \cup_{a'} V)} \arrow{s} \arrow{e}
    \node{\overline{w}(T \cup_a U) \smsh \overline{w}(V)} \arrow{s} \\
    \node{\overline{w}(T \cup_{a'} V) \smsh \overline{w}(U)} \arrow{e}
    \node{\overline{w}(T) \smsh \overline{w}(U) \smsh \overline{w}(V)}
\end{diagram} \)
\item
\( \qquad\qquad\qquad\begin{diagram}
    \node{\overline{w}(T \cup_a I)} \arrow{e} \arrow{se,b}{\isom}
    \node{\overline{w}(T) \smsh \overline{w}(I)} \arrow{s} \\
    \node[2]{\overline{w}(T)}
\end{diagram} \)
\item
\( \qquad\qquad\qquad\begin{diagram}
    \node{\overline{w}(I \cup_\basept T)} \arrow{e} \arrow{se,b}{\isom}
    \node{\overline{w}(I) \smsh \overline{w}(T)} \arrow{s} \\
    \node[2]{\overline{w}(T)}
\end{diagram} \)
\end{enumerate}
\end{lemma}
\begin{proof}
The argument for diagram (1) is contained in Figure \ref{fig:assoc1}. A
point in $\overline{w}(T \cup_a U \cup_b V)$ comes from a weighting of the
grafted tree $T \cup_a U \cup_b V$. The top-left corner of Figure
\ref{fig:assoc1} shows a generic version of such a tree with some
lengths labelled:
\begin{itemize}
\item $u$ is the length of the root edge.
\item $v$ is the distance from the root vertex to the lower vertex of
the edge that joins $U$ to $T$ (there may be intermediate vertices along
this route, we let $v$ denote the total distance).
\item $w$ is the length of the edge that joins $U$ to $T$.
\item $x$ is the distance from the upper vertex of that edge to the
lower vertex of the edge that joins $V$ to $U$.
\item $y$ is the length of the edge that joins $V$ to $U$.
\item $z$ is the remaining distance to any of the leaves of $V$.
\end{itemize}
Figure \ref{fig:assoc1} shows that whichever way we map our weighted
tree around diagram (1) we get the same result. (Note that if $y+z$ or
$w+x+y+z$ equals to zero, then $z = 0$ and we are the basepoint in every
corner of diagram (1).) We therefore conclude that diagram (1) commutes.

\begin{sidewaysfigure}[p]
\input{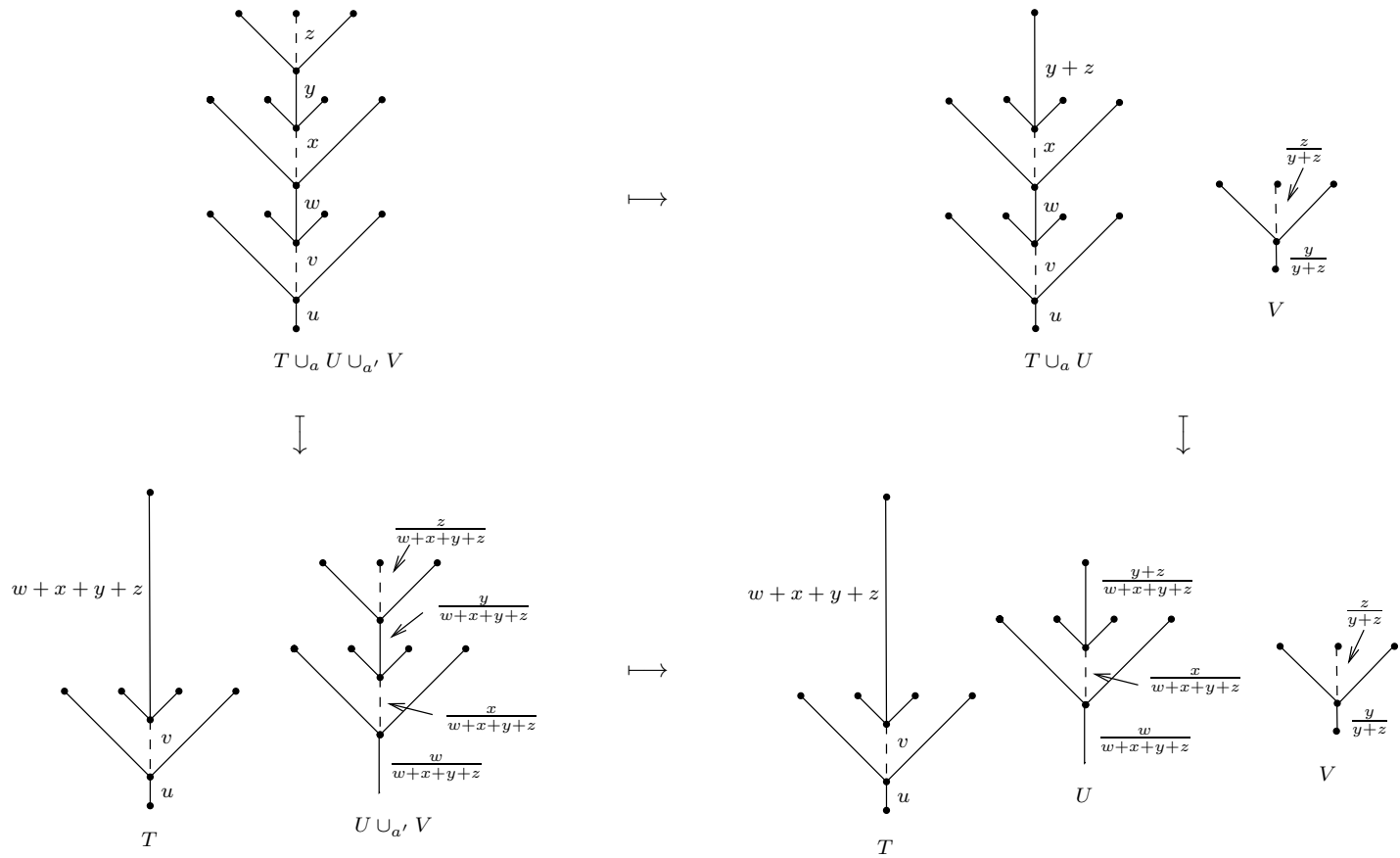}
\caption{Commutativity of diagram (1) of Lemma \ref{lem:keyassoc}. The
dashed lines represent sequences of possibly more than one edge.}
\label{fig:assoc1}
\end{sidewaysfigure}

Diagram (2) is similar to (1) but easier. For diagram (3), notice
that the image in $\overline{w}(T)$ of a weighting of $T \cup_a I$ will be
effectively the same weighting. The image in $\overline{w}(I) = S^0$ will
be the non-basepoint unless the leaf edge for $a$ has length zero. But
if this is the case our starting point was the basepoint in $\overline{w}(T
\cup_a I)$. This shows that the diagram commutes.

For diagram (4), the image in $\overline{w}(T)$ of a weighting of $I
\cup_\basept T$ will again be the very same weighting (no scaling
up is necessary). The image in $\overline{w}(I) = S^0$ will always be the
non-basepoint. Therefore this diagram also commutes.
\end{proof}

We are now in a position to state the main result of this paper.

\begin{thm} \label{thm:bar=cooperad}
Let $P$ be a reduced operad in the symmetric monoidal $\based$--category
$\cat{C}$. The maps of Definition \ref{def:cooperad_maps} make $B(P)$
into a reduced cooperad in $\cat{C}$.
\end{thm}
\begin{proof}
We give the formal argument for the maps of Definition
\ref{def:formal_cooperad_maps}. To fit the relevant diagrams onto a page
we need some new notation. Let's write
\[ \overline{w}(T,U) := \overline{w}(T) \smsh \overline{w}(U) \]
and
\[ P_{A,B}(T',U') := P_A(T') \barwedge P_B(U'). \]
Figure \ref{fig:assoc2} then shows the diagram that has to commute for
the dual of axiom (1) of Definition \ref{def:operad} to hold for $B(P)$.
\begin{sidewaysfigure}[p] \small
\[ \begin{diagram} \dgARROWLENGTH=.7em
    \node{\colim_{S \leq S' \in \mathsf{T}(A \cup_a B \cup_b C)_+}
    \overline{w}(S) \otimes P_{A \cup_a B \cup_b C}(S')} \arrow{s} \arrow{e}
    \node{\colim_{\substack{Q \leq Q' \in \mathsf{T}(A \cup_a B)_+ \\ V
    \leq V' \in \mathsf{T}(C)_+}} \begin{matrix} \overline{w}(Q,V) \otimes \\
    P_{A \cup_a B,C}(Q',V') \end{matrix}} \arrow{s} \arrow{e}
    \node{\colim_{\substack{Q \leq Q' \in \mathsf{T}(A \cup_a B)_+ \\
    V \leq V' \in \mathsf{T}(C)_+}} \begin{matrix} (\overline{w}(Q) \otimes
    P_{A \cup_a B}(Q')) \barwedge \\ (\overline{w}(V) \otimes P_C(V'))
    \end{matrix}} \arrow{s} \\
    \node{\colim_{\substack{T \leq T' \in \mathsf{T}(A)_+ \\ R \leq R'
    \in \mathsf{T}(B \cup_b C)_+}} \begin{matrix} \overline{w}(T,R) \otimes \\
    P_{A,B \cup_b C}(T',R') \end{matrix}} \arrow{s} \arrow{e}
    \node{\colim_{\substack{T \leq T' \in \mathsf{T}(A)_+ \\ U \leq U'
    \in \mathsf{T}(B)_+ \\ V \leq V' \in \mathsf{T}(C)_+}} \begin{matrix}
    \overline{w}(T,U,V) \otimes \\ P_{A,B,C}(T',U',V') \end{matrix}} \arrow{e}
    \arrow{s}
    \node{\colim_{\substack{T \leq T' \in \mathsf{T}(A)_+ \\ U \leq U'
    \in \mathsf{T}(B)_+ \\ V \leq V' \in \mathsf{T}(C)_+}} \begin{matrix}
    (\overline{w}(T,U) \otimes P_{A,B}(T',U')) \barwedge \\ (\overline{w}(V)
    \otimes P_C(V')) \end{matrix}} \arrow{s} \\
    \node{\colim_{\substack{T \leq T' \in \mathsf{T}(A)_+ \\ R \leq R'
    \in \mathsf{T}(B \cup_b C)_+}} \begin{matrix} (\overline{w}(T) \otimes
    P_A(T')) \barwedge \\ (\overline{w}(R) \otimes P_{B \cup_b C}(R'))
    \end{matrix}} \arrow{e}
    \node{\colim_{\substack{T \leq T' \in \mathsf{T}(A)_+ \\ U \leq U'
    \in \mathsf{T}(B)_+ \\ V \leq V' \in \mathsf{T}(C)_+}} \begin{matrix}
    (\overline{w}(T) \otimes P_A(T')) \barwedge \\ (\overline{w}(U,V) \otimes
    P_{B,C}(U',V')) \end{matrix}} \arrow{e}
    \node{\colim_{\substack{T \leq T' \in \mathsf{T}(A)_+ \\ U \leq U'
    \in \mathsf{T}(B)_+ \\ V \leq V' \in \mathsf{T}(C)_+}} \begin{matrix}
    (\overline{w}(T) \otimes P_A(T')) \barwedge \\ (\overline{w}(U) \otimes P_B(U'))
    \barwedge \\ (\overline{w}(V) \otimes P_C(V')) \end{matrix}}
\end{diagram} \]
\caption{The commutative diagram that verifies the first associativity
axiom for the cooperad structure on $B(P)$ in the proof of Theorem
\ref{thm:bar=cooperad}}
\label{fig:assoc2}
\end{sidewaysfigure}
The key to showing that this commutes is putting
\[ \colim_{\substack{T \leq T' \in \mathsf{T}(A)_+ \\ U \leq U' \in
\mathsf{T}(B)_+ \\ V \leq V' \in \mathsf{T}(C)_+}} \overline{w}(T,U,V)
\otimes P_{A,B,C}(T',U',V') \]
into the center of the square. We've connected this to the top and left
sides of the square using maps similar to the first map in Definition
\ref{def:formal_cooperad_maps}. We've connected it to the right and bottom
sides using maps of the form $d$ from Definition \ref{def:axiom}. It's
then enough to show that the four smaller squares commute.

The top-left square commutes because of diagram (1) in Lemma
\ref{lem:keyassoc}. The bottom-left and top-right squares commute because
of the naturality of the transformations $d$. The bottom-right square
commutes because it is an example of the associativity axiom we required
of our $d$ transformations in \ref{def:axiom}.

This completes the verification of the dual of axiom (1) of Definition
\ref{def:operad}. For axiom (2) the argument is similar, but using diagram
(2) of Lemma \ref{lem:keyassoc}. For the duals of axioms (3) and (4)
we use the unit axiom for the transformations $d$ together with diagrams
(3) and (4) of Lemma \ref{lem:keyassoc}. We leave the reader to fill in
the details of these proofs.
\end{proof}

\section{Cobar constructions for reduced cooperads} \label{sec:cobar}

We now dualize to cooperads. The cobar construction for a cooperad is
strictly dual to the bar construction for an operad. More precisely,
recall that a cooperad $Q$ in a category $\cat{C}$ is the same thing
as an operad $Q^\text{op}$ in the opposite category $\cat{C}^\text{op}$. The
cobar construction on $Q$ is then defined to be the bar construction
on $Q^\text{op}$. This bar construction is a cooperad in $\cat{C}^\text{op}$ and
hence an operad in $\cat{C}$. In symbols, the \emph{cobar construction
on $Q$} is
\[ \Omega(Q) := B(Q^\text{op})^\text{op}. \]
It can be useful to have a more explicit description of this.

\begin{definition}[(Cobar construction on a cooperad)]
\label{def:cobar(cooperad)}
The cobar construction,
being dual to the bar construction,
is defined as an end rather than a coend. Let $Q$ be a cooperad in
$\cat{C}$. Then
for each finite set $A$, $Q$ determines a functor
\[ Q_A(-) \co \mathsf{T}(A) \to \cat{C} \]
by
\[ Q_A(T) = Q(i(v_1)) \barwedge \dots \barwedge Q(i(v_n)) \]
where $v_1,\dots,v_n$ are the vertices of $T$. This is a functor because
the cocomposition maps for $Q$ give us maps
\[ Q_A(T/e) \to Q_A(T). \]
(Recall that the corresponding functor for an operad was defined on
$\mathsf{T}(A)^\text{op}$.) The \emph{cobar construction} $\Omega(Q)$ is then
the symmetric sequence with
\[ \Omega(Q)(A) := \Map_{\mathsf{T}(A)}(\overline{w}(-),Q_A(-)) =  \int_{T
\in \mathsf{T}(A)} \Map_{\cat{C}}(\overline{w}(T),Q_A(T)). \]
This is the end of the bifunctor
\[ \mathsf{T}(A)^\text{op} \times \mathsf{T}(A) \to \cat{C} \]
given by
\[ (T,U) \mapsto \Map_{\cat{C}}(\overline{w}(T),Q_A(U)) \]
where $\Map_{\cat{C}}$ denotes the cotensoring structure for $\cat{C}$
over $\based$ (and hence the tensoring structure for $\cat{C}^\text{op}$).
\end{definition}

\begin{remark}
The cobar construction $\Omega(Q)$ on a reduced cooperad $Q$ in based
spaces is isomorphic to the totalization of a cosimplicial cobar
construction that is dual to the simplicial bar construction. The
terms in this cosimplicial construction are iterated versions of the
dual composition product of Remark \ref{rem:cooperads}. The fact that
$\Omega(Q)$ is the totalization of this is dual to the result that $B(P)$
is the realization of the simplicial bar construction.
\end{remark}

The operad structure maps for $\Omega(Q)$ are dual to the cooperad
maps for $B(P)$. The following result is the dual of Proposition
\ref{thm:bar=cooperad}.

\begin{cor} \label{cor:cobar=operad}
Let $Q$ be a reduced cooperad in a symmetric monoidal $\based$--category
$\cat{C}$. Then the cobar construction $\Omega(Q)$ is a reduced operad
in $\cat{C}$.
\qed
\end{cor}

\section{Duality for operads and cooperads} \label{sec:dual}

In this section we examine how the bar and cobar constructions relate
to the `duality' functor
\[ \dual \co \based^\text{op} \to \cat{C}; \; X \mapsto \Map_{\cat{C}}(X,S) \]
where $S$ is the unit of the symmetric monoidal structure on $\cat{C}$.
The case to keep in mind is $\cat{C} = \spectra$ in which case $S$ is the
sphere spectrum and this duality functor is Spanier--Whitehead duality.

\begin{lemma} \label{lem:dual}
Let $Q$ be a cooperad of based spaces. Then $\dual Q$ is an operad in
the category $\cat{C}$.
\end{lemma}
\begin{proof}
The composition maps for $\dual Q$ are given by
\[ \begin{split} \Map_{\cat{C}}(Q(A),S) \barwedge \Map_{\cat{C}}(Q(B),S)
&\to \Map_{\cat{C}}(Q(A) \smsh Q(B),S) \\ &\to \Map_{\cat{C}}(Q(A \cup_a
B),S). \\ \end{split} \]
The first map is the natural transformation constructed in Proposition
\ref{prop:dual} (it's the distributive map $d$ for $\cat{C}^\text{op}$). The
second comes from the corresponding cocomposition map for $Q$.
\end{proof}

\begin{remark}
The dual of an operad need not in general be a cooperad because the map
$d$ need not in general have an inverse. However when it does we have
a nice duality result connecting the bar and cobar constructions. For
this to work we need to put the following condition on the spaces that
make up our operad.
\end{remark}

\begin{definition}
Two based spaces $X,Y$ are \emph{compatibly dualizable in $\cat{C}$}
if the map
\[ d\co \Map_{\cat{C}}(X,S) \barwedge \Map_{\cat{C}}(Y,S) \to
\Map_{\cat{C}}(X \smsh Y,S) \]
is an isomorphism.
\end{definition}

\begin{prop} \label{prop:duality}
Let $P$ be an operad in based spaces whose terms (that is, the $P(A)$
for finite sets $A$) are pairwise compatibly dualizable. Then $\dual P$
has a natural cooperad structure. Moreover, we have an isomorphism
\[ \dual B(P) \isom \Omega(\dual P) \]
of operads in $\cat{C}$.
\end{prop}
\begin{proof}
The cooperad structure maps for $\dual P$ are constructed in the same
way as the operad structure maps for $\dual Q$ in \ref{lem:dual} but
using the inverse of the relevant map $d$ provided by the `compatibly
dualizable' hypothesis.

The second part relies on the descriptions of the bar and cobar
constructions as coends and ends respectively. The coend $B(P)$ is
a colimit:
\[ B(P)(A) = \colim_{T \leq T'} \overline{w}(T) \smsh P_A(T') \]
where the colimit is taken over all inequalities of trees in
$\mathsf{T}(A)$. Therefore
\[\begin{split}
    \dual B(P)(A) & = \Map_{\cat{C}}( \colim \overline{w}(T) \smsh P_A(T'),
    S) \\
    & \isom \lim \Map_{\cat{C}}(\overline{w}(T) \smsh P_A(T'), S) \\
    & \isom \lim \Map_{\cat{C}}(\overline{w}(T), \Map_{\cat{C}}(P_A(T'),
    S)) \\
    & \isom \lim \Map_{\cat{C}}(\overline{w}(T), (\dual P)_A(T')). \\
\end{split} \]
The last identity again uses the `compatibly dualizable' hypothesis in
the form
\[ \Map_{\cat{C}}(P(i(v_1)) \smsh \ldots \smsh P(i(v_n)), S)
\isom \Map_{\cat{C}}(P(i(v_1)), S) \barwedge \ldots \barwedge
\Map_{\cat{C}}(P(i(v_n)),S). \]
The final line of this calculation is precisely the limit that defines
$\Omega(\dual P)$. We leave the reader to check that this is an
isomorphism of operads.
\end{proof}

\begin{remark}
The only case of this result we will use in this paper is when all
the terms of the operad $P$ are $S^0$. These are pairwise compatibly
dualizable in any $\cat{C}$ because
\[ \Map_{\cat{C}}(S^0,C) \isom C \]
for any $C \in \cat{C}$.
\end{remark}

\begin{remark}
Replacing $\cat{C}$ with $\cat{C}^\text{op}$ we obtain dual results. These
concern the functor $\mathbb{S}\co X \mapsto X \otimes S$, the `suspension
spectrum' functor. We find that if $Q$ is a cooperad in based spaces
then $\mathbb{S}Q$ is a cooperad in $\cat{C}$. If $P$ is an operad whose
terms are pairwise compatibly dualizable then $\mathbb{S}P$ is an operad
in $\cat{C}$ and $\mathbb{S}B(P) \isom B(\mathbb{S}P)$.
\end{remark}

We have now reached the stage where we can apply our constructions to
Goodwillie's calculus of functors (see
Section~\ref{sec:application}). Before doing so, we extend our bar and
cobar constructions to modules and comodules. This will then allow us
to construct modules over the derivatives of the identity.

\section{Bar constructions for modules and comodules}
\label{sec:bar(modules)}

In this section we extend the bar and cobar constructions to modules and
comodules. We show that there is a bar construction on left (respectively
right) modules over a reduced operad $P$ that yields left (respectively
right) comodules over the cooperad $B(P)$. Dually, there is a cobar
construction on left (respectively right) comodules over a reduced
cooperad $Q$ that yields left (respectively right) modules over the
operad $\Omega(Q)$. These are special cases of two-sided bar and cobar
constructions. Given a reduced operad $P$ with right module $R$ and left
module $L$, we will define a two-sided bar construction $B(R,P,L)$. Taking
either $R$ or $L$ to be the unit symmetric sequence $I$ will yield the
promised one-sided constructions for individual modules. The two-sided
construction is isomorphic to the standard simplicial two-sided bar
construction (see Definition \ref{def:two-sided_simp}) but, in order to
get the comodule structure, we have reinterpreted this in terms of trees.

Most of the material in this section is a straightforward generalization
of that of Sections~\ref{sec:trees}--\ref{sec:cobar}. First, in
Section~\ref{sec:gen_trees} we describe the more general species of tree
necessary for the definitions of the two-sided constructions. In
Section~\ref{sec:two-sided} we give these definitions and show that the
bar construction of Section~\ref{sec:bardef} is a special case. In
Section~\ref{sec:bar(modules)_maps} we construct the maps that make
the bar construction on a module into a comodule, and dually, the cobar
construction on a comodule into a module.

As previously, $\cat{C}$ denotes a symmetric monoidal $\based$--category
with null object $\basept$ and which has all necessary limits and colimits.

\subsection{Generalized trees} \label{sec:gen_trees}
To accommodate the presence of the $P$--modules $R$ and $L$ in the
two-sided bar construction, we need to make two changes to our notion
of tree, one at the root level and one at the leaf level:
\begin{enumerate}
\item We allow the root element of a tree to have more than one incoming
edge.
\item We allow the leaves of a tree to have repeated labels, that is, an
$A$--labelling is a surjection from $A$ to the set of leaves, rather
than a bijection.
\end{enumerate}
We will refer to this notion as a `generalized tree', or sometimes just
a `tree' if the context makes it clear that we mean the generalized
version. The following definition makes things precise.

\begin{definition} \label{def:gen_trees}
Let $A$ be a finite set. A \emph{generalized $A$--labelled tree}
consists of
\begin{itemize}
\item a poset $T$ with a unique minimal element $r$ (the \emph{root})
satisfying conditions (2) and (3) of Definition \ref{def:tree}, and
\item a surjection $\iota$ from the finite set $A$ to the set of maximal
elements (the \emph{leaves}) of $T$.
\end{itemize}
We use letters $T,U,\dots$ to denote generalized trees, usually taking the
labelling map $\iota$ for granted. We write $\mathsf{Tree}(A)$ for the set
of isomorphism classes of generalized $A$--labelled trees. All the
terminology of Definition \ref{def:tree} applies equally well to
generalized trees.
\end{definition}

Edge collapse for generalized trees is defined in exactly the same way as
for the trees of Section~\ref{sec:trees} except that now we allow
ourselves to collapse root edges as well as internal edges. To get the
right category structure on $\mathsf{Tree}(A)$ we need a way to collapse
leaf edges as well. The following definition provides this.

\begin{figure}[ht!]
\begin{center}
\input{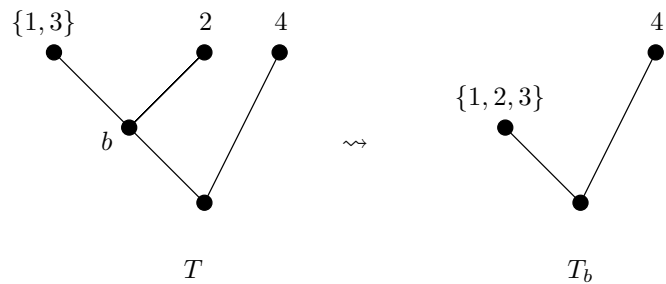}
\caption{An example of bud collapse for generalized
$\{1,2,3,4\}$--labelled trees}
\label{fig:bud_collapse}
\end{center}
\end{figure}

\begin{definition}[(Bud collapse)] \label{def:bud_collapse}
A \emph{bud} in a generalized tree $T$ is a vertex all of whose incoming
edges are leaf edges. Equivalently, a bud is a maximal vertex. If $b$ is a
bud in $T$, a \emph{$b$--leaf} is a leaf of $T$ that is attached to $b$.

Given a generalized $A$--labelled tree $T$ and a bud $b \in T$, we define
a generalized $A$--labelled tree $T_b$ which is obtained from $T$ by
\emph{bud collapse}. The underlying poset of $T_b$ is obtained from $T$ by
removing the $b$--leaves. This makes $b$ into a leaf in $T_b$. The
$A$--labelling on $T_b$ is that of $T$ for the leaves that still remain,
with $b$ inheriting the labels of its old leaves. Formally, we are
composing the $A$--labelling on $T$ with the surjection from the leaves of
$T$ to the leaves of $T_b$ that sends the $b$--leaves in $T$ to
$b$. Visually, we can think of this process as collapsing all the leaf
edges attached to $b$ (see Figure \ref{fig:bud_collapse}).
\end{definition}

\begin{figure}[ht!]
\begin{center}
\input{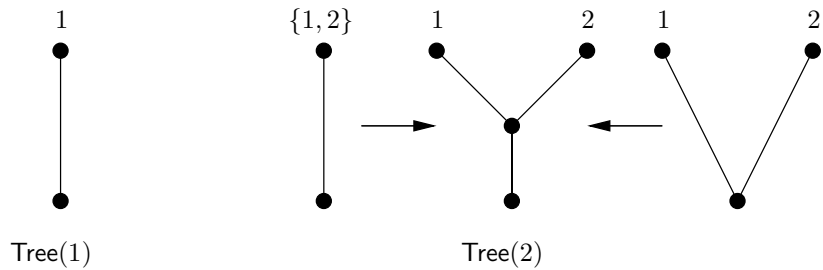}
\caption{$\mathsf{Tree}(1)$ and $\mathsf{Tree}(2)$ (the arrows represent
the direction of the morphisms in $\mathsf{Tree}(2)$)}
\label{fig:gen_trees}
\end{center}
\end{figure}

\begin{definition}[(The categories $\mathsf{Tree}(A)$)]
\label{def:Tree(A)}
If $T$ and $T'$ are generalized $A$--labelled trees, we say that $T \leq
T'$ if $T$ can be obtained from $T'$ by a sequence of edge collapses (of
either internal or root edges) or bud collapses. This makes the set
$\mathsf{Tree}(A)$ of isomorphism classes of generalized $A$--labelled
trees into a poset and hence a category. Standard $A$--labelled trees (as
defined in Section~\ref{sec:trees}) are also generalized $A$--labelled
trees and $\mathsf{T}(A)$ is a full subcategory of $\mathsf{Tree}(A)$. See
Figure \ref{fig:gen_trees} for pictures of $\mathsf{Tree}(1)$ and
$\mathsf{Tree}(2)$.
\end{definition}

\begin{definition} \label{def:gen_weighting}
We don't need to change the definition of a \emph{weighting} for
generalized trees: it is an assignment of lengths to the edges of a tree
such that the root-leaf distances all equal $1$. As before, we write
$w(T)$ for the space of weightings on the generalized tree $T$. The
following result generalizes Lemma \ref{lem:W}.
\end{definition}

\begin{lemma} \label{lem:gen_W}
Let $T$ be a generalized $A$--labelled tree with $n$ (internal)
vertices. Then $w(T)$ is homeomorphic to $D^n$ and the boundary $\partial
w(T) \isom S^{n-1}$ consists of those points in which some edge of $T$
has length zero.
\end{lemma}
\begin{proof}
The labelling plays no role in the space of weightings so we can ignore
it. Picture $T$ as a collection of (non-generalized) trees $T_1,\dots,T_k$
attached at their roots. Suppose $T_j$ has $n_j$ vertices so that $n =
\sum n_j$. Then we have
\[ w(T) \isom w(T_1) \times \dots \times w(T_k) \isom D^{n_1} \times
\dots \times D^{n_k} \isom D^n. \]
Under this decomposition, a point is in the boundary of $w(T)$ if and
only if any of it is in the boundary of any of the $w(T_j)$. That is,
if and only if any of the edges of $T$ has length zero.
\end{proof}

\begin{definition} \label{def:gen_w_functor}
The `space of weightings' functor $w(-)\co \mathsf{T}(A) \to
\spaces$ of Definition \ref{def:w_functor} can be extended to all of
$\mathsf{Tree}(A)$. To do this, we have to say what happens when we apply
$w(-)$ to a morphism $T_b \to T$ coming from a bud collapse (for $b$
a bud in a tree $T$). Given a weighting of $T_b$ we get a weighting of
$T$ by giving length zero to all the leaf edges attached to $b$. This
defines a map
\[ w(T_b) \to w(T) \]
and it is not hard to see that this does indeed give us a functor
\[ w(-)\co \mathsf{Tree}(A) \to \spaces \]
as claimed. Adding a disjoint basepoint we get a functor
\[ w(-)_+\co \mathsf{Tree}(A) \to \based. \]
\end{definition}

\subsection{The two-sided bar construction} \label{sec:two-sided}

We now update Definition \ref{def:formal_bar} to the two-sided case. Along
with the spaces of weightings the key parts of this definition were
functors
\[ P_A(-) \co \mathsf{T}(A)^\text{op} \to \cat{C}. \]
The appropriate generalizations of these to functors on
$\mathsf{Tree}(A)^\text{op}$ are as follows.

\begin{definition} \label{def:(R,P,L)}
Let $P$ be a reduced operad in $\cat{C}$ with right module $R$ and
left module $L$. We define functors $(R,P,L)_A \co \mathsf{Tree}(A)^\text{op}
\to \cat{C}$ by \footnote{It is a serendipitous fact of our terminology
for trees that the \textbf{r}ight module $\mathbf{R}$ relates to the
\textbf{r}oots of our trees and the \textbf{l}eft module $\mathbf{L}$
relates to the \textbf{l}eaves.}
\[ (R,P,L)_A(T) := R(i(r)) \barwedge \Barwedge_{\text{vertices }
v \in T} P(i(v)) \barwedge \Barwedge_{\text{leaves } l \in T}
L(\iota^{-1}l). \]
Recall that $i(v)$ denotes the set of incoming edges to a vertex $v \in
T$. Here $\iota$ denotes the labelling surjection from $A$ to the set
of leaves of $T$, so that $\iota^{-1}l$ is the set of labels attached
to the leaf $l$.

To complete the definition, we have to give the effect of $(R,P,L)_A(-)$
on morphisms in $\mathsf{Tree}(A)$. Notice that $\mathsf{Tree}(A)$
is generated by the morphisms corresponding to
\begin{enumerate}
\item collapse of root edges,
\item collapse of internal edges, and
\item bud collapse.
\end{enumerate}
We will describe the effect of $(R,P,L)_A(-)$ on each of these types of
generating morphism and then check that they are compatible.

(1)\qua Suppose first that $e$ is a root edge of the generalized
$A$--labelled tree $T$. Then we have a morphism $T/e \to T$ that collapses
$e$. The corresponding morphism
\[ (R,P,L)_A(T) \to (R,P,L)_A(T/e) \]
is given by the map
\[ R(i(r)) \barwedge P(i(v)) \to R(i(r \circ v)) \]
that comes from the right $P$--module structure on $R$. Here $v$ is the
upper endpoint of the edge $e$ in $T$. Notice that $r \circ v$ is the
root element in $T/e$.

(2)\qua Now suppose that $e$ is an internal edge of $T$. The morphism
\[ (R,P,L)_A(T) \to (R,P,L)_A(T/e) \]
is then given (as in Definition \ref{def:formal_bar}) by the partial
composition map
\[ P(i(u)) \barwedge P(i(v)) \to P(i(u \circ v)) \]
for the operad $P$ where $u,v$ are the endpoints of $e$.

(3)\qua Finally, suppose that $b$ is a bud in the generalized $A$--labelled
tree $T$. The required map
\[ (R,P,L)_A(T) \to (R,P,L)_A(T_b) \]
comes from the map
\[ P(i(b)) \barwedge L(\iota^{-1}l_1) \barwedge \dots \barwedge
L(\iota^{-1}l_r) \to L(\iota_1^{-1}b) \]
that is part of the left $P$--module structure on $L$. Here
$l_1,\dots,l_r$ are the $b$--leaves in $T$ and we have
\[ \iota^{-1}b = \coprod_{i=1}^{r} \iota^{-1}l_i \]
from the definition of bud collapse, where $\iota_1$ is the $A$--labelling
of $T_b$.

The associativity conditions for $P$ to be an operad and for $R$ and $L$
to be $P$--modules ensure that these choices indeed determine a
functor $\mathsf{Tree}(A)^\text{op} \to \cat{C}$.
\end{definition}

\begin{definition}[(Two-sided bar construction)] \label{def:two-sided_bar}
Let $P$ be a reduced operad in $\cat{C}$ with right module $R$ and left
module $L$ as above. The \emph{bar construction on $P$ with coefficients
in $R$ and $L$} is the symmetric sequence $B(R,P,L)$ defined by the coends
\[ B(R,P,L)(A) := \int^{T \in \mathsf{Tree}(A)} w(T)_+ \otimes
(R,P,L)_A(T) \]
for finite sets $A$. A bijection $A \to A'$ determines an isomorphism
of categories $\mathsf{Tree}(A) \to \mathsf{Tree}(A')$ under which the
pairs of functors $w_A(-)$, $w_{A'}(-)$ and $(R,P,L)_A$, $(R,P,L)_{A'}$
correspond. It therefore induces an isomorphism
\[ B(R,P,L)(A) \to B(R,P,L)(A'). \]
So we obtain a symmetric sequence $B(R,P,L)$.

There is a more informal description of this bar construction that
generalizes that of $B(P)$ from Definition \ref{def:bar(operad)}. For
a finite set $A$, a point in $B(R,P,L)(A)$ consists of
\begin{itemize}
\item a weighted generalized $A$--labelled tree $T$,
\item a point in $R(i(r))$ where $r$ is the root of $T$,
\item a point in $P(i(v))$ for each vertex $v \in T$, and
\item a point in $L(\iota^{-1}l)$ for each leaf $l \in T$.
\end{itemize}
These are subject to identifications that tell us what happens when the
lengths of some of the edges tend to zero. When a root edge tends to zero
we use the right $P$--module structure map for $R$. When an internal edge
tends to zero we use the composition map for $P$. When a collection of
leaf edges attached to a bud tend to zero (note that the leaf edges
attached to a particular bud must all have the same length in a weighting)
we use the left $P$--module structure for $L$. Finally, of course, we
identify to the basepoint in $B(R,P,L)(A)$ if any of the chosen points
in $R(i(r)),P(i(v)),L(\iota^{-1}l)$ are the basepoint there.
\end{definition}
We now recall the simplicial version of the two-sided bar construction
for an operads and modules over them.

\begin{definition}[(Simplicial two-sided bar construction)]
\label{def:two-sided_simp}
Let $P$ be an operad in $\cat{C}$ with right module $R$ and left module
$L$. The \emph{simplicial bar construction on $P$ with coefficients in
$L$ and $R$} is the simplicial object $\mathcal{B}_{\bullet}(R,P,L)$
in the category of symmetric sequences in $\cat{C}$ with
\[ \mathcal{B}_n(R,P,L) := R \circ \underset{n}{\underbrace{P \circ
\dots \circ P}} \circ L. \]
The face maps
\[ d_i \co \mathcal{B}_n(R,P,L) \to \mathcal{B}_{n-1}(R,P,L) \]
for $i = 1,\dots,n-1$ are given by the operad composition map $P \circ
P \to P$ applied to the \ord{$i$} and \ord{$i+1$} factors. The face
map $d_0$ is given by the right module structure $R \circ P \to R$
and $d_n$ is given by the left module structure $P \circ L \to L$. The
degeneracy map
\[ s_j \co \mathcal{B}_n(R,P,L) \to \mathcal{B}_{n+1}(R,P,L) \]
is given by using the unit map $I \to P$ to insert an extra copy of $P$
between the \ord{$j$} and \ord{$j+1$} factors.
\end{definition}

\begin{proposition} \label{prop:two-sided_simp}
Let $P$ be a reduced operad in $\cat{C}$ with right module $R$ and left
module $L$. The bar construction of Definition \ref{def:two-sided_bar}
is isomorphic to the geometric realization of the simplicial bar
construction:
\[ B(R,P,L) \isom |\mathcal{B}_{\bullet}(R,P,L)|. \]
\end{proposition}
\begin{proof}
This is a straightforward extension of the argument used to prove
Proposition~\ref{prop:simp}.
\end{proof}
Our first example of the two-sided bar construction is that the reduced
bar construction a lone operad is a special case.

\begin{example} \label{ex:old_bar=new}
Let $P$ be a reduced operad in $\cat{C}$ and take $R = L = I$ the unit
symmetric sequence. Recall that $I$ is a left and right module over any
augmented operad. It is easy to see from the definitions that for the
simplicial bar constructions we have
\[ \mathcal{B}_{\bullet}(I,P,I) \isom \mathcal{B}_{\bullet}(P). \]
This tells us that
\[ B(I,P,I) \isom B(P), \]
but we can see this directly as well. First notice that
\[ (I,P,I)_A(T) \isom \begin{cases} P_A(T) & \text{if $T \in
\mathsf{T}(A)$}; \\ \basept & \text{otherwise}. \end{cases} \]
The reason for this is as follows. Because $I(n) = \basept$ for $n > 1$,
we have $(I,P,I)_A(T) = \basept$ whenever $T$ has more than one root edge,
or when any leaf has more than one label. These are precisely the
generalized $A$--labelled trees not in $\mathsf{T}(A)$. For $T \in
\mathsf{T}(A)$ we have
\[ (I,P,I)_A(T) = I(1) \barwedge P_A(T) \barwedge I(1) \barwedge \dots
\barwedge I(1) \isom P_A(T). \]

This calculation means that only the objects $T \in \mathsf{T}(A)$
contribute to the calculation of the coend in Definition
\ref{def:two-sided_bar}. However, we still have to take into account
morphisms $U \to T$ with $U \notin \mathsf{T}(A)$. This amounts to
collapsing to the basepoint those weighted trees in which either the
root edge or a leaf edge has length zero (since these are the images of
the maps $w(U) \to w(T)$). All together this tells us that $B(I,P,I)(A)$
is equal to the coend
\[ \int^{T \in \mathsf{T}(A)} \overline{w}(T) \otimes P_A(T) \]
where $\overline{w}(T)$ is the quotient of $w(T)$ by the weightings which
have either root or leaf edge of length zero. This is precisely
$B(P)(A)$. Therefore we have $B(I,P,I) \isom B(P)$ as claimed.
\end{example}

\begin{example} \label{ex:two-sided_bar}
It is easy to see that $B(R,P,L)(1) = R(1) \barwedge L(1)$. We have
already seen (Figure \ref{fig:gen_trees}) that there are three objects in
$\mathsf{Tree}(2)$. From this we see that $B(R,P,L)(2)$ is the homotopy
pushout of the following diagram
\[ \begin{diagram} \node{R(1) \barwedge P(2) \barwedge L(1) \barwedge
L(1)} \arrow{s} \arrow{e} \node{R(1) \barwedge L(2)} \\
        \node{R(2) \barwedge L(1) \barwedge L(1)}
\end{diagram} \]
If $R = L = I$, the bottom-left and top-right objects are $\basept$
and the top-left object is $P(2)$. So we recover
\[ B(P)(2) = B(I,P,I)(2) = \Sigma P(2). \]
\end{example}

\begin{definition}[(Bar constructions for modules)]
Let $P$ be a reduced operad in $\cat{C}$ and let $R$ be a right
$P$--module. We define the \emph{bar construction on $R$} by
\[ B(R) := B(R,P,I) \]
where $I$, as previously, is the unit for the composition product of
symmetric sequences. If $L$ is a left $P$--module, its \emph{bar
construction} is
\[ B(L) := B(I,P,L). \]
We trust that it will not be confusing to use the same notation for the
bar construction of right and left modules.
\end{definition}

\begin{example} \label{ex:two-sided_modules}
Applying Example \ref{ex:two-sided_bar} to the one-sided case we see that
\[ B(R)(1) \isom R(1); \;\; B(R)(2) \isom \hocofib(R(1) \barwedge P(2)
\to R(2)) \]
and
\[ B(L)(1) \isom L(1); \;\; B(L)(2) \isom \hocofib(P(2) \barwedge L(1)
\barwedge L(1) \to L(2)). \]
\end{example}

\begin{definition}[(Cobar constructions for comodules)]
All the constructions of this section can be applied to operads and
modules in $\cat{C}^\text{op}$, that is, to cooperads and comodules in
$\cat{C}$. We summarize the results.

If $Q$ is a reduced cooperad in $\cat{C}$ with left comodule $L$ and right
comodule $R$, the formula of Definition \ref{def:(R,P,L)} defines functors
\[ (R,Q,L)_A(-)\co \mathsf{Tree}(A) \to \cat{C} \]
for each finite set $A$ and we define the \emph{cobar construction}
on $Q$ with coefficients in $R$ and $L$ to be the symmetric sequence
$\Omega(R,Q,L)$ with
\[ \Omega(R,Q,L)(A) := \int_{T \in \mathsf{Tree}(A)}
\Map_{\cat{C}}(w(T)_+,(R,Q,L)_A(T)). \]
This is isomorphic to the totalization of the two-sided cosimplicial cobar
construction on $Q$ with coefficients in $R$ and $L$. The \emph{cobar
construction on $R$} is
\[ \Omega(R) := \Omega(R,Q,I) \]
and the \emph{cobar construction on $L$} is
\[ \Omega(L) := \Omega(I,Q,L). \]
\end{definition}

\begin{example}
Taking $R = L = I$ we recover the cobar construction of
Section~\ref{sec:cobar}:
\[ \Omega(I,Q,I) \isom \Omega(Q). \]
\end{example}

\begin{example} \label{ex:two-sided_cobar}
Taking the duals of the results of Example \ref{ex:two-sided_bar} we
see that
\[ \Omega(R,Q,L)(1) \isom R(1) \barwedge L(1) \]
and that $\Omega(R,Q,L)(2)$ is the homotopy pullback of
\[ \begin{diagram}
    \node[2]{R(1) \barwedge L(2)} \arrow{s} \\
    \node{R(2) \barwedge L(1) \barwedge L(1)} \arrow{e} \node{R(1)
    \barwedge Q(2) \barwedge L(1) \barwedge L(1).}
\end{diagram} \]
In particular,
\[ \Omega(R)(1) \isom R(1); \;\; \Omega(R)(2) \isom \hofib(R(2) \to R(1)
\barwedge Q(2)) \]
and
\[ \Omega(L)(1) \isom L(1); \;\; \Omega(L)(2) \isom \hofib(L(2) \to Q(2)
\barwedge L(1) \barwedge L(1)). \]
\end{example}

\subsection{Structure maps for bar constructions on modules}
\label{sec:bar(modules)_maps}

In this section we use similar methods to Section~\ref{sec:cooperad} to
show that the bar construction on a $P$--module (that is, a single
left or right module) is a comodule over the cooperad $B(P)$. In fact,
we will construct maps of the form
\begin{equation} B(R,P,L) \to B(R,P,I) \varcirc B(I,P,L)
\label{eq:general} \end{equation}
where $\varcirc$ is the composition of symmetric sequences defined using
the product in $\cat{C}$ rather than the coproduct (see Remark
\ref{rem:cooperads}). Taking $R = I$ and recalling that $B(I,P,I) = B(P)$
we obtain a left $B(P)$--comodule structure on $B(L) = B(I,P,L)$.
Similarly, taking $L = I$ we get a right $B(P)$--comodule structure on
$B(R)$. Notice that taking $R = L = I$ we recover the cooperad structure
on $B(P)$.

The definition of the map (\ref{eq:general}) is a relatively
straightforward generalization of the cooperad structure on $B(P)$. We
start by describing the grafting and ungrafting processes for generalized
trees.

\begin{figure}[htbp]
\begin{center}
\input{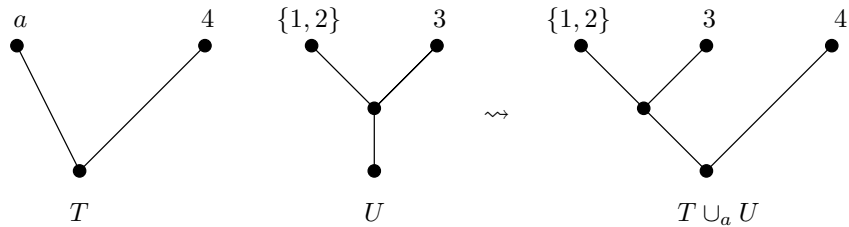}
\caption{Grafting generalized labelled trees}
\label{fig:ungraft}
\end{center}
\end{figure}

\begin{definition}[(Grafting for generalized trees)]
Let $T$ be a generalized $A$--labelled tree and $U$ a generalized
$B$--labelled tree and let $a$ be an element of $A$. We will define
the \emph{grafting of $U$ onto $T$} only if $T$ and $U$ satisfy the
following conditions:
\begin{itemize}
\item The root of $U$ has only one incoming edge.
\item The leaf of $T$ labelled by $a$ is labelled only by $a$ and no
other elements of $A$.
\end{itemize}
In this case, the grafted tree $T \cup_a U$ is defined exactly as in
Definition \ref{def:grafting} by identifying the root edge of $U$ to the
$a$--leaf edge of $T$. Figure \ref{fig:ungraft} gives an example.
\end{definition}

To define the maps (\ref{eq:general}) we will need to graft trees onto
all of the leaf edges of the base tree $T$. To do this, we must assume
that all the leaves of $T$ only have one label, so that $T$ satisfies
the stronger condition for a labelling we required in Definition
\ref{def:labelling}. Notice also that the trees $U$ we are to graft onto
$T$ satisfy the stronger root condition of Definition \ref{def:tree}.
The following definitions will help us talk about trees of these types.

\begin{definition} \label{def:more_trees}
For a finite set $A$, we define the following full subcategories of
$\mathsf{Tree}(A)$:
\begin{gather*} \mathsf{T}_\text{root}(A) := \{ T \in \mathsf{Tree}(A) |
\text{ the root of $T$ has only one incoming edge}\} \\
     \mathsf{T}_\text{leaf}(A) := \{ T \in \mathsf{Tree}(A) | \text{ the
     leaves of $T$ are labelled bijectively by $A$}\}. \end{gather*}
Notice that $\mathsf{T}(A) = \mathsf{T}_\text{root}(A) \cap
\mathsf{T}_\text{leaf}(A)$.
\end{definition}

\begin{definition}
Let $A = \coprod_{j \in J} A_j$ be a partition of $A$ into nonempty
subsets. Given trees $U_j \in \mathsf{T}_\text{root}(A_j)$ and $T \in
\mathsf{T}_\text{leaf}(J)$, we denote the tree obtained by grafting all the
$U_j$ onto $T$ at the appropriate places by
\[ T \cup_J U_j. \]
We say that a generalized $A$--labelled tree is \emph{of type}
$\{A_j\}$ if it is of the form $T \cup_J U_j$ for some such $T$
and $U_j$. The correct generalization of the functor of Proposition
\ref{prop:category_cooperad} is then a functor
\[ \mathsf{Tree}(A)_+ \to \mathsf{T}_\text{leaf}(J)_+
\smsh \mathsf{T}_\text{root}(A_{j_1})_+ \smsh \dots \smsh
\mathsf{T}_\text{leaf}(A_{j_r})_+ \]
that breaks the tree $(T \cup_J U_j)$ into its components $T$ and the
$U_j$ and sends a tree not of type $\{A_j\}$ to the initial object on
the right-hand side. This `ungrafting' functor is the basis of the map
(\ref{eq:general}).
\end{definition}

Our new categories of trees can be used as the base categories for
defining the one-sided bar constructions. For this we need the appropriate
spaces of weightings.

\begin{definition}
For each finite set $A$ we define a functor
\[ w_\text{leaf}(-)\co \mathsf{T}_\text{leaf}(A) \to \based \]
where $w_\text{leaf}(T)$ is the quotient of $w(T)$ by the space of weightings
in which some leaf edge has length zero, and a functor
\[ w_\text{root}(-)\co \mathsf{T}_\text{root}(A) \to \based \]
where $w_\text{root}(T)$ is the quotient of $w(T)$ by the space of weightings
in which the root edge has length zero.
\end{definition}

\begin{lemma}
Let $P$ be a reduced operad in $\cat{C}$ with right module $R$ and left
module $L$. Then the one-sided bar constructions are given by
\[ B(R)(A) = B(R,P,I)(A) \isom \int^{T \in \mathsf{T}_\text{leaf}(A)}
w_\text{leaf}(T) \otimes (R,P,I)_A(T) \]
and
\[ B(L)(A) = B(I,P,L)(A) \isom \int^{T \in \mathsf{T}_\text{root}(A)}
w_\text{root}(T) \otimes (I,P,L)_A(T). \]

\end{lemma}
\proof
These calculations are similar to that in Example \ref{ex:old_bar=new}
where we showed that $B(P) = B(I,P,I)$. They use the facts that
\[ (R,P,I)_A(T) = \basept \text{ for $T \notin \mathsf{T}_\text{leaf}(A)$} \]
and
\[ (I,P,L)_A(T) = \basept \text{ for $T \notin \mathsf{T}_\text{root}(A)$}.
  \eqno{\Box} \]

The final piece of the puzzle is the construction of a map analogous
to (\ref{eq:key}) that tells us how to weight the trees obtained from
ungrafting.

\begin{definition} \label{def:gen_key}
Let $A = \coprod_{j \in J} A_j$ be a partition of the finite set $A$
into nonempty subsets. Given trees $T \in \mathsf{T}_\text{leaf}(J)$ and $U_j
\in \mathsf{T}_\text{root}(A_j)$ we define a map
\[ w(T \cup_J U_j)_+ \to w_\text{leaf}(T) \smsh w_\text{root}(U_{j_1}) \smsh \dots
\smsh w_\text{root}(U_{j_r}) \]
by the obvious generalization of the construction of the maps
$\overline{w}(T \cup_a U) \to \overline{w}(T) \smsh \overline{w}(U)$ in Definition
\ref{def:cooperad_maps}.
\end{definition}

\begin{definition} \label{def:general}
Putting together all these ingredients we construct maps
\[ B(R,P,L)(A) \to B(R,P,I)(J) \barwedge B(I,P,L)(A_{j_1}) \barwedge
\dots \barwedge B(I,P,L)(A_{j_r}).\]
In an analogous way to Definition \ref{def:formal_cooperad_maps},
these come from the maps of Definition \ref{def:gen_key} together with
the isomorphisms
\[ (R,P,L)_A(T \cup_J U_j) \to (R,P,I)_J(T) \barwedge
(I,P,L)_{A_{j_1}}(U_{j_1}) \barwedge \dots \barwedge
(I,P,L)_{A_{j_r}}(U_{j_r}). \]
Together these maps make up the map of symmetric sequences
\[ B(R,P,L) \to B(R,P,I) \varcirc B(I,P,L) \]
as promised.
\end{definition}

\begin{prop} \label{prop:general}
Let $P$ be a reduced operad in $\cat{C}$ with right module $R$ and left
module $L$. The maps of Definition \ref{def:general} determine a right
$B(P)$--comodule structure on $B(R)$ and a left $B(P)$--comodule structure
on $B(L)$.
\end{prop}
\begin{proof}
Taking $L = I$ in \ref{def:general} we get the right comodule structure on
$B(R)$. Taking $R = I$ we get the left comodule structure on $B(L)$. We
have to check the appropriate associativity and unit axioms. This is a
generalization of the work of Section~\ref{sec:cooperad}. We leave the
reader to write out all the details, including the diagrams corresponding
to Figure \ref{fig:assoc2}.
\end{proof}

\begin{cor}
Dually, suppose that $Q$ is a reduced cooperad in $\cat{C}$ with right
comodule $R$ and left comodule $L$. Then there is a map
\[ \Omega(R,Q,I) \circ \Omega(I,Q,L) \to \Omega(R,Q,L) \]
that makes $\Omega(R)$ into a right $\Omega(Q)$--module (by taking $L =
I$) and $\Omega(L)$ into a left $\Omega(Q)$--module (by taking $R = I$).
\end{cor}
\begin{proof}
Apply Proposition \ref{prop:general} to $Q$ considered as an operad
in $\cat{C}^\text{op}$.
\end{proof}

This completes our descriptions of the bar and cobar constructions
for operads, cooperads, modules and comodules. We turn now to our main
application of this theory -- the Goodwillie derivatives of the identity
functor.

\section{Application to the calculus of functors} \label{sec:application}

In this section we describe our application of bar and cobar constructions
to Goodwillie's calculus of homotopy functors. The main result is that the
derivatives of the identity form an operad in spectra. We now assume that
$\cat{C}$ is a suitable category $\spectra$ of spectra, for example, the
$S$--modules of EKMM \cite{elmendorf/kriz/mandell/may:1997} (see Example
\ref{ex:categories}(2)).

Let $I\co \based \to \based$ be the identity functor on based spaces. The
Goodwillie derivatives of $I$ can be described in terms of the partition
poset complexes \cite{arone/mahowald:1999}. We recall one of the ways
to define these.

\begin{definition}
A \emph{partition} of a finite set $A$ is an equivalence relation
on $A$. Let $K(A)$ be the poset formed by the partitions of $A$ with
$\lambda \leq \mu$ if $\lambda$ is \emph{finer} than $\mu$, that is, if
the set of relations for $\lambda$ is contained in the set of relations
for $\mu$. The category $K(A)$ has an initial object $\widehat{0}$ and a
terminal object $\widehat{1}$. Let $K_0(A) = K(A) - \widehat{0}$, the category
of \emph{proper} partitions, and $K_1(A) = K(A) - \widehat{1}$, the category
of \emph{non-trivial} partitions. Note that the group $\Sigma_A$ of
permutations of $A$ acts on all of these categories in an obvious way.
\end{definition}

\begin{definition}[(Partition poset complexes)]
For a finite set $A$, the \emph{partition poset complex} $\Delta(A)$ is
the geometric realization of the following simplicial set $T(A)_{\bullet}$
formed from the nerves of these categories of partitions:
\[ T(A)_{\bullet} = \frac{N_{\bullet}K(A)}{N_{\bullet}K_0(A) \union
N_{\bullet}K_1(A)} \]
So the $n$--simplices in $T(A)_{\bullet}$ are sequences of $n+1$
partitions
\[ \lambda_0 \leq \lambda_1 \leq \dots \leq \lambda_n \]
with a sequence identified to the basepoint if it does not have
both $\lambda_0 = \widehat{0}$ and $\lambda_n = \widehat{1}$. The face and
degeneracy maps are given by respectively removing partitions from
the sequence and repeating terms in the usual way for the nerve of a
category. The simplicial set $T(A)_{\bullet}$ is pointed and so its
geometric realization $\Delta(A)$ is a based space. A bijection $A \to
A'$ induces an isomorphism $\Delta(A) \to \Delta(A')$ that makes $\Delta$
into a symmetric sequence in $\based$.
\end{definition}

\begin{remark}
What we are calling the partition poset complex is the suspension of
the complex $K_n$ of \cite{arone/mahowald:1999}. The simplicial set
$T(n)_{\bullet}$ is isomorphic to that called $T_n$ in Definition 1.1
of \cite{arone/mahowald:1999}.
\end{remark}

\begin{prop}[Arone--Mahowald, \protect\cite{arone/mahowald:1999}]
\label{prop:AM}
The derivatives of the identity are modelled by the dual spectra of the
finite complexes $\Delta(n) = \Delta({\{1,\dots,n\}})$:
\[ \partial_nI \simeq \Map_{\spectra}(\Delta(n),S) \]
The action of the symmetric group $\Sigma_n$ on $\Delta(n)$ induces an
action on the dual spectrum and this agrees with the action that comes
with the spectrum $\partial_nI$.
\end{prop}

The key observation (apparently due to Greg Arone) is that the partition
poset complexes can be described as spaces of trees. We can interpret
these as the spaces of a bar construction.

\begin{definition}
Let $\underline{S^0}$ be the operad in based spaces with
\[ \underline{S^0}(A) := S^0 \]
for all finite sets $A$ and with all composition maps equal to the
identity on $S^0$. This is the operad for commutative monoids of based
spaces.
\end{definition}

\begin{lemma} \label{lem:ppc=bar}
The partition poset complex $\Delta(A)$ is homeomorphic to the bar
construction $B(\underline{S^0})(A)$.
\end{lemma}
\begin{proof}
We have already seen that $B(\underline{S^0})$ is homeomorphic to the
realization of the simplicial bar construction on $\underline{S^0}$. It
is therefore enough to show that the simplicial set $T(A)_{\bullet}$ used
to define $\Delta(A)$ is also given by this simplicial bar construction.

A non-basepoint $n$--simplex in $T(A)$ is an increasing sequence of
partitions of $A$ of length $n-1$. On the other hand the based set of
$n$--simplices in the simplicial bar construction is
\[ \underset{n}{\underbrace{\underline{S^0} \circ \dots \circ
\underline{S^0}}}(A). \]
But this is equal to the wedge over increasing sequences of partitions of
length $n-1$ of $S^0$. Hence we see that the two sets of $n$--simplices
are the same. The face and degeneracy maps in each case correspond to
removing a partition and repeating a partition respectively. We therefore
have isomorphic simplicial sets.
\end{proof}

\begin{remark} \label{rem:vallette}
In \cite{vallette:2004}, Bruno Vallette describes the notion of a
\emph{$P$--partition} for an operad $P$ in $\mathsf{Set}$. The
$P$--partitions form a poset whose nerve (or \emph{order complex} in
\cite{vallette:2004}) is isomorphic to the bar construction $B(P_+)$
(where we are considering $P$ as a discrete operad in unbased spaces and
adding a disjoint basepoint). Lemma \ref{lem:ppc=bar} is the special
case of this fact when $P$ is the `commutative operad' in $\mathsf{Set}$,
that is, with $P(n) = \basept$ for all $n$.
\end{remark}

\begin{cor} \label{cor:operad}
Let $\partial_nI$ denote the model of the \ord{$n$} derivative of the
identity given by
\[ \partial_nI = \Map_{\spectra}(\Delta(n),S). \]
Then we have
\[ \partial_nI = \Omega(\dual \underline{S^0})(n). \]
In particular, the derivatives of the identity form an operad in
spectra. We denote this operad by $\partial_*I$.
\end{cor}
\begin{proof}
We have
\[ \partial_nI = \Map(\Delta(n),S) = \dual B(\underline{S^0})(n) =
\Omega(\dual \underline{S^0})(n) \]
by Lemma \ref{lem:ppc=bar} and Proposition \ref{prop:duality} (which
applies since all the spaces in $\underline{S^0}$ are $S^0$).
\end{proof}

\begin{remark}
The derivatives of the identity are the cobar construction on the cooperad
$\underline{S}$ in spectra with
\[ \underline{S}(A) = \dual \underline{S^0}(A) = S \]
where $S$ is the sphere spectrum, for all finite sets $A$ and with all
cocomposition maps the canonical isomorphisms. This is the analogue for
spectra of the cooperad for cocommutative coalgebras.
\end{remark}

\begin{remark} \label{rem:modules}
We can use the constructions of Section~\ref{sec:bar(modules)} to get
modules over the operad $\partial_*I$. If $C$ is a comodule over
$\underline{S}$ then its cobar construction $\Omega(C)$ is a
$\partial_*I$--module. We give two examples:

(1)\qua Let $X$ be a based space. Then the suspension spectrum
$\Sigma^{\infty}X$ is a $\underline{S}$--coalgebra (that is, just a
commutative coalgebra) with comultiplication given by the (reduced)
diagonal map on $X$:
\[ \Sigma^{\infty}X \to \Sigma^{\infty}(X \smsh X) \isom \Sigma^{\infty}X
\smsh \Sigma^{\infty}X. \]
As remarked in Definition \ref{def:comodule}, a coalgebra over a cooperad
$Q$ determines a left $Q$--comodule. Thus we obtain a left
$\underline{S}$--comodule $\underline{\Sigma^{\infty}X}$. We now take the
cobar construction to get a left $\partial_*I$--module
\[ M_X := \Omega(\underline{\Sigma^{\infty}X}) =
\Omega(I,\underline{S},\underline{\Sigma^{\infty}X}) \]
(where $I$ in this formula denotes the unit symmetric sequence
of Definition \ref{def:compprod}). From the calculations of
\ref{ex:two-sided_cobar} we find that
\[ M_X(1) = \Sigma^{\infty}X \]
and
\[ \begin{split}  M_X(2) &\isom \hofib(\Sigma^{\infty}X \to
\Sigma^{\infty}X \smsh \Sigma^{\infty}X) \\
    &\simeq \Sigma^{-1} \hocofib(\Sigma^{\infty}X \to \Sigma^{\infty}X
    \smsh \Sigma^{\infty}X) \\
    &\simeq \Sigma^{-1} \Sigma^{\infty} \hocofib(X \to X \smsh X) \\
\end{split} \]
So $M_X(2)$ is (up to homotopy and a desuspension) the mapping cone of the
reduced diagonal on $X$. Further work is needed to analyze the spectra
$M_X(n)$ for larger $n$. In Section~\ref{sec:modules_specseq} we will
look at ways to calculate the homology of these spectra.

(2)\qua A moment's thought will reveal that a right
$\underline{S}$--comodule is precisely the same thing as a functor
\[ (\mathsf{FinSets},\twoheadrightarrow) \longrightarrow \spectra \]
where the left-hand side is the category of finite sets with morphisms
given by the surjections. Work in progress by Greg Arone has demonstrated
a relationship between such functors and the Goodwillie calculus of
homotopy functors $F$ from based spaces to spectra.
\end{remark}

\begin{remark}[(Derivatives of general homotopy functors)]
The derivatives of any homotopy functor $F$ form a symmetric sequence in
spectra and it is natural to ask how these symmetric sequences might be
related for different functors. We conjecture that there is in general
a map of symmetric sequences
\[ \partial_*F \circ \partial_*G \to \partial_*(FG) \]
for any two homotopy functors $F,G\co \based \to \based$ such that
$F(\basept) = \basept$, where $FG$ denotes the composite of $F$ and $G$.
These maps should have suitable associativity properties that taking $F =
G = I$ would recover an operad structure on $\partial_*I$ equivalent to
the one we have constructed in this section. Similarly, taking $F = I$
would yield the structure of a left $\partial_*I$--module on $\partial_*G$
and taking $G = I$ a right $\partial_*I$--module structure on
$\partial_*F$. The main obstacle at present for constructing these maps
is finding good models for the derivatives of a general functor in a
symmetric monoidal category $\spectra$ of spectra. In the case of the
identity functor we were fortunate that such models naturally arose from
the partition poset complexes.
\end{remark}

\section{Homology of the bar and cobar constructions and Koszul duality}
\label{sec:alg}

In this section we look at spectral sequences for calculating the homology
of the bar and cobar constructions on operads and cooperads in based
spaces or spectra. It turns out that we can relate the $E^1$--term
of these spectral sequences to the algebraic bar and cobar
constructions described in, for example, \cite{getzler/jones:1994} and
\cite{fresse:2004}. This leads to a link with Koszul duality which says,
briefly, that if the homology of the reduced operad $P$ is Koszul, then
the homology of $B(P)$ is its Koszul dual cooperad, and dually, if the
homology of the cooperad $Q$ is Koszul then the homology of $\Omega(Q)$
is its Koszul dual operad. This supports the point-of-view that the bar
construction for an operad in based spaces or spectra is the analogue
of the Koszul dual for an algebraic operad.

Here is a summary of this section. We start in Section~\ref{sec:homology}
by recalling how the homology (with coefficients in the commutative ring
$k$) of an operad in based spaces or spectra has the structure of an
operad in graded $k$--modules. Then in Section~\ref{sec:filter}, the main
work of the chapter begins and we describe the filtration of the bar
construction that gives rise to our spectral sequence and identify the
`filtration quotients'. This filtration is based on the number of vertices
in the trees that underlie the bar construction. We deal immediately with
the two-sided construction of Section~\ref{sec:two-sided}, recalling that
the construction for a lone operad is a special case of this. As usual,
for the cobar construction, we just dualize everything. That is, we get a
cofiltration, or tower, whose inverse limit is the cobar construction and
we identify the fibres of the stages in this tower. In
Section~\ref{sec:cofibrations} we give conditions under which the
inclusion maps of the filtrations are cofibrations, thus ensuring that our
`filtration quotients' are actually the homotopy cofibres of filtration.
This will allow us later to use our identification of these quotients to
calculate the $E^1$ term in the spectral sequence. This $E^1$ term turns
out to be given by the algebraic bar construction which we describe in
Section~\ref{sec:alg_bar}. We give a definition of this that emphasizes
its similarity to the topological version and show that this definition is
equivalent to that given by Getzler and Jones \cite{getzler/jones:1994}
and Fresse \cite{fresse:2004}. Then in Section~\ref{sec:specseq} we
finally set up the spectral sequence  and identify its $E^1$ term with the
algebraic bar construction as claimed. In Section~\ref{sec:koszul} we look
at Koszul operads and prove the result identifying the homology of the bar
construction on $P$ with the Koszul dual of the homology of $P$. Finally,
in Section~\ref{sec:modules_specseq} we use our spectral sequences to
investigate the homology of the $\partial_*I$--modules $M_X$ constructed
in Remark \ref{rem:modules}(1).

\subsection{Homology of topological operads} \label{sec:homology}

Throughout the chapter we fix a commutative ring $k$ and consider the
categories $\mathsf{Mod}_k$ of graded $k$--modules and $\mathsf{Ch}_k$
of chain complexes over $k$. First we describe the symmetric monoidal
structure on these categories.

\begin{definition} \label{def:ch_k}
The tensor product determines a symmetric monoidal structure on graded
$k$--modules with
\[ (M \otimes N)_r := \bigoplus_{p + q = r} M_p \otimes N_q \]
where the \emph{graded} symmetry isomorphism
\[ M \otimes N \to N \otimes M \]
is given by
\[ m \otimes n \mapsto (-1)^{|m||n|} n \otimes m \]
and the unit object is the graded module $k$ concentrated in degree
$0$. If $M$ and $N$ are chain complexes with differentials $d_M$ and $d_N$
respectively, we define a differential on $M \otimes N$ by
\[ d_{M \otimes N}(m \otimes n) := d_M(m) \otimes n + (-1)^{|m|}m \otimes
d_N(n). \]
This makes $\otimes$ into a symmetric monoidal structure on
$\mathsf{Ch}_k$ with the same unit $k$ endowed with the trivial
differential.
\end{definition}

Throughout this section we will use $H_*(-)$ to denote the homology with
coefficients in the commutative ring $k$ of an object in $\cat{C}$ when
$\cat{C}$ is either $\based$ or $\spectra$. If $\cat{C}$ is the category
$\based$ of based spaces, this is the \emph{reduced} homology.\footnote{We
stress that any homology group of a based space in this paper is meant
to be the \emph{reduced} homology.} If $\cat{C}$ is a category $\spectra$
of spectra, it is the spectrum homology $H_*(E) = \pi_*(Hk \smsh E)$. We
recall the K\"{u}nneth maps for these homology theories.

\begin{prop} \label{prop:kunneth}
Let $\cat{C} = \based$ or $\spectra$ and take $C,D \in \cat{C}$. Then
there is a natural map
\[ H_*(C) \otimes H_*(D) \to H_*(C \barwedge D) \]
that is an isomorphism if either $H_*(C)$ or $H_*(D)$ consists of flat
$k$--modules. These maps are symmetric monoidal in the sense that they
commute with the associativity and commutativity isomorphisms in the
categories $\cat{C}$ and $\mathsf{Mod}_k$.
\end{prop}

\begin{definition}
Let $M$ be a symmetric sequence in $\based$ or $\spectra$. Then we denote
by $H_*M$ the symmetric sequence of graded $k$--modules given by
\[ H_*M(A) := H_*(M(A)). \]
\end{definition}

The main result of this section is that the homology of a topological
operad or cooperad is, under suitable conditions, an operad or cooperad
in $\mathsf{Mod}_k$.

\begin{lemma}
Let $P$ be an operad in $\based$ or $\spectra$. Then $H_*P$ is an operad
of graded $k$--modules. If $P$ is reduced then so is $H_*P$. If $M$ is a
left (respectively, right) $P$--module, then $H_*M$ is a left
(respectively, right) $H_*P$--module.

Let $Q$ be a cooperad in $\based$ or $\spectra$ such that the homology
groups $H_*(Q(A))$ are flat $k$--modules. Then $H_*(Q)$ is a cooperad of
graded $k$--modules that is reduced if $Q$ is. If $C$ is a right
$Q$--comodule then $H_*(C)$ is a right $H_*(Q)$--comodule. If $C$ is a
left $Q$--comodule such that the $H_*(C(A))$ are flat $k$--modules then
$H_*(C)$ is a left $H_*(Q)$--comodule.
\end{lemma}
\begin{proof}
The operad structure maps are given by the composites
\[ H_*(P(A)) \otimes H_*(P(B)) \to H_*(P(A) \smsh P(B)) \to H_*(P(A
\cup_a B)) \]
and the unit by the map
\[ k \isom H_*(S) \to H_*(P)(1) \]
where $S$ denotes either $S^0$, the unit of $\based$, or the unit
of $\spectra$. To check the operad axioms we use the associativity
and commutativity of the K\"{u}nneth formula as stated in Proposition
\ref{prop:kunneth}. Clearly, if $P$ is reduced (so that the unit map $S
\to P(1)$ is an isomorphism) then so is $H_*P$. The structure maps for
$H_*M$ are defined similarly.

In the cooperad case we need the flatness condition. It allows us to
define cocomposition maps by
\[ H_*(Q(A \cup_a B)) \to H_*(Q(A) \smsh Q(B)) \isom H_*(Q(A)) \otimes
H_*(Q(B)) \]
using the inverse of the K\"{u}nneth map. The counit map is the composite
\[ H_*Q(1) \to H_*(S) \isom k \]
and again, if $Q$ is reduced, so is $H_*Q$. In the case of a right
comodule $C$ we similarly get comodule structure maps
\[ H_*(C(A \cup_a B)) \to H_*(C(A) \smsh Q(B)) \isom H_*(C(A)) \otimes
H_*(Q(B)) \]
where the K\"{u}nneth map is an isomorphism without any condition on
$H_*(C(A))$ (we are still assuming that the $H_*(Q(B))$ are flat). For
a left comodule, we do still need the flatness assumption.
\end{proof}

\begin{remark}
We can consider cohomology instead of homology in which case the
K\"{u}nneth isomorphism also requires a finiteness hypothesis. We get the
following results. If $Q$ is a cooperad in based spaces or spectra then
$H^*(Q)$ is an operad of graded $k$--modules. If $P$ is an operad with the
cohomology groups $H^*(P)$ finitely-generated flat $k$--modules then
$H^*(P)$ is a cooperad of graded $k$--modules. Similar results hold for
comodules and modules.
\end{remark}

\subsection{Filtering the bar construction} \label{sec:filter}
The spectral sequence we want to construct comes from a filtration
on the bar construction by the number of vertices in the underlying
trees. In this section we construct this filtration and calculate the
filtration quotients.

\begin{definition}[(Filtration on the category of trees)]
Write $\mathsf{Tree}_s(A)$ for the subcategory of $\mathsf{Tree}(A)$
whose objects are the (isomorphism classes of) trees with less than or
equal to $s$ (internal) vertices. We then have
\[ \mathsf{Tree}_0(A) \subset \mathsf{Tree}_1(A) \subset \dots \subset
\mathsf{Tree}_{|A|-1}(A) = \mathsf{Tree}(A). \]
Each $\mathsf{Tree}_s(A)$ is an \emph{initial} subcategory of
$\mathsf{Tree}(A)$. That is, if $U \leq T$ and $T \in \mathsf{Tree}_s(A)$
then $U \in \mathsf{Tree}_s(A)$. The filtration `quotients' are the
discrete categories
\[ \mathsf{Q}_s(A) := \mathsf{Tree}_s(A) - \mathsf{Tree}_{s-1}(A) \]
whose objects are the trees with precisely $s$ vertices. For each tree
$T \in \mathsf{Tree}(A)$ we write $|T|$ for the number of vertices of $T$.
\end{definition}

\begin{definition}[(Filtration on the two-sided bar construction)]
\label{def:filter}
For a reduced operad $P$ in $\cat{C}$ with right module $R$ and left
module $L$, define
\[ B(R,P,L)_s(A) := \int^{T \in \mathsf{Tree}_s(A)} w(T)_+ \otimes
(R,P,L)_A(T). \]
For varying finite sets $A$ these form a symmetric sequence in
$\cat{C}$. From the inclusion of categories $\mathsf{Tree}_{s-1}(A)
\subset \mathsf{Tree}_s(A)$ we get natural maps
\[ B(R,P,L)_{s-1}(A) \to B(R,P,L)_s(A). \]
In the case $\cat{C} = \based$, it is easy to see that the resulting
sequence of maps is a filtration of $B(R,P,L)(A)$ by subspaces. The
subspace $B(R,P,L)_s(A)$ consists of those points in $B(R,P,L)(A)$
that can be represented by trees with less than or equal to $s$ vertices.
\end{definition}

\begin{example}
The generalized $A$--labelled trees with no vertices (i.e. only a root
and some leaves) correspond one-to-one with (unordered) partitions of
$A$. We therefore see that
\[ B(R,P,L)_0 = R \circ L \]
where $\circ$ is the composition product of symmetric sequences.
\end{example}

\begin{example}
Take $R = L = I$ so that $B(R,P,L) = B(P)$. We then have $B(P)_0 = I$ by
the previous example. If $|A| > 1$ there is precisely one
(non-generalized) $A$--labelled tree with only one vertex and we
therefore get
\[ B(P)_1(A) = \begin{cases} S^1 \otimes P(A) & \text{if $|A| > 1$}; \\
B(P)(1) \isom S & \text{if $|A| = 1$}; \end{cases} \]
where $S$ is the unit of the symmetric monoidal category $\cat{C}$.
\end{example}

We can think of the sequence
\[ B(R,P,L)_0(A) \to B(R,P,L)_1(A) \to \dots \to B(R,P,L)(A) \]
as a kind of `cellular' filtration. That is, we obtain $B(R,P,L)_s(A)$ by
attaching `cells' to $B(R,P,L)_{s-1}(A)$, one for each generalized
$A$--labelled tree $T$ with exactly $s$ vertices. The following
proposition makes this precise.

\begin{proposition} \label{prop:cells}
There is a pushout square in $\cat{C}$ of the form
\[ \begin{diagram}
    \node{\bigvee_{T \in \mathsf{Q}_s(A)} \partial w(T)_+ \otimes
    (R,P,L)_A(T)} \arrow{s} \arrow{e} \node{B(R,P,L)_{s-1}(A)} \arrow{s}
    \\
    \node{\bigvee_{T \in \mathsf{Q}_s(A)} w(T)_+ \otimes (R,P,L)_A(T)}
    \arrow{e} \node{B(R,P,L)_s(A)}
\end{diagram} \]
where $\partial w(T)$ denotes the boundary of the space $w(T)$.
\end{proposition}

To identify the top horizontal map in this diagram we use the following
simple but important lemma.

\begin{lemma} \label{lem:boundary}
Let $T$ be a generalized $A$--labelled tree. Then
\[ \partial w(T)_+ \isom \colim_{U < T} w(U)_+. \]
The indexing category of the colimit is the full subcategory of $U \in
\mathsf{Tree}(A)$ with $U < T$.
\end{lemma}
\begin{proof}[Proof of Lemma]
This is a categorical reflection of that fact (Lemma \ref{lem:gen_W})
that the boundary $\partial w(T)$ consists precisely of those weightings
of $T$ in which some edge has length zero.
\end{proof}
\begin{proof}[Proof of Proposition \ref{prop:cells}]
The top horizontal map in the diagram is given by
\begin{eqnarray*} \bigvee_{T \in \mathsf{Q}_s(A)} \partial w(T)_+ \otimes
(R,P,L)_A(T) & \isom & \bigvee_{T \in \mathsf{Q}_s(A)} \colim_{U < T}
\left[w(U)_+ \otimes (R,P,L)_A(T)\right] \\
        & \longrightarrow & \bigvee_{T \in \mathsf{Q}_s(A)} \colim_{U <
        T} \left[w(U)_+ \otimes (R,P,L)_A(U)\right] \\
        & \longrightarrow & B(R,P,L)_{s-1}(A).
\end{eqnarray*}
Here we've used the fact that $- \otimes C$ is a left adjoint so
commutes with colimits. If $T \in \mathsf{Q}_s(A)$ and $U < T$ then $U
\in \mathsf{Tree}_{s-1}(A)$ so there are compatible maps from $w(U)_+
\otimes (R,P,L)_A(U)$ to the coend defining $B(R,P,L)_{s-1}(A)$.

With this definition, it is a simple exercise in naturality and colimits
to see that the square commutes. To see that it is a pushout, take a
commutative diagram
\[ \begin{diagram}
    \node{\bigvee_{T \in \mathsf{Q}_s(A)} \partial w(T)_+ \otimes
    (R,P,L)_A(T)} \arrow{s} \arrow{e} \node{B(R,P,L)_{s-1}(A)} \arrow{s}
    \\
    \node{\bigvee_{T \in \mathsf{Q}_s(A)} w(T)_+ \otimes (R,P,L)_A(T)}
    \arrow{e} \node{X}
\end{diagram} \tag{$*$} \]
We have to show that this factors via a unique map
\[ B(R,P,L)_s(A) \to X. \]
Since $B(R,P,L)_s(A)$ is a coend and hence a colimit, it is enough to
get a unique set of compatible maps
\[ w(U)_+ \otimes (R,P,L)_A(T) \to X \]
for $U \leq T$ in $\mathsf{Tree}_s(A)$. If $T \notin \mathsf{Q}_s(A)$
the required map comes from the right-hand edge of diagram $(*)$. So
suppose that $T \in \mathsf{Q}_s(A)$. Then we have
\[ w(U)_+ \otimes (R,P,L)_A(T) \to w(T)_+ \otimes (R,P,L)_A(T) \to X \]
where the second map comes from the bottom edge of diagram $(*)$. We leave
the reader to check that these maps are compatible in the appropriate
way and suitably unique. We conclude that $B(R,P,L)_s(A)$ is the claimed
pushout.
\end{proof}

We use this result to identify the quotients of our filtration of the
bar construction.

\begin{corollary} \label{cor:pushout}
Let $P$ be a reduced operad in $\cat{C}$ with right module $R$ and left
module $L$. The following is a pushout square in $\cat{C}$:
\[ \begin{diagram}
    \node{B(R,P,L)_{s-1}(A)} \arrow{s} \arrow{e} \node{B(R,P,L)_s(A)}
    \arrow{s} \\
    \node{\basept} \arrow{e} \node{\bigvee_{T \in \mathsf{Q}_s(A)}
    w(T)_+/\partial w(T)_+ \otimes (R,P,L)_A(T)}
\end{diagram} \]
\end{corollary}
\begin{proof}
Since $- \otimes C$ preserves colimits, the following is a pushout
square in $\cat{C}$:
\[ \begin{diagram}
    \node{\partial w(T)_+ \otimes (R,P,L)_A(T)} \arrow{s} \arrow{e}
    \node{w(T)_+ \otimes (R,P,L)_A(T)} \arrow{s} \\
    \node{\basept} \arrow{e} \node{w(T)_+/\partial w(T)_+ \otimes
    (R,P,L)_A(T)}
\end{diagram} \]
The corollary now follows from Proposition \ref{prop:cells} and the
universal properties of colimits.
\end{proof}

\begin{remark} \label{rem:w_sphere}
Recall from Lemma \ref{lem:gen_W} that for any generalized $A$--labelled
tree $T$ with $s$ vertices, $w(T) \isom D^s$. Therefore, $w(T)_+/\partial
w(T)_+ \isom S^s$. We will be talking a lot about these spaces in the
coming sections, so we will give them some more compact notation:
\[ \underline{w}(T) := w(T)_+/\partial w(T)_+ \isom w(T)/\partial w(T)
\isom S^s \]
\end{remark}

The results for the cobar construction are, as usual, just the duals of
those for the bar construction. We summarize these briefly.

\begin{definition}[(Cofiltration of the cobar construction)]
\label{def:cofilter}
Let $Q$ be a reduced cooperad in a symmetric monoidal $\based$--category
$\cat{C}$ with right comodule $R$ and left comodule $L$. Then the
two-sided cobar construction $\Omega(R,Q,L)$ has a `cofiltration',
that is, there is a sequence
\[ \Omega(R,Q,L)(A) \to \dots \to \Omega(R,Q,L)^s(A) \to
\Omega(R,Q,L)^{s-1}(A) \to \cdots \]
where
\[ \Omega(R,Q,L)^s(A) := \int_{T \in \mathsf{Tree}_s(A)}
\Map_{\cat{C}}(w(T)_+,(R,Q,L)_A(T)), \]
and the `projection' map
\[ \Omega(R,Q,L)^s(A) \to \Omega(R,Q,L)^{s-1}(A) \]
comes from the inclusion of categories $\mathsf{Tree}_{s-1}(A) \to
\mathsf{Tree}_s(A)$ for $s \geq 1$.
\end{definition}
\begin{corollary}
With $Q,R,L$ as in Definition \ref{def:cofilter}, the following is a
pullback square:
\[ \begin{diagram}
    \node{\Omega(R,Q,L)^s(A)} \arrow{s} \arrow{e} \node{\prod_{T \in
    \mathsf{Q}_s(A)} \Map_{\cat{C}}(w(T)_+,(R,Q,L)_A(T))} \arrow{s} \\
    \node{\Omega(R,Q,L)^{s-1}(A)} \arrow{e} \node{\prod_{T \in
    \mathsf{Q}_s(A)} \Map_{\cat{C}}(\partial w(T)_+,(R,Q,L)_A(T))}
   \end{diagram} \]
We can identify the fibres of the projections by the pullback
squares
\[ \begin{diagram}
    \node{\prod_{T \in \mathsf{Q}_s(A)}
    \Map_{\cat{C}}(\underline{w}(T),(R,Q,L)_A(T))} \arrow{e} \arrow{s}
    \node{\Omega(R,Q,L)^s(A)} \arrow{s} \\
    \node{\basept} \arrow{e} \node{\Omega(R,Q,L)^{s-1}(A),}
   \end{diagram} \]
where $\underline{w}(T) = w(T)/\partial w(T)$. \qed
\end{corollary}

\subsection{Conditions for the inclusion maps of the filtration to be
cofibrations} \label{sec:cofibrations}

In the case that $\cat{C}$ is either $\based$ or $\spectra$, the
filtration of Section~\ref{sec:filter} allows us to construct a spectral
sequence converging to the homology of $B(R,P,L)$. The $E^1$ term of this
spectral sequence is given by the homologies of the homotopy cofibres of
the inclusion maps of the filtration. In this section we give conditions
under which these inclusions are cofibrations (in the standard model
category structures on $\based$ and $\spectra$) and which therefore
ensure that the homotopy cofibres are given by the strict cofibres,
or filtration quotients, that we have already calculated.

We state the main result of this section (Proposition
\ref{prop:cofibrations} below) for a general symmetric monoidal
$\based$--category $\cat{C}$ with a compatible model structure. Definition
\ref{def:model} says what we mean by `compatible' here. We use Mark
Hovey's book \cite{hovey:1999} as our basic reference for model
categories.

\begin{definition} \label{def:model}
A \emph{symmetric monoidal $\based$--model category} is a symmetric
monoidal $\based$--category $\cat{C}$ (as in Section~\ref{sec:monoidal})
together with a model structure (in the sense of
\cite[Definition~1.1.3]{hovey:1999}) such that the tensoring makes
$\cat{C}$ into a
$\based$--model category (in the sense of
\cite[Definition~4.2.18]{hovey:1999}). That is, if $X \to Y$ is a cofibration in
$\based$ and $C \to D$ is a cofibration in $\cat{C}$ then the induced map
\[ (X \otimes D) \amalg_{X \otimes C} (Y \otimes C) \to Y \otimes D \]
is a cofibration in $\cat{C}$ that is trivial if either of our original
cofibrations is. (The domain of this map is the pushout of $X \otimes D$
and $Y \otimes C$ over $X \otimes C$.)
\end{definition}

\begin{remark} \label{rem:model}
We should say a few words about this definition. Firstly, we are
\emph{not} requiring that $\cat{C}$ be a monoidal model category in
its own right (in the sense of \cite[Section~4.2.6]{hovey:1999}). That is, we
are not insisting that the symmetric monoidal structure $\barwedge$
on $\cat{C}$ in any way respect the model structure. Our reason for
doing this is to preserve the self-duality of Definition \ref{def:model}
(see Lemma \ref{lem:model_dual} below). In general, the opposite category
of a monoidal model category is not another monoidal model category and
we wish to dualize our theory to obtain results on the cobar construction.

On the other hand, the hypotheses we need to prove Proposition
\ref{prop:cofibrations} are natural consequences of the assumption
that $\cat{C}$ \emph{is} a monoidal model category, \emph{and} the
categories we are most interested in, $\based$ and $\spectra$, satisfy
this assumption. This suggests that a breaking of the symmetry between
bar and cobar is necessary when we come to study the homotopy theory
of these constructions. In this paper, we do not pretend to give the
beginnings of such a theory and, in particular, we do not claim that
Definition \ref{def:model} is the philosophically correct way to mix
model category theory into this paper. For us, it serves the purposes
of allowing us to make calculations with our spectral sequence in cases
that are of interest.
\end{remark}

\begin{lemma} \label{lem:cofibration}
Let $\cat{C}$ be a symmetric monoidal $\based$--model category. If $C
\in \cat{C}$ is cofibrant and $X \to Y$ is a cofibration in $\based$
then $X \otimes C \to Y \otimes C$ is a cofibration in $\cat{C}$.
\end{lemma}
\begin{proof}
Apply the definition of $\based$--model category to the cofibrations $X
\to Y$ and $\basept \to C$.
\end{proof}

\begin{prop} \label{prop:cofibrations}
Let $\cat{C}$ be a symmetric monoidal $\based$--model category such that
if $C,D$ are cofibrant then $C \barwedge D$ is also cofibrant. Let $P$
be a reduced operad in $\cat{C}$ with right module $R$ and left module $L$
such that, for all $A$, the objects $P(A),R(A),L(A)$ are cofibrant. Then,
for all $s \geq 1$ and all finite sets $A$, the map
\[ B(R,P,L)_{s-1}(A) \to B(R,P,L)_s(A) \]
of Definition \ref{def:filter} is a cofibration in $\cat{C}$.
\end{prop}
\begin{proof}
The cofibrancy conditions on the $P(A),L(A),R(A)$ together with the
extra condition on $\cat{C}$ ensure that the objects $(R,P,L)_A(T)$
are all cofibrant. For any generalized tree $T$, the map
\[ \partial w(T)_+ \to w(T)_+ \]
is a cofibration in $\based$ (it is the inclusion of the boundary of a
ball). Therefore, by Lemma \ref{lem:cofibration},
\[ \partial w(T)_+ \otimes (R,P,L)_A(T) \to w(T)_+ \otimes (R,P,L)_A(T) \]
is a cofibration. Proposition \ref{prop:cells} tells us that the
filtration map
\[ B(R,P,L)_{s-1}(A) \to B(R,P,L)_s(A) \]
is a pushout of a coproduct of such maps so it too is a cofibration.
\end{proof}

\begin{remark}
As we commented in Remark \ref{rem:model} above, if $\cat{C}$ is a
symmetric monoidal model category in its own right, we get for free that
$C$ and $D$ cofibrant imply $C \barwedge D$ cofibrant. In particular this
is the case for $\based$ and $\spectra$ (that is, the $S$--modules of
EKMM \cite{elmendorf/kriz/mandell/may:1997}).
\end{remark}

As promised, our definition of symmetric monoidal $\based$--model category
is self-dual.

\begin{lemma} \label{lem:model_dual}
Let $\cat{C}$ be a symmetric monoidal $\based$--model category. Then
$\cat{C}^\text{op}$ is also a symmetric monoidal $\based$--model category
with the standard dual symmetric monoidal and model structures.
\end{lemma}
\begin{proof}
We already know from Proposition \ref{prop:dual} that $\cat{C}^\text{op}$ is a
symmetric monoidal $\based$--category. Recall that the tensoring for
$\cat{C}^\text{op}$ is given by the cotensoring for $\cat{C}$, the cofibrations
in $\cat{C}^\text{op}$ are the fibrations in $\cat{C}$ and a pushout in
$\cat{C}^\text{op}$ is a pullback in $\cat{C}$. The weak equivalences in
$\cat{C}^\text{op}$ are the same as those in $\cat{C}$.

To see that $\cat{C}^\text{op}$ is a $\based$--model category we have to show
that if $X \to Y$ is a cofibration in $\based$ and $D \to C$ a fibration
in $\cat{C}$ then
\[ \Map_{\cat{C}}(Y,D) \to \Map_{\cat{C}}(Y,C)
\times_{\Map_{\cat{C}}(X,C)} \Map_{\cat{C}}(X,D) \]
is a fibration in $\cat{C}^\text{op}$ that is trivial if either of our
original maps is a weak equivalence. This result is given by Lemma 4.2.2
of \cite{hovey:1999}.
\end{proof}

The result dual to Proposition \ref{prop:cofibrations} is then the
following.

\begin{cor} \label{cor:fibrations}
Let $\cat{C}$ be a symmetric monoidal $\based$--model category such that
if $C,D$ are fibrant then $C \barwedge D$ is also fibrant. Let $Q$ be a
reduced cooperad in $\cat{C}$ with right comodule $R$ and left comodule
$L$ such that all the objects $Q(A),R(A),L(A)$ are fibrant. Then the map
\[ \Omega(R,Q,L)^s(A) \to \Omega(R,Q,L)^{s-1}(A) \]
of Definition \ref{def:cofilter} is a fibration in $\cat{C}$. \qed
\end{cor}
In these circumstances, then, the fibres of the maps in the tower for
$\Omega(R,Q,L)$ are also the homotopy fibres and so can be used to
calculate the $E^1$ term of the associated spectral sequence.

\begin{remark}
In our categories of interest, $\based$ and $\spectra$, all objects
are fibrant and so the conditions of Corollary \ref{cor:fibrations}
hold for any cooperad and any comodules over it.
\end{remark}

\subsection{The algebraic bar and cobar constructions} \label{sec:alg_bar}

So far we have constructed (under suitable conditions) a filtration of the
two-sided bar construction $B(R,P,L)$ by a sequence of cofibrations. This
filtration yields a homology spectral sequence whose $E^1$ term turns out
to be given by an algebraic version of our bar construction. In fact, it
was this algebraic version, previously studied by Ginzburg--Kapranov
\cite{ginzburg/kapranov:1994}, Getzler--Jones \cite{getzler/jones:1994}
and Fresse \cite{fresse:2004} among others, that inspired our definition
of the bar construction for operads in topological settings. This section
is devoted to the description of this algebraic bar construction. As
in the topological case, we will only deal with \emph{reduced} operads,
that is, those for the unit map $k \to P(1)$ is an isomorphism.

Our definition of the algebraic bar construction emphasizes its similarity
to the topological versions of Section~\ref{sec:bar} and
Section~\ref{sec:bar(modules)} and it will follow the same pattern.

\begin{definition}
Let $P$ be a reduced operad in the category $\mathsf{Ch}_k$ of chain
complexes over the commutative ring $k$ (with the symmetric monoidal
structure of Definition \ref{def:ch_k}). Let $R$ be a right $P$--module
and $L$ a left $P$--module. More or less repeating Definition
\ref{def:(R,P,L)}, we define a functor
\[ (R,P,L)_A \co \mathsf{Tree}(A)^\text{op} \to \mathsf{Ch}_k \]
for each nonempty finite set $A$ by the formula
\[ (R,P,L)_A(T) := R(i(r)) \otimes \bigotimes_{\text{vertices $v \in T$}}
P(i(v)) \otimes \bigotimes_{\text{leaves $l \in T$}} L(\iota^{-1}l). \]
The composition maps for $R$, $P$ and $L$ make $(R,P,L)_A$ into a functor
as claimed. In making explicit calculations we have to be careful with
the signs involved in the symmetry isomorphism for $\otimes$ but for
theoretical purposes we can treat $(R,P,L)_A(T)$ as an unordered tensor
product (see Remark \ref{rem:symmetry}).
\end{definition}

We now wish to define the bar construction $B(R,P,L)$ by the same coend
formula as in Definition \ref{def:two-sided_bar}. For this we need chain
complex versions of the spaces $w(T)$ of weightings on trees $T \in
\mathsf{Tree}(A)$. As in the topological case, the structures of these
`spaces', and how they fit together for different trees, are the key
parts of the definition of the bar construction.

\begin{figure}[htbp]
\begin{center}
\input{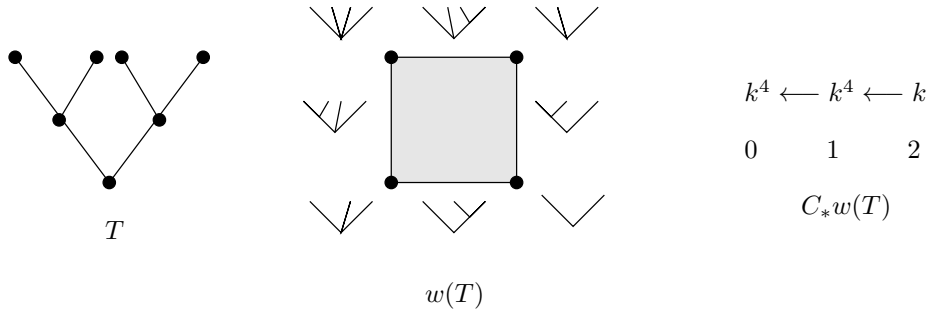}
\caption{Example of the chain complex $C_*w(T)$ showing the cellular
decomposition of $w(T)$}
\label{fig:chains}
\end{center}
\end{figure}

\begin{definition}
Let $T$ be a generalized $A$--labelled tree. The chain complex $C_*w(T)$
representing the space of weightings on $T$ will be the cellular chain
complex for a certain cellular decomposition of the space $w(T)$. The
cells in this decomposition correspond one-to-one with the trees $U \in
\mathsf{Tree}(A)$ with $U \leq T$. The $r$--skeleton of $w(T)$ is given by
%
%
%
\[ \operatorname{sk}_r w(T) := \colim_{U < T : \; U \in
\mathsf{Tree}_r(A)}w(U). \]
The attaching map for the cell corresponding to the tree $U$ with $r+1$
vertices is the map
\[ S^r \isom \partial w(U) \isom \colim_{V < U} w(V) \to \colim_{V <
T : \; V \in \mathsf{Tree}_r(A)} w(V) = \operatorname{sk}_r w(T). \]
The cellular chain complex for this cell structure then has
\[ C_r w(T) = \bigoplus_{U \leq T : \; U \in \mathsf{Q}_r(A)}
H_r(w(U),\partial w(U)) \isom \bigoplus_{U \leq T : \; U \in
\mathsf{Q}_r(A)} \widetilde{H}_r(\underline{w}(U)). \]
Recall from Remark \ref{rem:w_sphere} that $\underline{w}(U)$ denotes
the quotient $w(U)/\partial w(U)$. The differential
\[ C_r w(T) \to C_{r-1} w(T) \]
is given by summing the maps\footnote{The last part of this composite
comes from the map
\[ \partial w(U)_+ \isom \colim_{V < U} w(V)_+ \to w(V)_+/\partial w(V)_+
= \underline{w}(V) \]
given by collapsing to the basepoint everything except the interior of
the `face' $w(V)$ of $\partial w(U)$.}
\[ \widetilde{H}_r(\underline{w}(U)) \isom H_r(w(U),\partial w(U)) \to
\widetilde{H}_{r-1}(\partial w(U)_+) \to \widetilde{H}_{r-1}(\underline{w}(V)) \]
for pairs $(U,V)$ with $V < U$, $|U| = r$ and $|V| = r-1$. An example
of this chain complex for a particular tree is shown in Figure
\ref{fig:chains}.

The inclusion $w(U) \to w(T)$ is cellular and so we have inclusions
\[ C_*w(U) \to C_*w(T) \]
for $U < T$. These make $C_*w(-)$ into a functor
\[ C_*w(-) \co \mathsf{Tree}(A) \to \mathsf{Ch}_k. \]
This is the chain complex analogue of the functor $w(-)$ of Definition
\ref{def:gen_w_functor}.
\end{definition}

\begin{definition} \label{def:alg_bar}
With our `chain complexes of weighted trees' $C_*w(T)$, we now define the
\emph{two-sided algebraic bar construction} on the reduced operad $P$ with
coefficients in $R$ and $L$ to be the symmetric sequence $B(R,P,L)$ with
\[ B(R,P,L)(A) := \int^{T \in \mathsf{Tree}(A)} C_*w(T) \otimes
(R,P,L)_A(T). \]
This coend is calculated in the category of chain complexes on $k$ and
results in a chain complex $B(R,P,L)(A)$. However, it will be useful to
consider a bicomplex structure on $B(R,P,L)(A)$ for which this chain
complex is the total complex. The bicomplex structure comes about by
considering the tensor product of the chain complexes $C_*w(T)$ and
$(R,P,L)_A(T)$ as a bicomplex with gradings and differentials coming
from these separate terms. We will write
\[ B(R,P,L)_{*,*}(A) \]
to emphasize this bigrading with the first index denoting the grading
that comes from $C_*w(T)$ (we'll call this the \emph{tree grading}) and
the second the grading that comes from $(R,P,L)_A(T)$ (which we'll call
the \emph{internal grading}). We then have two separate differentials
on $B(R,P,L)_{*,*}$:
\[ \partial \co B(R,P,L)_{*,*} \to B(R,P,L)_{*-1,*} \]
coming from the differentials on the chain complexes $C_*w(T)$ which
will refer to as the \emph{tree differential} on the bar construction, and
\[ d\co B(R,P,L)_{*,*} \to B(R,P,L)_{*,*-1} \]
coming from the differentials on the $(R,P,L)_A(T)$ which we will call
the \emph{internal differential}.

In later sections, we will be applying the algebraic bar construction to
operads of graded $k$--modules, that is, chain complexes with zero
differential. In this case, the internal differential of $B(R,P,L)(A)$
will be zero.
\end{definition}

We can give a more explicit description of $B(R,P,L)$ as follows.

\begin{lem} \label{lem:explicit_bar}
Let $P$ be a reduced operad in $\mathsf{Ch}_k$ with right module $R$ and
left module $L$. Then we have\footnote{Here, as elsewhere, the reduced
homology of the quotient $\underline{w}(T) = w(T)/\partial w(T)$ can
be replaced with the homology of the pair $(w(T),\partial w(T))$. Both
of these are isomorphic to the graded module $k$ concentrated in degree
$|T|$.}
\[ B(R,P,L)_{s,*}(A) \isom \bigoplus_{T \in \mathsf{Q}_s(A)}
\widetilde{H}_s(\underline{w}(T)) \otimes (R,P,L)_A(T) \]
as chain complexes of $k$--modules with respect to the internal grading
and differential.

Under these isomorphisms, the tree differential
\[ \partial\co B(R,P,L)_{s,*}(A) \to B(R,P,L)_{s-1,*}(A) \]
is given by summing, over all pairs $(T,U)$ with $U < T$, $|T| = s$
and $|U| = s-1$, the maps
\[ \widetilde{H}_s(\underline{w}(T)) \otimes (R,P,L)_A(T) \to
\widetilde{H}_{s-1}(\underline{w}(U)) \otimes (R,P,L)_A(U) \]
obtained by combining the maps
\[ (R,P,L)_A(T) \to (R,P,L)_A(U) \]
with the terms
\[ \widetilde{H}_s(\underline{w}(T)) \to \widetilde{H}_{s-1}(\underline{w}(U)) \]
from the top differential of the chain complex $C_*w(T)$.
\end{lem}
\begin{proof}
We consider a filtration of the algebraic bar construction analogous to
that of Section~\ref{sec:filter} for the topological version. Virtually
the same analysis applies and we get short exact sequences of chain
complexes\footnote{The notation here is probably rather confusing. We
are using $B(R,P,L)_s(A)$ to denote the part of the filtration of
$B(R,P,L)(A)$ obtained via the chain complex version of Definition
\ref{def:filter}. This is not to be confused with $B(R,P,L)_{s,*}$ which
is the graded summand of tree degree $s$. In fact, it's a consequence
of the proof of this lemma that
\[ B(R,P,L)_s(A) \isom \bigoplus_{r \leq s} B(R,P,L)_{r,*}(A). \]}
\[ B(R,P,L)_{s-1}(A) \rightarrow B(R,P,L)_s(A) \rightarrow \bigoplus_{T
\in \mathsf{Q}_s(A)} C_*w/C_*\partial w(T) \otimes (R,P,L)_T(A). \]
where $C_*\partial w(T)$ is the cellular chain complex for the subcomplex
$\partial w(T) \subset w(T)$ (that is, everything except the top-dimension
cell).
Notice that
\[ C_*w(T)/C_*\partial w(T) \isom \widetilde{H}_s(\underline{w}(T)) \]
for $T \in \mathsf{Q}_s(A)$. We construct a splitting of the above short
exact sequence (with respect to the internal differential but \emph{not}
the tree differential) using the obvious splittings (as $k$--modules)
of the sequences
\[ C_*\partial w(T) \to C_*w(T) \to \widetilde{H}_*(\underline{w}(T)). \]
We get by induction on $s$ that
\[ B(R,P,L)(A) \isom \bigoplus_{T \in \mathsf{Tree}(A)}
\widetilde{H}_{|T|}(\underline{w}(T)) \otimes (R,P,L)_A(T) \]
which splits, by tree degree, into the isomorphisms of the lemma. We leave
the reader to check that the tree differential has the promised formula.
\end{proof}

\begin{remark} \label{rem:fresse}
Choosing generators of the groups $H_{|T|}(w(T),\partial w(T)) \isom k$
determines an isomorphism
\[ B(R,P,L)(A) \isom \bigoplus_{T \in \mathsf{Tree}(A)} (R,P,L)_A(T) \]
which is the definition of the algebraic bar construction given by Fresse
in \cite[Section~4.4]{fresse:2004}. Such choices determine choices of the
coefficients (in fact, signs) for the maps that make up the differential
$\partial$ on $B(R,P,L)(A)$.

Fresse shows that this bar construction is a representative of the derived
composition product of $R$ and $L$ as $P$--modules, that is,
\[ B(R,P,L) \simeq R \circ^{\mathbb{L}}_{P} L \]
and so the homology groups of $B(R,P,L)$, with respect to the tree
differential, are $\operatorname{Tor}$ groups of $P$--modules.
\end{remark}

The relationship between this algebraic bar construction and the
simplicial bar construction was analyzed by Fresse. His proof of the
following proposition uses a `levelization' process analogous to that
we used in the proof of Proposition \ref{prop:simp}.

\begin{prop} {\rm(\cite[Theorem~4.1.8]{fresse:2004})}\qua
The algebraic two-sided bar construction $B(R,P,L)$ is quasi-isomorphic
to the normalized chain complex of the simplicial bar construction on $P$
with coefficients in $R$ and $L$ (the algebraic version of Definition
\ref{def:two-sided_simp}).
\end{prop}

As usual, we have the dual constructions and results.

\begin{definition}
Let $Q$ be a reduced cooperad of chain complexes of $k$--modules
with right comodule $R$ and left comodule $L$. Then there are functors
$(R,Q,L)_A$ from $\mathsf{Tree}(A)$ to $\mathsf{Ch}_k$ and we can define
the \emph{algebraic cobar construction on $Q$ with coefficients in $R$
and $L$} by the same formula
\[ \Omega(R,Q,L)(A) := \int_{T \in \mathsf{Tree}(A)}
\Hom(C_*w(T),(R,Q,L)_A(T)) \]
as in Definition \ref{def:cobar(cooperad)}, where, for chain complexes
$M,N$, $\Hom(M,N)$ denotes the chain complex of maps of graded modules $M
\to N$. The cobar construction is a bicomplex with an internal grading
and differential coming from the $(R,Q,L)_A(T)$ and a tree grading
and differential $\partial^*$ coming from the $C_*w(T)$. We follow the
convention that $\Hom(M,N)_{s,t} = \Hom(M_{-s},N_t)$ so that the tree
grading on the cobar construction is concentrated in negative degrees.
\end{definition}
There is an explicit description of the cobar construction analogous to
that of Lemma \ref{lem:explicit_bar} for the bar construction.

\begin{lem} \label{lem:explicit_cobar}
With $Q,R,L$ as above:
\[ \Omega(R,Q,L)_{-s,*}(A) := \bigoplus_{T \in \mathsf{Q}_s(A)}
\Hom(\widetilde{H}_s(\underline{w}(T)),(R,Q,L)_A(T)) \]
which again is just isomorphic to
\[ \bigoplus_{T \in \mathsf{Tree}(A)} (R,Q,L)_A(T) \]
after choosing generators of the groups
$\widetilde{H}_s(\underline{w}(T))$. The internal grading and differential
correspond in the obvious way under this isomorphism. The explicit form
of the tree differential $\partial^*$ is given by the maps
\[ (R,Q,L)_A(U) \to (R,Q,L)_A(T) \]
with coefficients again given by the components
\[ \widetilde{H}_s(\underline{w}(T)) \to \widetilde{H}_{s-1}(\underline{w}(U)) \]
of the top differential on the chain complex $C_*w(T)$. \qed
\end{lem}

\begin{definition}
When $P$ is a reduced operad in the category of graded $k$--modules,
the unit symmetric sequence $I$ defined by
\[ I(A) := \begin{cases} k & \text{if $|A| = 1$}; \\ 0 &
\text{otherwise}. \end{cases} \]
is both a left and right $P$--module. The \emph{reduced bar construction
on $P$} is then given by the two-sided bar construction with coefficients
in $I$ on both sides:
\[ B(P) := B(I,P,I) \]
The definition of the algebraic bar construction reduces in this case to
\[ B(P)(A) = \int^{T \in \mathsf{T}(A)} C_*\overline{w}(T) \otimes P_A(T). \]
Recall that the space $\overline{w}(T)$ is the quotient of $w(T)$ by the
subspace $w_0(T)$ of weightings that give length $0$ to either the root
edge or a leaf edge of $T$. This subspace is in fact a subcomplex with
respect to our chosen cellular structure on $w(T)$.\footnote{It is the
union of the cells corresponding to $U < T$ that are not in the original
category $\mathsf{T}(A)$, that is, that are generalized trees, but not
trees in the sense of Section~\ref{sec:trees}.} Therefore we
obtain a cell structure on $\overline{w}(T)$ and in the above formula,
$C_*\overline{w}(T)$ denotes the relative cellular chain complex for
the pair $(\overline{w}(T),\basept)$, or equivalently, for the pair
$(w(T),w_0(T))$.\footnote{The tensor product of chain complexes is
here playing the role of the smash product of based spaces so we need
the reduced version of the cellular chain complex. Strictly speaking,
the chain complex $C_*w(T)$ is the relative chain complex of the pair
$(w(T)_+,\basept)$.} It is clear that $C_*\overline{w}(T)$ is a quotient
of $C_*w(T)$.

It's also easy to check that by Lemma \ref{lem:explicit_bar} we have
\[ B(P)(A) \isom \bigoplus_{T \in \mathsf{T}(A)}
\widetilde{H}_{|T|}(\underline{w}(T)) \otimes P_A(T) \]
which (after choosing isomorphisms $\widetilde{H}_{|T|}(\underline{w}(T))
\isom k$) is the original definition of the algebraic bar construction
given in Getzler--Jones \cite[Section~2.1]{getzler/jones:1994}.

Similarly, if $Q$ is a reduced cooperad then $I$ is both a left and right
$Q$--comodule and the \emph{reduced cobar construction on $Q$} is
\[ \Omega(Q) := \Omega(I,Q,I) \]
and is given by a formula analogous to that of Definition
\ref{def:cobar(cooperad)}.
\end{definition}

As in the topological case, the reduced algebraic bar construction on a
reduced operad $P$ of chain complexes has a cooperad structure. We now
describe this.

\begin{definition} \label{def:alg_cooperad(bar)}
The required maps
\[ B(P)(A \cup_a B) \to B(P)(A) \otimes B(P)(B) \]
are defined in exactly the same way as the corresponding maps in the
topological case (Definition \ref{sec:cooperad}). To do this we must
construct the algebraic versions of the key maps (\ref{eq:key}):
\[ C_*\overline{w}(T \cup_a U) \to C_*\overline{w}(T) \otimes C_*\overline{w}(U) \]
for $A$--labelled trees $T$ and $B$--labelled trees $U$. We get this by
taking the map of cellular chain complexes induced by the topological map
\[ \overline{w}(T \cup_a U) \to \overline{w}(T) \smsh \overline{w}(U) \]
of Definition \ref{def:cooperad_maps}. For this to work, we need the
following lemma.

\begin{lem}
Let $T$ be an $A$--labelled tree, $U$ a $B$--labelled tree and let $a
\in A$. The map
\[ \overline{w}(T \cup_a U) \to \overline{w}(T) \smsh \overline{w}(U) \]
is cellular, that is, it preserves skeleta.
\end{lem}
\begin{proof}
A point $p$ in $\overline{w}(T \cup_a U)$ is in the $s$--skeleton if and only
if it is in the subspace $\overline{w}(V)$ for some tree $V$ with $s$ vertices.
If this tree $V$ is not of type $(A,B)$ then $p$ is mapped to the
basepoint which is certainly in the $s$--skeleton of the right-hand side.
If $V$ is of type $(A,B)$ (that is, obtained by grafting an $A$--labelled
tree $T'$ and a $B$--labelled tree $U'$) then the point $p$ maps to a pair
consisting of a point in some $\overline{w}(T') \subset \overline{w}(T)$ and a point
in some $\overline{w}(U') \subset \overline{w}(U)$. The first point is in the
$s'$--skeleton of $\overline{w}(T)$ where $T'$ has $s'$ vertices. The second
point is in the $s''$--skeleton of $\overline{w}(U)$ where $U'$ has $s''$
vertices. Therefore the pair is in the $s'+s''$--skeleton of $\overline{w}(T)
\smsh \overline{w}(U)$. However, since $V$ only had $s$ vertices, we must have
$s'+s'' \leq s$. So the image of $p$ is in the $s$--skeleton of
$\overline{w}(T) \smsh \overline{w}(U)$ as required.
\end{proof}

It is easy to describe explicitly the resulting map of chain complexes
\[ C_*\overline{w}(T \cup_a U) \to C_*\overline{w}(T) \otimes C_*\overline{w}(U). \]
Recall that the left-hand side is given by the direct sum over $V \in
\mathsf{Tree}(A \cup_a B)$ with $V \leq T \cup_a U$ of the homology
groups $\widetilde{H}_*(\underline{w}(V))$. The above map sends the term
corresponding to a tree $V$ that is not of type $(A,B)$, to zero. If $V
= T' \cup_a U'$, then $T' \leq T$ and $U' \leq U$ and the corresponding
term maps to the right-hand side via the isomorphism
\[ \widetilde{H}_*(\underline{w}(T' \cup_a U')) \to
\widetilde{H}_*(\underline{w}(T')) \otimes
\widetilde{H}_*(\underline{w}(U')), \]
which is induced by the homeomorphism\footnote{It is easy to check that
this map is a bijection. The spaces involved are all spheres which are
compact Hausdorff, so it is a homeomorphism.}
\[ \underline{w}(T' \cup_a U') \to \underline{w}(T') \smsh
\underline{w}(U'), \]
which in turn is a quotient of the map
\[ \overline{w}(T' \cup_a U') \to \overline{w}(T') \smsh \overline{w}(U'). \]

With this key map in place, the rest of the formal definition of the
cooperad structure maps for the topological bar construction (Definition
\ref{def:formal_cooperad_maps}) carries over to the algebraic case.
\end{definition}

\begin{lem} \label{lem:explicit_cooperad}
Let $P$ be a reduced operad in $\mathsf{Ch}_k$. Under the isomorphism
of Lemma \ref{lem:explicit_bar}, the cooperad structure on $B(P)$
corresponds to the cooperad structure on the chain complexes
\[ \bigoplus_{T \in \mathsf{T}(A)} \widetilde{H}_{|T|}(\underline{w}(T))
\otimes P_A(T) \]
whose cocomposition maps are given by summing over the maps obtained by
combining the isomorphisms
\[ H_*(\underline{w}(T \cup_a U)) \to H_*(\underline{w}(T)) \otimes
H_*(\underline{w}(U)) \]
with the isomorphisms
\[ P_{A \cup_a B}(T \cup_a U) \to P_A(T) \otimes P_B(U). \]
\end{lem}
\begin{proof}
This is a simple check using the definition of the isomorphism in Lemma
\ref{lem:explicit_bar} by splittings of short exact sequences.
\end{proof}

\begin{remark}
Choosing generators for the groups $\widetilde{H}_*(\underline{w}(T))$, we see
that this is equivalent to the cooperad structure defined by
Getzler--Jones \cite{getzler/jones:1994} and by Fresse \cite{fresse:2004}.
\end{remark}

\begin{definition}
Dually, if $Q$ is a reduced cooperad of chain complexes, there
is an operad structure on the reduced algebraic cobar construction
$\Omega(Q)$. The corresponding operad structure under the isomorphism
of Lemma \ref{lem:explicit_cobar} is built from the isomorphisms
\[ Q_A(T) \otimes Q_B(U) \to Q_{A \cup_a B}(T \cup_a U) \]
and the same maps
\[ \widetilde{H}_*(\underline{w}(T \cup_a U)) \to
\widetilde{H}_*(\underline{w}(T)) \otimes \widetilde{H}_*(\underline{w}(U)). \]
\end{definition}

\begin{remark} \label{rem:extend}
It does not take much more effort to extend the cooperad and operad
structure above to maps
\[ B(R,P,L) \to B(R,P,I) \varcirc B(I,P,L) \]
and
\[ \Omega(R,Q,I) \circ \Omega(I,Q,L) \to \Omega(R,Q,L) \]
following the same sort of generalization that we did in
Section~\ref{sec:bar(modules)_maps}.
\end{remark}

\subsection{A spectral sequence for the homology of the bar construction}
\label{sec:specseq}
We now turn our attention directly to the homology spectral sequences
born from the filtration of the bar construction and cofiltration of
the cobar construction.\footnote{If $\cat{C}$ is the category of based
spaces, we only get a spectral sequence for the bar construction and not
for the cobar construction. This is because a fibre sequence in $\based$
does not immediately yield a long exact sequence in homology.} The work
we have done in the last few sections allows us to identify the $E^1$
terms of these spectral sequences, under suitable conditions, with the
algebraic bar and cobar constructions.

A quick word on notation: from now on, the only topological categories
$\cat{C}$ we are interested in are $\based$ and $\spectra$. We will
therefore drop the notation $\barwedge$ for the monoidal product and
$\otimes$ for the tensoring over $\based$, replacing both with the
standard notation $\smsh$. We will reserve $\otimes$ for the tensor
product of graded $k$--modules.

\begin{prop} \label{prop:specseq}
Let $P$ be a reduced operad in $\based$ or $\spectra$ with right module
$R$ and left module $L$ such that all the objects $P(A),R(A),L(A)$ are
cofibrant and all homology groups $H_*P,H_*R,H_*L$ flat $k$--modules. Then
for each finite set $A$ there is a spectral sequence converging to
$H_*B(R,P,L)(A)$ with $E^1$--term and first differential given by the
algebraic bar construction:
\[ (E^1,d^1) \isom (B(H_*R,H_*P,H_*L)(A),\partial) \implies
H_*B(R,P,L)(A). \]
Let $Q$ be a reduced cooperad in $\spectra$ with right comodule $R$ and
left comodule $L$ such that all the objects $Q(A),R(A),L(A)$ are
fibrant\footnote{This is really automatic since we have chosen $\spectra$
to be the category of $S$--modules of EKMM
\cite{elmendorf/kriz/mandell/may:1997} in which all objects are fibrant.
If we want to work with other categories of spectra, however, we need this
condition.} and all the homology groups $H_*Q,H_*R,H_*L$ are flat
$k$--modules. Then for each finite set $A$ there is a spectral sequence
converging to $H_*(\Omega(R,Q,L)(A))$ with $E^1$--term and first
differential given by the algebraic cobar construction:
\[ (E^1,d^1) \isom (\Omega(H_*R,H_*Q,H_*L)(A),\partial^*) \implies
H_*\Omega(R,Q,L)(A). \]
\end{prop}
\begin{remark}
By the comments of Remark \ref{rem:fresse}, the work of Fresse allows
us to identify the $E^2$ terms of these spectral sequences as suitable
$\operatorname{Tor}$ groups. That is, our spectral sequence take the form
\[ E^2 = \operatorname{Tor}^{H_*P}(H_*R,H_*L) \implies H_*B(R,P,L) \]
and
\[ E^2 = \operatorname{Tor}^{H_*Q}(H_*R,H_*L) \implies
H_*\Omega(R,Q,L). \]
This suggests that the topological bar and cobar constructions should have
an interpretation as topological $\operatorname{Tor}$ objects. We have
not yet studied the homotopy theory of these constructions sufficiently
to make this precise.
\end{remark}

\begin{proof}[Proof of Proposition \ref{prop:specseq}]
By Proposition \ref{prop:cofibrations}, we have cofibre sequences
\[ B(R,P,L)_{s-1}(A) \to B(R,P,L)_s(A) \to \bigvee_{T \in \mathsf{Q}_s(A)}
\underline{w}(T) \smsh (R,P,L)_A(T) \tag{$*$} \]
Summing these over $s$ we obtain an exact couple and hence a spectral
sequence. The $E^1$ term of this spectral sequence is
\[ \begin{split}
    E^1_{s,t} &:= H_{s+t} (\bigvee_{T \in \mathsf{Q}_s(A)}
    \underline{w}(T) \smsh (R,P,L)_A(T)) \\
    &\isom \bigoplus_{T \in \mathsf{Q}_s(A)} H_{s+t} (\underline{w}(T)
    \smsh (R,P,L)_A(T)) \\
    &\isom \bigoplus_{T \in \mathsf{Q}_s(A)}
    \widetilde{H}_s(\underline{w}(T)) \otimes H_t((R,P,L)_A(T)) \\
    &\isom \bigoplus_{T \in \mathsf{Q}_s(A)}
    \widetilde{H}_s(\underline{w}(T)) \otimes (H_*R,H_*P,H_*L)_A(T)_t \\
    &\isom B(H_*R,H_*P,H_*L)_{s,t}(A). \\
\end{split} \]
where we have made plentiful use of the K\"{u}nneth formula. In
particular, we need the flatness assumptions to get
\[ H_*((R,P,L)_A(T)) \isom (H_*R,H_*P,H_*L)_A(T). \]
The final isomorphism is that of Lemma \ref{lem:explicit_bar}. Since
the filtration of each individual $B(R,P,L)(A)$ is finite, this spectral
sequence certainly converges to $H_*B(R,P,L)(A)$. It remains to be shown
that $d^1$ is given by the differential $\partial$ of the algebraic
bar construction.

The differential $d^1$ is the composite
\[ \begin{split} H_*\left(\bigvee_{T \in \mathsf{Q}_s(A)} \underline{w}(T)
\smsh (R,P,L)_A(T)\right) & \to H_{*-1}B(R,P,L)_{s-1}(A) \\
     \to & H_{*-1}\left(\bigvee_{U \in \mathsf{Q}_{s-1}(A)}
     \underline{w}(U) \smsh (R,P,L)_A(U)\right) \end{split} \]
of the boundary map in the long exact sequence associated to one of the
cofibre sequences $(*)$, with the projection map from another one. To
analyze this, fix for the moment a generalized $A$--labelled tree $T$
with $s$ vertices and consider the following map of cofibre sequences:
\small\[ \begin{diagram}  \dgARROWLENGTH=.84em
    \node{\partial w(T)_+ \smsh (R,P,L)_A(T)} \arrow{s} \arrow{e}
    \node{w(T)_+ \smsh (R,P,L)_A(T)} \arrow{s} \arrow{e}
    \node{\underline{w}(T) \smsh (R,P,L)_A(T)} \arrow{s} \\
    \node{B(R,P,L)_{s-1}(A)} \arrow{e}
    \node{B(R,P,L)_s(A)} \arrow{e}
    \node{\bigvee_{T \in \mathsf{Q}_s(A)} \underline{w}(T) \smsh
    (R,P,L)_A(T)}
\end{diagram} \] \normalsize
This induces a map of long exact sequences in homology, and in particular
we have a commutative diagram
\[ \begin{diagram}
    \node{H_*(\underline{w}(T) \smsh (R,P,L)_A(T))} \arrow{s} \arrow{e}
    \node{H_{*-1}(\partial w(T)_+ \smsh (R,P,L)_A(T))} \arrow{s} \\
    \node{H_*(\bigvee_{T \in \mathsf{Q}_s(A)} \underline{w}(T) \smsh
    (R,P,L)_A(T))} \arrow{e}
    \node{H_{*-1}B(R,P,L)_{s-1}(A).}
\end{diagram} \]
On the other hand, using the identity
\[ \partial w(T)_+ \isom \colim_{U < T} w(U)_+ \]
we also have a commutative diagram
\[ \begin{diagram}
    \node{\partial w(T)_+ \smsh (R,P,L)_A(T)} \arrow{s} \arrow{e}
    \node{\bigvee_{U \in \mathsf{Q}_{s-1}(A):\; U < T}
    \underline{w}(U) \smsh (R,P,L)_A(U)} \arrow{s} \\
    \node{B(R,P,L)_{s-1}(A)} \arrow{e}
    \node{\bigvee_{U \in \mathsf{Q}_{s-1}(A)} \underline{w}(U) \smsh
    (R,P,L)_A(U),}
   \end{diagram} \]
where the top horizontal map is constructed from the quotient maps
\[ \partial w(T)_+ \to \underline{w}(U), \]
for $U \in \mathsf{Tree}_{s-1}(A)$ such that $U < T$, together with the
operad composition maps
\[ (R,P,L)_A(T) \to (R,P,L)_A(U). \]
Taking the homology of this diagram, combining it with our other
diagram of homology groups, throwing in the K\"{u}nneth formula,
summing the top lines over all $T \in \mathsf{Q}_s(A)$ and using Lemma
\ref{lem:explicit_bar}, we get the big commutative diagram of Figure
\ref{fig:big_diagram} in which the top row is the differential $\partial$
on the algebraic bar construction $B(H_*R,H_*P,H_*L)(A)$ (under the
isomorphism of Lemma \ref{lem:explicit_bar}) and the bottom row is the
differential $d^1$ of our spectral sequence. The left and right sides
of the diagram are the isomorphisms described at the beginning of this
proof that identify $E^1$ with the algebraic bar construction.

\begin{sidewaysfigure}[p] \small \dgARROWLENGTH=.64em
\[ \begin{diagram}
 \node{B(H_*R,H_*P,H_*L)_{s,*}(A)} \arrow[2]{e,t}{\partial}
 \arrow{s,l}{\isom}
 \node[2]{B(H_*R,H_*P,H_*L)_{s-1,*}(A)} \arrow{s,r}{\isom} \\
 \node{\bigoplus_{T \in \mathsf{Q}_s(A)} \widetilde{H}_s(\underline{w}(T))
 \otimes (H_*...)_A(T)} \arrow{s,l}{\isom} \arrow{e}
 \node{\bigoplus_{T \in \mathsf{Q}_s(A)} \widetilde{H}_{s-1}(\partial w(T)_+)
 \otimes (H_*...)_A(T)} \arrow{s,l}{\isom} \arrow{e}
 \node{\bigoplus_{U \in \mathsf{Q}_{s-1}(A)}
 \widetilde{H}_{s-1}(\underline{w}(U)) \otimes (H_*...)_A(U)}
 \arrow{s,r}{\isom} \\
 \node{\bigoplus_{T \in \mathsf{Q}_s(A)} H_{s+*}(\underline{w}(T) \smsh
 (...)_A(T))} \arrow{s,l}{\isom} \arrow{e}
 \node{\bigoplus_{T \in \mathsf{Q}_s(A)} H_{s+*-1}(\partial w(T)_+ \smsh
 (...)_A(T))} \arrow{e} \arrow{s}
 \node{\bigoplus_{U \in \mathsf{Q}_{s-1}(A)} H_{s+*-1}(\underline{w}(U)
 \smsh (...)_A(U))} \arrow{s,r}{\isom} \\
 \node{H_{s+*}(\bigvee_{T \in \mathsf{Q}_s(A)} \underline{w}(T) \smsh
 (R,P,L)_A(T))} \arrow{s,l}{\isom} \arrow{e}
 \node{H_{s+*-1}B(R,P,L)_{s-1}(A)} \arrow{e}
 \node{H_{s+*-1}(\bigvee_{U \in \mathsf{Q}_{s-1}(A)} \underline{w}(U)
 \smsh (R,P,L)_A(U))} \arrow{s,r}{\isom} \\
 \node{E^1_{s,*}} \arrow[2]{e,b}{d^1}
 \node[2]{E^1_{s-1,*}}
\end{diagram} \]
\caption{Another big commutative diagram. This shows that the tree
differential on $B(H_*R,H_*P,H_*L)$ is the same as the $d^1$ differential
in the spectral sequence of Proposition \ref{prop:specseq}. Here $P$ is a
reduced operad in $\based$ or $\spectra$ with right module $R$ and left
module $L$. In the second row, $(H_*...)$ stands for $(H_*R,H_*P,H_*L)$
and in the middle row $(...)$ stands for $(R,P,L)$.}
\label{fig:big_diagram}
\end{sidewaysfigure}

The argument for the cobar construction is dual but only applies when
we are working in a category of spectra. The sequence of isomorphisms
that identifies the $E^1$ term then takes the form
\[ \begin{split}
    E^1_{-s,t} &:=
    H_{-s+t}(\prod_{T\in\mathsf{Q}_s(A)}\Map_{\spectra}(\underline{w}(T),(R,Q,L)_A(T)))\\
    &\isom
    \bigoplus_{T\in\mathsf{Q}_s(A)}H_{-s+t}\Map_{\spectra}(\underline{w}(T),(R,Q,L)_A(T)))
    \\
    &\isom
    \bigoplus_{T\in\mathsf{Q}_s(A)}\Hom(H_s(\underline{w}(T)),H_t(R,Q,L)_A(T))
    \\
    &\isom \bigoplus_{T\in\mathsf{Q}_s(A)}
    \Hom(H_s(\underline{w}(T)),(H_*R,H_*Q,H_*L)_A(T)_t) \\
    &\isom \Omega(H_*R,H_*Q,H_*L)_{s,t}(A). \\
\end{split} \]
In particular we use the fact that we are working with spectra and not
based spaces to get the isomorphism
\[ H_{-s+t}\Map(\underline{w}(T),X) \isom H_{-s+t}(\Sigma^{-s}X) \isom
H_t X \isom \Hom(H_s(\underline{w}(T)),H_t X) \]
that replaces an application of the K\"{u}nneth formula in the bar
construction case.
\end{proof}

\begin{remark}
Notice that the spectral sequence for the bar construction lies in the
right half-plane (and the first quadrant if the objects $R(A),P(A),L(A)$
only have non-negative homology). That for the cobar construction
lies in the left half-plane (and the second quadrant if the objects
$R(A),Q(A),L(A)$ only have non-negative homology).
\end{remark}

\subsection{The link to Koszul duality} \label{sec:koszul}

We now use our spectral sequence to look at the relationship between the
bar construction on an operad in based spaces or spectra and Koszul
duality. Koszul duality for operads initially appeared in
Ginzburg--Kapranov \cite{ginzburg/kapranov:1994}. Further references
include Getzler--Jones \cite{getzler/jones:1994} and Fresse
\cite{fresse:2004}.

The main result of this section is that if $P$ is a reduced operad in
based spaces or spectra such that $H_*P$ is a Koszul operad in graded
$k$--modules, then the spectral sequence for calculating $H_*B(P)$
collapses at the $E^2$--term and we conclude that $H_*B(P)$ is the
Koszul dual cooperad of $H_*P$. This result is a simple consequence of
the definitions of a Koszul operad and its Koszul dual cooperad. The
dual result holds for cooperads in spectra.

\begin{definition}[(Koszul operads)]
Let $P$ be a reduced operad in the category $\mathsf{Mod}_k$ of graded
$k$--modules. We say $P$ is \emph{Koszul} if the homology of the reduced
bar construction on $P$ is concentrated in the top tree degree. We
explain exactly what we mean by this. The reduced bar construction $B(P)$
is given by
\[ B(P)_{s,*}(A) \isom \bigoplus_{T \in \overline{\mathsf{Q}}_s(A)}
\widetilde{H}_s(\underline{w}(T)) \otimes P_A(T). \]
where $\overline{\mathsf{Q}}_s(A) = \mathsf{Q}_s(A) \cap \mathsf{T}_s(A)$ is
the set of $A$--labelled trees (in the sense of Section~\ref{sec:trees},
that is, \emph{not} generalized trees) with exactly $s$ vertices. If $|A|
= 1$, this is concentrated in the column $s = 0$. If $|A| > 1$, it is
concentrated in $1 \leq s \leq |A|-1$. We say that $P$ is \emph{Koszul}
if, for all $A$,
\[ H_{s,*}(B(P)(A),\partial) = 0 \text{ for } s \neq |A|-1 \]
where $\partial$ denotes the tree differential on $B(P)$.
\end{definition}

\begin{definition}[(Koszul duals)]
Let $P$ be a Koszul operad in graded $k$--modules. The \emph{Koszul
dual} of $P$ is the symmetric sequence $K(P)$ given by the homology of
the reduced bar construction on $P$. We grade $K(P)$ according to the
total degree (that is, internal degree plus tree degree) of $B(P)$:
\[ K(P)_r(A) = H_{|A|-1,r+1-|A|}(B(P)(A),\partial). \]
Notice that $K(P)(A)$ is the kernel of the differential
\[ B(P)_{|A|-1,*}(A) \to B(P)_{|A|-2,*}(A), \]
so there is a natural inclusion
\[ K(P) \to B(P). \]
\end{definition}

\begin{prop} \label{prop:koszul_cooperad}
Let $P$ be a Koszul operad in graded $k$--modules such that each $K(P)(A)$
is a flat $k$--module. Then the Koszul dual $K(P)$ has a natural cooperad
structure.
\end{prop}
\begin{proof}
We already know from Definition \ref{def:alg_cooperad(bar)} that the bar
construction $B(P)$ has a cooperad structure. We get the structure for
$K(P)$ by taking homology. So cocomposition maps for $K(P)$ are given by
\[ H(B(P)(A \cup_a B)) \to H(B(P)(A) \otimes B(P)(B)) \isom H(B(P)(A))
\otimes H(B(P)(B)) \]
where we use the flatness assumption to get the isomorphism.
\end{proof}

We dually define the Koszul property and Koszul dual for cooperads of
graded $k$--modules.

\begin{definition}[(Koszul cooperads and Koszul duals)]
Let $Q$ be a reduced cooperad of graded $k$--modules. Then $Q$ is
\emph{Koszul} if the homology of the reduced cobar construction is
concentrated in the lowest\footnote{Recall that the tree grading for the
cobar construction is concentrated in negative degrees. `Lowest' here
means `most negative'.} tree degree. In this case, the \emph{Koszul dual}
of $Q$ is the symmetric sequence $K(Q)$ of graded $k$--modules with
\[ K(Q)_r(A) := H_{1-|A|,r+|A|-1}(\Omega(Q)(A),\partial^*), \]
where $\partial^*$ is the tree differential on $\Omega(Q)$. Since $K(Q)$
is the bottom homology group of $\Omega(Q)$ there is a natural surjection
\[ \Omega(Q) \to K(Q). \]
\end{definition}

\begin{proposition}
Let $Q$ be a Koszul cooperad of graded $k$--modules. Then the Koszul
dual $K(Q)$ has a natural operad structure.
\end{proposition}
\begin{proof}
The composition maps for $K(Q)$ are given by
\[ H(\Omega(Q)(A)) \otimes H(\Omega(Q)(B)) \to H(\Omega(Q)(A) \otimes
\Omega(Q)(B)) \to H(\Omega(Q)(A \cup_a B). \]
Notice that we don't need a flatness assumption here.
\end{proof}

Fresse \cite{fresse:2004} gives various fundamental results for Koszul
duality of operads and cooperads, in particular, the following.

\begin{lem}[Fresse,\cite{fresse:2004}, Lemma 5.2.10]
Let $P$ be a Koszul operad of graded $k$--modules such that the
$k$--modules $P(A)$ and $K(P)(A)$ are flat. Then $K(P)$ is a Koszul
cooperad and
\[ K(K(P)) \isom P \]
as operads. Dually, let $Q$ be a Koszul cooperad of graded $k$--modules
such that the modules $Q(A)$ and $K(Q)(A)$ are flat. If $Q$ is Koszul
then its Koszul dual operad $K(Q)$ is also Koszul and
\[ K(K(Q)) \isom Q \]
as cooperads. \qed
\end{lem}

We now give the main result of this section.

\begin{prop} \label{prop:koszul}
Let $P$ be a reduced operad in $\based$ or $\spectra$ such that each
object $P(A)$ is cofibrant and all homology groups $H_*P(A)$ and
$H_*B(P)(A)$ are flat $k$--modules. If $H_*P$ is a Koszul operad then
\[ H_*B(P) \isom K(H_*P) \]
as cooperads.

Dually, let $Q$ be a reduced cooperad in $\spectra$ such that each object
$Q(A)$ is fibrant and the homology groups $H_*Q(A)$ are flat
$k$--modules. If $H_*Q$ is a Koszul cooperad then
\[ H_*\Omega(Q) \isom K(H_*Q) \]
as operads.
\end{prop}
\begin{proof}
The cofibrancy and flatness conditions ensure that the spectral sequence
of Proposition \ref{prop:specseq} exists for each finite set $A$ and
that $H_*B(P)$ is a cooperad in $\mathsf{Mod}_k$. We have already seen
that the spectral sequence has the form
\[ (E^1_{*,*},d^1) = (B(H_*P)_{*,*}(A),\partial) \implies H_*B(P). \]
Because $H_*P$ is Koszul, the homology of the bar construction is
concentrated in the $s = |A|-1$ column. Therefore, the $E^2$--term is
concentrated in this column and so the spectral sequence collapses. We
then see that
\[ H_r B(P)(A) \isom E^2_{|A|-1,r-|A|+1} \isom
H_{|A|-1,r-|A|+1}(B(H_*P)(A),\partial) \isom K(H_*P)_r(A) \]
and so
\[ H_*B(P) \isom K(H_*P) \]
as claimed. It follows that the modules $K(H_*P)(A)$ are flat so,
by Proposition \ref{prop:koszul_cooperad}, $K(H_*P)$ has a cooperad
structure. It remains to show that this cooperad structure agrees with
that on $H_*B(P)$.

The first thing to notice is that the above identification of $H_*B(P)(A)$
with the submodule $K(H_*P)(A)$ on $B(H_*P)(A)$ is realized by an edge
homomorphism of our spectral sequence. This edge homomorphism comes from
applying homology to the quotient map
\[ B(P)(A) \to \bigvee_{T \in \mathsf{Q}_s(A)} \underline{w}(T) \smsh
P_A(T) \]
where $s = |A|-1$.
The key property of these maps is that they fit into commutative diagrams
\[ \begin{diagram} \dgARROWLENGTH=2.4em
    \node{B(P)(A \cup_a B)} \arrow{s} \arrow{e}
    \node{\bigvee_{V \in \mathsf{Q}_{s+s'}(A \cup_a B)} \underline{w}(V)
    \smsh P_{A \cup_a B}(V)} \arrow{s} \\
    \node{B(P)(A) \smsh B(P)(B)} \arrow{e}
    \node{\bigvee_{T \in \mathsf{Q}_s(A)} \bigvee_{U \in
    \mathsf{Q}_{s'}(B)} \underline{w}(T) \smsh \underline{w}(U) \smsh
    P_A(T) \smsh P_B(U)}
\end{diagram} \]
where the map on the right-hand side is built from the familiar maps
\[ \underline{w}(T \cup_a U) \to \underline{w}(T) \smsh \underline{w}(U)
\]
and the isomorphisms
\[ P_{A \cup_a B}(T \cup_a U) \to P_A(T) \otimes P_B(U) \]
with terms for trees $V$ not of type $(A,B)$ mapping to the basepoint.

Taking homology of this diagram, the right-hand side map
gives the cooperad structure on $B(H_*P)$ as described in Lemma
\ref{lem:explicit_cooperad}. This shows that the edge homomorphisms
of the spectral sequence identify the cooperad structure on $H_*B(P)$
with the restriction of that on $B(H_*P)$. Since the cooperad structure
on $K(H_*P)$ is also the restriction of that on $B(H_*P)$, it follows that
\[ H_*B(P) \isom K(H_*P) \]
is an isomorphism of cooperads. The dual result is proved similarly.
\end{proof}

\begin{remark}
Proposition \ref{prop:koszul} extends a result of Vallette
\cite{vallette:2004} for discrete operads. Recall from Remark
\ref{rem:vallette} that he constructs the `order complex' for an operad
$P$ in $\mathsf{Set}$. His main result then is that $H_*P$ is Koszul
if and only if the homology of the order complex is concentrated in top
degree. This follows immediately from our spectral sequence argument by
identifying the order complex with the bar construction.
\end{remark}

\begin{example} \label{ex:derivatives_specseq}
We finally return to the Goodwillie derivatives of the identity
functor. Recall that
\[ \partial_* I \isom \Omega(\underline{S}) \]
where $\underline{S}$ is the cooperad of spectra with $\underline{S}(A)
= S$ for all $A$. The homology of this cooperad is given by
\[ H_*(\underline{S})(A) = \begin{cases} k & \text{if $* = 0$}; \\ 0 &
\text{otherwise}; \end{cases} \]
for all finite sets $A$. This is the cooperad of commutative coalgebras in
the category of graded $k$--modules. Fresse shows in
\cite[Section~6]{fresse:2004} (by updating a result of Ginzburg and
Kapranov \cite{ginzburg/kapranov:1994}) that this cooperad is Koszul
(for $k = \mathbb{Q},\mathbb{F}_p,\mathbb{Z}$) with Koszul dual given by
a suspension of the Lie operad. Proposition \ref{prop:koszul} therefore
applies and we recover the homology of the derivatives of the identity:
\[ H_*(\partial_n I) = \begin{cases} \mathsf{Lie}(n) \otimes sgn_n &
\text{if $* = 1-n$}; \\ 0 & \text{otherwise}. \end{cases} \]
Moreover, we now know that the induced operad structure on this homology
of the derivatives is equal to the operad structure on the Koszul dual
of the commutative cooperad, that is, the desuspended Lie structure. This
completes the main goal set out in the introduction to this paper.
\end{example}

\subsection{Homology of modules over the derivatives of the identity}
\label{sec:modules_specseq}

In this final section, we use our spectral sequence to investigate the
homology of the left $\partial_*I$--module $M_X$ associated to a based
space $X$ as described in Remark \ref{rem:modules}(1). Recall that this
module is given by a cobar construction:
\[ M_X := \Omega(I,\underline{S},\underline{\Sigma^{\infty}X}). \]
We can describe explicitly the spectral sequence for calculating
$H_*M_X(2)$. The cobar construction is one-sided and we only have to
consider trees for which the root has a single incoming edge. There are
two $2$--labelled trees of this type with zero and one vertices
respectively and a morphism between them. The $E^1$ term in the spectral
sequence therefore only has nonzero entries in the columns $s = 0$ and
$s = -1$. These entries are respectively $H_*X$ and $H_*(X \smsh X) \isom
H_*X \otimes H_*X$ with the differential given by the reduced diagonal $X
\to X \smsh X$. The spectral sequence therefore takes the following form.

\begin{center}
\input{M_X2.pstex_t}
\end{center}
This reduces to the long exact sequence of homology determined by the
cofibre sequence
\[ X \to X \smsh X \to \hocofib(X \to X \smsh X) \]
which is consistent with the calculation of $M_X(2)$ made in Remark
\ref{rem:modules}.

Things become more interesting (and much more complicated) for $M_X(n)$
when $n > 2$. For $n = 3$ there are eight trees of interest:
\begin{center}
\input{T_root3.pstex_t}
\end{center}
and the $E^1$ term of the spectral sequence takes the form
\begin{center}
\input{M_X3.pstex_t}
\end{center}
The differential $d^1$ is built from the reduced diagonal (between pairs
of terms corresponding to bud collapse) and isomorphisms (between pairs
of terms corresponding to collapse of an internal edge).

We will close the paper by looking at $X = S^r$, the $r$--sphere (for
$r \geq 2$). In this situation the reduced diagonal is zero on homology
and there can be no higher differentials or extensions in the spectral
sequence. This allows us to calculate $H_*M_{S^r}$ with $\mathbb{Z}$
coefficients in its entirety.

\begin{prop}
Let $S^r$ denote the $r$--sphere for $r \geq 2$. Then we have
\[ H_*(M_{S^r}) \isom H_*(\partial_*I) \circ H_*(\underline{S^r}) \]
where $H_*(\underline{S^r})$ is the symmetric sequence with
\[ H_*(\underline{S^r})(n) = \begin{cases} \mathbb{Z} & \text{if $* =
r$}; \\ 0 & \text{otherwise}. \end{cases} \]
The left action of $H_*(\partial_*I)$ on $H_*(M_{S^r})$ is given by the
operad structure on $H_*(\partial_*I)$.
\end{prop}
\begin{proof}
The $E^1$ term of the spectral sequence for the homology of $M_X$ is in
this case the algebraic cobar construction
\[ \Omega(I,H_*(\underline{S}),H_*(\underline{S^r})). \]
The coaction of $H_*(\underline{S})$ on $H_*(\underline{S^r})$ is trivial
in the sense that the only nonzero cocomposition maps are
\[ H_*(\underline{S^r})(n) \to H_*(\underline{S})(1) \otimes
H_*(\underline{S^r})(n). \]
This is equivalent to saying that
\[ H_*(\underline{S^r}) \isom I \circ H_*(\underline{S^r}) \]
as left $H_*(\underline{S})$--comodules, where the coaction of
$H_*(\underline{S})$ on the right-hand side is via the coaugmentation
action on $I$. It follows that the $E^1$ term of our spectral sequence
can be written
\[ \Omega(I,H_*(\underline{S}),I \circ H_*(\underline{S^r})) \isom
\Omega(I,H_*(\underline{S}),I) \circ H_*(\underline{S^r}) \]
where the differential on the right-hand side comes solely from the cobar
construction and not from $H_*(\underline{S^r})$. This isomorphism can be
seen by working through the definition of the algebraic bar construction
in this case.

It now follows that the $E^2$ term of our spectral sequence is given by
\[ H_*(\partial_*I) \circ H_*(\underline{S^r}). \]
In the $E^2$ term for calculating $H_*M_{S^r}(n)$, we only have nonzero
entries in bidegrees $(-k,r(k+1))$ for integers $k \geq 0$. Since $r \geq
2$ there can be no further differentials or extensions and so we see that
\[ H_*(M_{S^r}) \isom H_*(\partial_*I) \circ H_*(\underline{S^r}). \]
The proof of Proposition \ref{prop:koszul} extends to show that the left
action of $H_*(\partial_*I)$ is as claimed.
\end{proof}

\begin{remark}
The functor $P \circ -$ from symmetric sequences to left $P$--modules is
left adjoint to the forgetful functor and so can rightfully be called the
\emph{free left $P$--module} functor. Hence the homology of $M_{S^r}$ is
the free left $P$--module on $H_*(\underline{S^r})$.

Explicitly, there is a generator $x_A$ in $H_r(M_{S^r})(A)$ for each
finite set $A$. The entire homology group $H_*(M_{S^r})(A)$ then has a
basis given by all possible iterated brackets of the form
\[ [\dots [[x_{A_1},x_{A_2}],x_{A_3}] \dots,x_{A_k}] \]
where $A_1,\dots,A_k$ is a partition of $A$ into nonempty finite subsets,
and $[-,-]$ is a Lie bracket of degree $-1$. This Lie bracket also
represents the action of $H_*(\partial_*I)$ on $H_*(M_{S^r})$.
\end{remark}


\end{document}